\documentclass[reqno]{amsart}
\usepackage{geometry}
\usepackage{amsmath, amssymb, amsthm, amsfonts}
\usepackage{physics}
\usepackage{mathtools}
\usepackage{mathrsfs}
\usepackage{graphicx}
\usepackage{setspace}
\usepackage{latexsym}
\usepackage{mleftright}
\usepackage{cite}
\usepackage[hidelinks]{hyperref}
\usepackage{enumerate}

\allowdisplaybreaks[4]

\makeatletter
\providecommand{\@secnumpunct}{.}
\makeatother

\numberwithin{equation}{section}

\newtheorem{Theorem}{Theorem}[section]
\newtheorem{Definition}{Definition}[section]
\newtheorem{Lemma}{Lemma}[section]
\newtheorem{Remark}{Remark}[section]

\newtheorem{Proposition}{Proposition}[section]

\newcommand{\Om}{\Omega}
\newcommand{\ud}{\mathrm{d}}
\newcommand{\bm}{\boldsymbol}

\makeatletter
\def\@setaddresses{}

\renewcommand{\@maketitle}{%
  \normalfont\normalsize
  \let\@makefnmark\relax  \let\@thefnmark\relax
  \ifx\@empty\@date\else \@footnotetext{\@setdate}\fi
  \ifx\@empty\@subjclass\else \@footnotetext{\@setsubjclass}\fi
  \ifx\@empty\@keywords\else \@footnotetext{\@setkeywords}\fi
  \ifx\@empty\thankses\else \@footnotetext{%
    \def\par{\let\par\@par}\@setthanks}\fi
  \@mkboth{\@nx\shortauthors}{\@nx\shorttitle}%
  \global\topskip42\p@\relax
  \@settitle
  \ifx\@empty\authors \else \@setauthors \fi
  \ifx\@empty\@dedicatory
  \else
    \baselineskip18\p@
    \vtop{\centering{\footnotesize\itshape\@dedicatory\@@par}%
      \global\topskip30\p@\relax}\relax
  \fi
  \vskip 1.5ex
  {\small
    \begin{center}
    \textsc{Yinghua Li$^*$}\\
    School of Mathematical Sciences, South China Normal University, Guangzhou, 510631, China\\
    \textit{E-mail:} \texttt{yinghua@scnu.edu.cn}\\[1.2em]
    \textsc{Wenlin Ye}\\
    School of Mathematical Sciences, South China Normal University, Guangzhou, 510631, China\\
    \textit{E-mail:} \texttt{15807979775@163.com}
    \end{center}
  }
  \vskip 2ex
  \@setabstract
  \normalsize
}
\makeatother

\begin{document}
\title[]{Well-posedness of the nonhomogeneous incompressible Navier-Stokes/Allen-Cahn system}

\thanks{$^*$Corresponding author. Email: yinghua@scnu.edu.cn}
\thanks{{\it Keywords}: Nonhomogeneous Navier-Stokes equations; Allen-Cahn equations; weak solutions; strong solutions}
\thanks{{\it AMS Subject Classification}: 35Q35, 76D03, 76T05}


\begin{abstract}
In this paper, we investigate a diffuse interface model for two-phase flow of viscous, incompressible fluids with different densities. The model is described by a coupled nonhomogeneous Navier-Stokes/Allen-Cahn system in a bounded domain $\Omega\subset\mathbb R^d~(d=2, 3)$. The mobility may depend on the phase variable but is assumed to be non-degenerate. We first establish the existence of global weak solutions in both two and three dimensions. When the viscosity is constant, we then obtain local strong solutions in 3D and global strong solutions in 2D. Furthermore, under smallness conditions on the initial data, the 3D local strong solution is extended globally in time, and the solution converges exponentially to the equilibrium. Finally, we prove a weak-strong uniqueness theorem, which ensures that the weak solutions coincide with the strong solutions whenever the latter exist. These results constitute a comprehensive well-posedness theory for the nonhomogeneous incompressible Navier-Stokes/Allen-Cahn system with the density fully coupled into the Allen-Cahn dynamics.

\end{abstract}

\maketitle
\vspace{-5mm}

\section{Introduction}
In diffuse interface models for two-phase flows, the phase separation is commonly described by either the Cahn-Hilliard equation or the Allen-Cahn equation. The former captures conserved dynamics, where the total mass of each phase is preserved and the evolution is driven by diffusive fluxes of the chemical potential. The latter captures non-conserved dynamics, where the phase evolution is driven by local relaxation of the free energy -- a scenario relevant to phenomena like cavitation, non-conserved domain coarsening, and certain solid-liquid phase transitions.

This paper is devoted to the corresponding nonhomogeneous incompressible Navier-Stokes/ Allen-Cahn (NSAC) system, a diffuse interface model for two-phase flows of viscous fluids with different densities. The system reads
\begin{equation} \label{1.1}    
\begin{cases}
\partial_t\rho+{\bm u}\cdot\nabla\rho=0,
\\
\rho\partial_t{\bm u}+\rho({\bm u}\cdot\nabla){\bm u}-{\rm div}(\eta(\chi){\mathbb D}{\bm u})
+\nabla p=-{\rm div}(\nabla\chi\otimes\nabla\chi),
\\
{\rm div}{\bm u}=0,
\\
\rho\partial_t\chi+\rho{\bm u}\cdot\nabla\chi=-m(\chi)\mu,
\\
\rho\mu=-\Delta\chi+\rho F^{\prime}(\chi)
\end{cases}
\end{equation}
in $\Omega\times(0,T)$, where $\Omega$ is a smooth bounded domain in $\mathbb R^d~(d = 2, 3)$.
Here $\rho$, $u$, $p$, $\chi$, $\mu$ denote the density function, the velocity field, the pressure,
the phase variable, and the chemical potential, respectively.
We assume that the viscosity $\eta(s)$ and the mobility $m(s)$ are smooth on $\mathbb R$ satisfying
$0<\eta_\ast\leq \eta(s) \leq \eta^\ast$ and $0<m_\ast\leq m(s) \leq m^\ast$ for all $s\in\mathbb R$
and positive constants $\eta_\ast, \eta^\ast, m_*, m^*$.
${\mathbb D}{\bm u}=\frac12(\nabla{\bm u}+\nabla{\bm u}^{T})$.
We supplement system (\ref{1.1}) with the following initial and boundary value conditions
\begin{align}    
({\bm u},\, \partial_{\mathbf{n}}\chi)&=(0,\,0) \qquad {\rm on}\; \partial\Omega\times(0,T),
\\
\label{1.2}
\left(\rho,\, {\bm u}, \, \chi\right)(\cdot,0)&=\left(\rho_0, {\bm u}_0, \chi_0\right) \qquad{\rm in}\;\Omega.
\end{align}
We take $F$ to be the Landau double-well potential, namely,
\begin{align}  \label{1.4}   
F(s)=\frac{1}{4}(s^2-1)^2, \quad s\in\mathbb R.
\end{align}
For sufficiently regular solutions of (\ref{1.1})--(\ref{1.2}), 
the following energy balance holds
\begin{equation*} 
\frac{\rm d}{{\rm d}t}\int_{\Omega}\left(\frac{1}{2}\rho|{\bm u}|^2+\frac{1}{2}|\nabla\chi|^2+\rho F(\chi)\right){\rm d}x
+{\int_{\Omega}\left({\eta}(\chi )|\mathbb{D} \boldsymbol{u}|^2+m(\chi )|\mu|^2\right)} \mathrm{d}x =0.
\end{equation*}

\vskip1mm
The corresponding Navier-Stokes/Allen-Cahn system was first proposed by Blesgen \cite{B-99} in the compressible setting. In that framework, Feireisl et al. \cite{F-10} established weak solutions and Kotschote \cite{K-12} obtained local strong solutions; further results include time-periodic solutions, spherically symmetric solutions and decay estimates \cite{S-Z-W-20, S-Z-W-24, Z-22}. In the incompressible setting, Li and Huang \cite{L-H} proved local strong solutions for the nonhomogeneous system, and blow-up criteria were derived \cite{L-D-H, Z-16, F-L}. A simplified version \cite{F-P-P}, where the density dependence in the Allen-Cahn equation is removed, was recently studied by Li, Xie and Yan \cite{L-X-Y}, who obtained global existence and exponential decay.
However, global strong solvability for the fully coupled nonhomogeneous incompressible NSAC system remains largely open when the viscosity is nonconstant, even in two spatial dimensions.

\vskip1mm
{ Another important direction concerns mass-conserving Allen-Cahn approximations for incompressible binary fluids. Giorgini, Grasselli, and Wu \cite{G-G-W} studied a mass-conserving NSAC system with unmatched densities, nonconstant viscosities, and a Flory--Huggins potential. They established global weak solutions in 2D and 3D, and global strong solutions in 2D, together with further uniqueness and separation results under additional assumptions. For the corresponding two-dimensional model, Sun and Tan \cite{S-T} subsequently studied the long-time dynamics and established a global attractor in the strong-solution phase space. }

\vskip1mm
{ More recently, Grasselli, Hurm, and Poiatti \cite{G-H-P} investigated the conserved NSAC system with constant density, nonconstant viscosity, and a singular potential. They established global weak solutions, local strong solutions in 3D, a conditional weak-strong uniqueness principle based on relative energy, and unconditional uniqueness of weak solutions in 2D. They also proved asymptotic strict separation from the
pure phases and, under an analyticity assumption on the potential, convergence of each weak trajectory to a single equilibrium. Their convergence analysis allows merely continuous positive viscosity coefficients; under stronger viscosity regularity, they additionally obtained eventual regularization and convergence in stronger norms. }

\vskip1mm
{ The systems studied in these works differ structurally
from \eqref{1.1}.
In the conserved Allen--Cahn dynamics, a spatially
nonlocal correction enforces conservation of the mean
phase.
In the present system, no such correction is imposed:
the density is an independent unknown satisfying a
transport equation, and it also appears in the
phase-field evolution and chemical-potential relation.
Accordingly, the results for conserved NSAC systems
provide an important context for the present analysis,
while the density coupling in \eqref{1.1} requires
separate estimates. }

The Navier-Stokes/Cahn-Hilliard (NSCH) system, widely known as Model H following Hohenberg and Halperin \cite{H-H-77}, has been extensively studied since its derivation by Gurtin, Polignone and Vi\~{n}als \cite{G-P-V} from rational continuum mechanics. For mixtures with different densities, various formulations have been proposed \cite{A-F, A-G-G, B-01, L-T-98}. In the nonhomogeneous incompressible framework, Giorgini and Temam \cite{G-T} recently established a comprehensive well-posedness theory: they proved the existence of global weak solutions in both 2D and 3D, global strong solutions in 2D, and local strong solutions in 3D.
The uniqueness of strong solutions was later settled by Jiang, Wu and Xu \cite{J-W-X}.
More recently, Fang, Nei and Guo \cite{F-N-G}
proved the global existence of strong solutions in 3D
and obtained a decay-in-time estimate under the
smallness condition
\begin{equation}\label{eq:JDE-smallness}
\|\nabla u_0\|_{L^2} + \|\nabla\mu_0\|_{L^2} + \|\rho_0\|_{L^\infty} \leq \varepsilon_0,
\end{equation}
where the analysis relies crucially on the fourth-order Cahn-Hilliard structure that yields $\nabla\mu \in L^2(0,T;L^2(\Omega))$, which further implies an $L^4(0,T;H^2(\Omega))$-norm bound on $\chi$ and provides sufficient regularity to control the capillary stress term $\mathrm{div}(\nabla\chi\otimes\nabla\chi)$.

\vskip1mm
However, the above approach for the NSCH system does not directly apply to the NSAC system.
The essential obstruction lies in the lower regularity of the chemical potential: for the NSCH system, the energy inequality provides an $L^2(0,T;L^2(\Omega))$-norm bound on $\nabla\mu$, whereas for NSAC, one only obtains an $L^2(0,T;L^2(\Omega))$-norm bound on $\mu$ itself. This loss of gradient information prevents a direct transfer of the estimates in \cite{G-T} and makes the control of the capillary stress term substantially more delicate.
The aim of this paper is to address this gap. Specifically, we establish the following results:
\begin{enumerate}[(i)]
\item {\bf(Global weak solutions)} Under the assumptions that
    $(\rho_0, \boldsymbol{u}_0, \chi_0)\in L^\infty(\Omega)\times\mathbf{H}_{\sigma}\times H^1(\Omega)$
    is given and satisfies $0 <\rho_{*} \leq \rho_{0} \leq \rho^{*}$, $|\chi_0|\le1$ a.e. in $\Omega$,
   the system $(\ref{1.1})$--$(\ref{1.2})$ has a global weak solution for $d=2,3$.
   The same existence result holds for the physically relevant Flory-Huggins logarithmic potential; see Appendix \ref{B-F-H} for the necessary modifications.
   \vskip1mm
\item {\bf(Strong solutions)} Under the assumptions that
   $(\rho_0, \boldsymbol{u}_0, \chi_0)\in L^\infty(\Omega)\times\mathbf{V}_{\sigma}\times H^2(\Omega)$
   is given and satisfies $0 <\rho_{*} \leq \rho_{0} \leq \rho^{*}$, $|\chi_0|\le1$ a.e. in $\Omega$,
   the system $(\ref{1.1})$--$(\ref{1.2})$ has a strong solution.
   When $d=3$, the strong solution is local in time; when $d=2$ and the viscosity coefficient $\eta \equiv 1$,
   the strong solution is global for any $T>0$.
   \vskip1mm
\item  {\bf(Global strong solutions with exponential decay for small data)} When $d=3$ and $\eta \equiv 1$, if
    $\|\boldsymbol{u}_0\|_{\mathbf{V}_{\sigma}}+\|\chi_0^2-1\|_{L^2(\Omega)}+\|\nabla\chi_0\|_{H^1(\Omega)}$ is small enough,
    the local strong solution obtained in (ii) can be extended to {the} global one.
    Moreover, the solution converges exponentially to the equilibrium.
    \vskip1mm
\item {\bf(Weak-strong uniqueness)} When $d=2,3$, any weak solution obtained in (i) coincides with a strong solution originating from the
same initial data, provided that the latter satisfies the additional regularity assumptions of Theorem~\ref{Th7.1}. This ensures that the weak solutions are admissible and that the solution is unique within the class of weak solutions when strong solutions exist.
\end{enumerate}

\begin{Remark}
For the Flory-Huggins logarithmic potential
\begin{align} \label{1.5}
\tilde F(s)=\frac{\theta}{2}\Big[(1+s)\ln(1+s)+(1-s)\ln(1-s)\Big]-\frac{\theta_0}{2}s^2, \quad s\in[-1,1],
\end{align}
where $0<\theta<\theta_0$,
the existence of global weak solutions for the NSAC system is established by an adaptation of the Landau case; the details are provided in Appendix \ref{B-F-H}.
\end{Remark}

\vskip1mm
We point out that the results above constitute a comprehensive well-posedness theory for the nonhomogeneous incompressible NSAC system. In particular, the global strong solution in 2D and the global extension of the 3D strong solution with exponential decay rates (results (ii)-(iii)) are the first global existence results for the fully coupled NSAC system with a density-dependent Allen-Cahn equation. The exponential decay quantifies the relaxation rate of the system toward equilibrium under the smallness assumption, revealing the dissipative structure of the nonhomogeneous NSAC dynamics. Moreover, the weak-strong uniqueness (result (iv)) ensures that the weak solutions constructed in (i) are physically relevant, as they must coincide with the strong solutions originating from
the same initial data and satisfying the regularity
assumptions of Theorem~\ref{Th7.1}. Taken together, these results provide a rigorous and self-contained mathematical theory of the NSAC system within the considered solution classes, establishing a solid foundation for further investigations into qualitative and asymptotic properties of the dynamics.

\vskip1mm

The remainder of this article is organized as follows.
In Section 2, we collect the notation and preliminary tools.
In Section 3, we establish the existence of global weak solutions (Theorem \ref{Th4.1}) via a Galerkin approximation.
In Section 4, we derive higher-order a priori estimates and establish the existence of strong solutions (Theorem \ref{Th6.1}), which are local in 3D and global in 2D.
In Section 5, we extend the 3D strong solution globally in time under a smallness condition and derive exponential decay rates (Theorem \ref{Th8.1}). In Section 6, we prove the weak-strong uniqueness theorem (Theorem \ref{Th7.1}).
Appendix \ref{A-A} details the existence of approximate solutions via a fixed-point argument, and Appendix \ref{B-F-H} outlines how to extend Theorem \ref{Th4.1} to the Flory-Huggins logarithmic potential.

\section{Preliminaries}
Let $\Omega$ be a bounded domain with smooth boundary in $\mathbb R^d~(d = 2, 3)$.
$L^p(\Omega)$ and $W^{k,p}(\Omega)$ with $k$ being an integer and $p>1$ represent the usual Lebesgue spaces and Sobolev spaces,
respectively. As usual, $W^{k,2}(\Omega)$ is denoted by $H^{k}(\Omega)$.
We denote by $H_{0}^{1}(\Omega)$  the closure of $\mathcal{C}_{0}^{\infty}\left( \Omega \right)$ in $H^{1}(\Omega)$
and by $H^{-1}(\Omega)$ its dual space.
The inner product and norm in $L^{2}(\Omega)$ are denoted by $\left( \cdot ,\cdot \right)$ and  $\|\cdot\|_{L^{2}(\Omega)}$,  respectively.
By the generalized {Poincar\'{e}'s} inequality (\cite{T-97}, Chapter II, Section 1.4), we recall that
 \begin{equation}  \label{3.1}  
u\rightarrow \Big( \| \nabla u \| _{L^2\left( \Omega \right)}^{2}+\Big| \frac{1}{\left| \Omega \right|}\int_{\Omega}{u \ \ud x} \Big|^2 \Big) ^{\frac{1}{2}} \quad \text{and}\quad  u\rightarrow \Big( \| \nabla u \| _{L^2\left( \Omega \right)}^{2}+\Big| \int_{\Omega}{\xi u \ \ud x} \Big|^2 \Big) ^{\frac{1}{2}},
\end{equation}
%
 where $\xi \in L^ \infty \left( \Omega \right)$ is such that $0<\xi _{\ast}\le \xi \left( x \right) \le \xi ^{\ast}$ for almost every $x\in \Omega$,  are norms on $H^1\left( \Omega \right)$ equivalent to $\|u\|_{H^{1}(\Omega)}$. In particular, there exists a positive constant\  $C=C(\Om ,\xi_\ast, \xi^\ast)$\  such that
\begin{equation} \label{3.2}  
\|u\|_{H^{1}(\Omega)}\le C \Big( \| \nabla u \| _{L^2\left( \Omega \right)}^{2}+\Big| \int_{\Omega}{\xi u\  \ud x} \Big|^2 \Big) ^{\frac{1}{2}},\quad u\in H^{1}(\Omega).
\end{equation}

Now, we introduce the Hilbert spaces of vector-valued functions. We denote
\begin{align*}
\mathbf{H}_{\sigma}=\overline{\left\{\boldsymbol{u}\in C_0^{\infty}(\Omega): {\rm div}\boldsymbol{u}=0\;\;\text{in}\;\Omega\right\}}^{L^2}
\end{align*}
with scalar product and the norm of $L^2(\Omega)$
and
\begin{align*}
\mathbf{V}_{\sigma}=\overline{\left\{\boldsymbol{u}\in C_0^{\infty}(\Omega): {\rm div}\boldsymbol{u}=0\;\;\text{in}\;\Omega\right\}}^{H^1_0}
\end{align*}
endowed with the product
$\left( \boldsymbol{u},\boldsymbol{v} \right) _{\mathbf{V}_{\sigma}}=\left( \nabla \boldsymbol{u},\nabla \boldsymbol{v} \right) $
and norm $\left\| \boldsymbol{u} \right\| _{\mathbf{V}_{\sigma}}=\left\| \nabla \boldsymbol{u} \right\| _{L^2\left( \Omega \right)}$.
We consider the Hilbert space $\mathbf{W}_{\sigma}=\mathbf{H}^{{2}}\left( {\Omega } \right) \cap \mathbf{V}_{\sigma}$ with inner product and norm
$\left( \boldsymbol{u},\boldsymbol{v} \right) _{\mathbf{W}_{\sigma}}=\left( \mathbf{A}\boldsymbol{u},\mathbf{A}\boldsymbol{v} \right) $
and $\left\| \boldsymbol{u} \right\| _{\mathbf{W}_{\sigma}}=\left\| \mathbf{A}\boldsymbol{u} \right\| $,  where $\mathbf{A}$ is the Stokes operator. We recall that there exists $C>0$ such that
\begin{equation}  \label{3.9}   
\left\| \boldsymbol{u} \right\| _{H^2\left( \Omega \right)}\leq C\left\| \boldsymbol{u} \right\| _{\mathbf{W}_{\sigma}},\quad
\boldsymbol{u}\in \mathbf{W}_{\sigma}.
\end{equation}

Let $X$ be a (real) Banach or Hilbert space with norm $\|\cdot\|_X$. We denote by $X^{\prime}$ the dual space of X
and by $\left< \cdot ,\cdot \right>$,  the duality product between $X$ and $X'$.
Given $1\leq p\leq \infty$ and an interval $I\subseteq \left[ 0,\infty \right) $,  the set $L^p\left( I;X \right)$ consists of all Bochner measurable $p$-integrable functions defined on $I$ with values in $X$. We define the space $W^{1,p}\left( 0,T;X \right) $ to consist of all
functions $f\in L^p\left( 0,T;X \right) $ with the vector-valued distributional derivative $\partial _tf$ in $L^p\left( 0,T;X \right) $.
In particular,
we set $H^1\left( 0,T;X \right) =W^{1,2}\left( 0,T;X \right) $. The set of continuous functions $f:I\rightarrow X$ is denoted by
$\mathcal{C} \left( I;X \right)$. The space $\mathcal{C} \left( I;X_w \right) $ consists of all functions $f\in L^{\infty}\left( I;X \right) $ such that the map $t\in I\mapsto \left< \chi ,f\left( t \right) \right> $ is continuous, for all $\chi \in X ^ \prime$.
The Besov spaces denoted by $B_{p,\infty}^{s}\left(I; X \right) $,  with $s\in \left( 0,1 \right)$, are defined as
\begin{equation*}
B_{p,\infty}^{s}\left(I; X \right)=\left\{f\in L^p\left( I;X \right): \left\| f \right\| _{B_{p,\infty}^{s}\left( I;X \right)}<\infty\right\},
\end{equation*}
\begin{equation*}
\left\| f \right\| _{B_{p,\infty}^{s}\left( I;X \right)}=\left\| f \right\| _{L^p\left( I;X \right)}+\underset{0<h\leq 1}{\mathrm{sup}}h^{-s}\left\| \Delta _hf \right\| _{L^p\left( I_h;X \right)},
\end{equation*}
where $\Delta _hf\left( t \right) =f\left( t+h \right) -f\left( t \right) $ and $I_h=\left\{ t\in I:t+h\in I \right\} $.

\vskip2mm
Moreover, we need the following compactness and embedding results.
\begin{Lemma}{[cf. \cite{G-T, S-90}]}   \label{L3.1}  
Let $X, Y, Z$ be three Banach spaces such that $X\subset Y \subset Z$. Assume that $X\overset{c}{\hookrightarrow}Y$ and $Y\hookrightarrow Z$.
Then, for all $1\leq q \leq \infty $ and $0<\sigma <1$,  we have
\begin{align*}
&L^q\left( 0,T;X \right) \cap B_{q,\infty}^{\sigma}\left( 0,T;Z \right) \overset{c}{\hookrightarrow}L^q\left( 0,T;Y \right) ,
~~~\quad 1\le q<\infty,
\\
&\mathcal{C} \left( \left[ 0,T \right] ;X \right) \cap B_{q,\infty}^{\sigma}\left( 0,T;Z \right) \overset{c}{\hookrightarrow}\mathcal{C} \left( \left[ 0,T \right] ;Y \right), ~~\qquad  q=\infty .
\end{align*}
\end{Lemma}

\begin{Lemma}{[cf. \cite{G-T, S-66}]}  \label{L3.2}   
Let $X,Y$ be two Banach spaces such that $X\hookrightarrow Y$ and $Y^\prime \hookrightarrow X^\prime $ densely. Then,
$L^{\infty}\left( 0,T;X \right) \cap \mathcal{C} \left( \left[ 0,T \right] ;Y_w \right) \hookrightarrow \mathcal{C} \left( \left[ 0,T \right] ;X_w \right)$.
\end{Lemma}

\begin{Lemma}{[cf. \cite{L-17, L-96}]}  \label{Lions}
Let $T$ be a positive time, and assume that $\{(\rho_m, \boldsymbol{u}_m)\}_{m=1}^\infty$ satisfies
\begin{align*}
&\rho_m\in C([0, T]; L^1(\Omega)), \quad 0\le\rho_m\le C~~ {\rm a.e.~on}~\Omega\times(0, T),
\\
&{\rm div}\boldsymbol{u}_m=0~~ {\rm a.e.~on}~\Omega\times(0, T), \quad \|\boldsymbol{u}_m\|_{L^2(0, T; \mathbf{V}_{\sigma})}\le C,
\\
&\partial_t\rho_m+{\rm div}(\rho_m \boldsymbol{u}_m)=0 ~~ {\rm in}~\mathcal{D}^{'}(\Omega\times(0, T)),
\\
&\rho_m(0)\rightarrow \rho_0 ~~{\rm in}~L^1(\Omega), \qquad \boldsymbol{u}_m\rightharpoonup \boldsymbol{u}~~{\rm in}~L^2(0,T; \mathbf{V}_{\sigma})),
\end{align*}
where $C$ is a positive constant independent of $m$.
Then $\rho_m$ converges in $C([0, T]; L^r(\Omega))$, for $1\le r<\infty$, to the unique solution $\rho$, bounded on $\Omega\times(0, T)$, of
\begin{align*}
&\quad\partial_t\rho+{\rm div}(\rho\boldsymbol{u})=0 \quad {\rm in}~\mathcal{D}^{'}(\Omega\times(0, T)),
\\
&\rho\in C([0, T]; L^1(\Omega)), \quad \rho(0)=\rho_0~~ {\rm a.e.~in}~\Omega.
\end{align*}
\end{Lemma}

\section{Existence of weak solutions}
\noindent
In this section we prove the existence of global weak solutions to the problem (\ref{1.1})--(\ref{1.2}) with the Landau-type potential (\ref{1.4}) or the Flory-Huggins potential (\ref{1.5}).
\begin{Theorem}       
\label{Th4.1}
Let $\Omega$ be a bounded domain with smooth boundary in ${\mathbb R}^d~(d=2,3)$, and let $T$ be a positive time.
Assume that $(\rho_0, \boldsymbol{u}_0, \chi_0)\in L^\infty(\Omega)\times\mathbf{H}_{\sigma}\times H^1(\Omega)$
is given and such that $0 \leq \rho_{*} \leq \rho_{0} \leq \rho^{*}$, $|\chi_0|\le1$ a.e. in $\Omega$~
for some positive constants $\rho_{*}, \rho^{*}$.
Then, there exists a weak solution $(\rho,\boldsymbol{u},\chi)$ to the system \;$(\ref{1.1})$--$(\ref{1.2})$\; in the following sense:
\begin{enumerate}[(i)]  
\item The solution \;$(\rho,\boldsymbol{u},\chi)$\; satisfies  
\begin{align*}  
\begin{split}
&\rho \in \mathcal{C}\left( \left[ 0,T \right] ;L^r\left( \Omega \right) \right) \cap L^{\infty}\left(\Omega \times \left( 0,T \right) \right) \cap W^{1,\infty}( 0,T;H^{-1}\left( \Omega \right)  ),
\\
&\boldsymbol{u} \in \mathcal{C}\left( \left[ 0,T \right] ;\left( \mathbf{H}_{\sigma} \right) _{w} \right) \cap L^2\left( 0,T;\bf{V}_{\sigma} \right) \cap B_{2,\infty}^{\frac{1}{4}}\left( 0,T;\mathbf{H}_{\sigma}\right),
\\
&\chi \in \mathcal{C}( \left[ 0,T \right] ;( H^{1}\left( \Omega ) \right) _{w} ) \cap L^2( 0,T;H^{2}\left( \Omega \right)  )
 \cap H^1( 0,T;L^{\frac{3}{2}}\left( \Omega \right)),
\\
&
\mu \in L^2( 0,T;L^{2}\left( \Omega \right) )
\end{split}
\end{align*}
for any $ r \in \left[ 1,\infty \right)$. In addition,
\begin{equation*} 
0<\rho _{\ast}\le \rho (x,t)\le \rho ^{\ast}, \quad \left| \chi\left(x,t \right) \right|\le 1  \quad a.e.\ in~~\Omega \times (0,T).
\end{equation*}
\item The system \;$(\ref{1.1})$\; is satisfied as follows:
\begin{equation*}
\int_0^T{\int_{\Omega}{\rho}}\partial _t\psi \ \mathrm{d}x\mathrm{d}t+\int_0^T{\int_{\Omega}{\rho \boldsymbol{u}\cdot \nabla \psi \ \mathrm{d}x\mathrm{d}t=0}}
\end{equation*}
for all $\psi \in \mathcal{C} _{c}^{\infty}\left( \Omega \times \left( 0,T \right) \right)$;
\begin{equation*} 
\begin{split}
&-\int_0^T{\int_{\Omega}{\rho \boldsymbol{u}\cdot}}\partial _t\boldsymbol{w}\ \mathrm{d}x\mathrm{d}t-\int_0^T{\int_{\Omega}{\rho \boldsymbol{u}\otimes \boldsymbol{u}}}: \nabla \boldsymbol{w}\ \mathrm{d}x\mathrm{d}t
+\int_0^T{\int_{\Omega}{\eta \left( \chi \right) \mathbb{D} \boldsymbol{u}:\nabla}}\boldsymbol{w}\ \mathrm{d}x\mathrm{d}t
\\
&=\int_{\Omega}{\rho _0\boldsymbol{u}_0\cdot \boldsymbol{w}\left( 0 \right) \ \mathrm{d}x}+\int_0^T{\int_{\Omega}{\nabla \chi \otimes \nabla \chi}}: \nabla \boldsymbol{w}\ \mathrm{d}x\mathrm{d}t
\end{split}
\end{equation*}
for all $\boldsymbol{w}\in \mathcal{C} _{c}^{1}\left( \left[ 0,T \right) ;\mathbf{V}_{\sigma} \right) $;
\begin{equation*}  
\rho\partial _t\chi+\rho \boldsymbol{u}\cdot \nabla \chi =-m(\chi )\mu,
\\   
\quad \rho \mu =-\Delta \chi +\rho F^{\prime}\left( \chi \right) \quad a.e.\;in\;\Omega \times \left( 0,T \right) .
\end{equation*}
Furthermore, $\partial _{\boldsymbol{n}}\chi =0$ almost everywhere on $\partial \Omega \times \left( 0,T \right) $.
\item The initial data are attained in the following sense: $\rho \left( t \right) \rightarrow \rho \left( 0 \right) =\rho _0$ in $L^r\left( \Omega \right)$,\; for all $r\in \left[ 1,\infty \right) $, $ \boldsymbol{u}\left( t \right) \rightarrow \boldsymbol{u}\left( 0 \right) =\boldsymbol{u}_0$ in $L^2\left( \Omega \right)$, $\chi \left( t \right) \rightarrow \chi \left( 0 \right) =\chi _0$ in $H^1\left( \Omega \right)$, as $t\rightarrow 0^+$.
\item  The energy inequality
\begin{equation} \label{4.7}  
E_0\left(\rho (t),\boldsymbol{u}(t),\chi (t) \right )+\int_0^t{\int_{\Omega}{\left(\eta(\chi )|\mathbb{D} \boldsymbol{u}|^2+{m(\chi )|\mu|^2}\right)}}\mathrm{d}x\mathrm{d}\tau \le E_0(\rho _0,\boldsymbol{u}_0,\chi _0)
\end{equation}
holds for all $t\in \left[ 0,T \right] $,\ where
\begin{equation*}
E_0(\rho ,\boldsymbol{u},\chi )=\int_{\Omega}{\left(\frac{1}{2}\rho |\boldsymbol{u}|^2+\frac{1}{2}|\nabla \chi |^2+\rho F(\chi )\right)}\mathrm{d}x.
\end{equation*}

\end{enumerate}  

\end{Theorem}   
\begin{Remark}   \label{R4.2}
The pressure $p$ can be recovered in a classical way, see \cite{S-90 , T-77, L-96} for details.
\end{Remark}

\begin{Remark}
The same existence result holds for the Flory-Huggins logarithmic potential
\begin{align*}
\tilde F(s)=\frac{\theta}{2}\Big[(1+s)\ln(1+s)+(1-s)\ln(1-s)\Big]-\frac{\theta_0}{2}s^2, \quad s\in[-1,1],
\end{align*}
where $0<\theta<\theta_0$; the necessary modifications, which rely on the pointwise algebraic structure of the Allen-Cahn equation rather than the monotonicity arguments required in \cite{G-T}, are outlined in Appendix \ref{B-F-H}.
\end{Remark}
\noindent
{\bf Proof of Theorem \ref{Th4.1}.} The proof consists of four steps.
{We follow the semi-Galerkin framework of \cite[proof of Theorem 4.1]{G-T}, detailing only the estimates and limiting arguments specific to the Allen--Cahn coupling. }
\vskip1.8mm
\noindent
{\it\bfseries Step 1. Approximate problem.}\quad
For given $\rho_0 \in L^\infty (\Omega)$ satisfying
$0<\rho _{\ast}\le \rho_0\le \rho ^{\ast}$ a.e. in $\Omega$,
there exists a sequence $\left\{ \rho _{0m} \right\} _{m=1}^{\infty}\subset\mathcal{C} ^{\infty}( \overline{\Omega} )$ such that
$0<\rho _{\ast}\le \rho _{0m}\le \rho ^{*}$ in $\overline{\Omega}$ for all $m\in \mathbb{N}$, and
\begin{equation} \label{4.8} 
\rho _{0m}\rightarrow \rho _0~ \mathrm{strongly}\;\mathrm{in}\;L^r\left( \Omega \right) ~(r\in[1,\infty)), 
\qquad
\rho _{0m}\rightharpoonup \rho _0~{\rm weakly^{*}}\;\mathrm{in}~L^{\infty}\left( \Omega \right).
\end{equation}
Let $\{\boldsymbol{w}_i\} _{i=1}^{\infty}$ be the eigenfunctions of the Stokes operator $\mathbf{A}$,
$\mathbf{V}_m=\mathrm{span}\left\{ \boldsymbol{w}_1,\dots \boldsymbol{w}_m \right\}$ and
let $\mathbb{P} _{\mathrm{m}}$ be the orthogonal projection onto $\mathbf{V}_m$ in $\mathbf{H}_{\sigma}$.
For $\boldsymbol{u}_{0m}=\mathbb{P} _{\mathrm{m}}\boldsymbol{u}_{0}$ there holds
\begin{equation} \label{4.9} 
\boldsymbol{u} _{0m}\rightarrow \boldsymbol{u} _0~~ \mathrm{strongly}\;\mathrm{in}~\;\mathbf{H}_{\sigma}.
\end{equation}

We are going to solve the following approximate problem
\begin{equation} \label{A.1}   
\begin{cases}
\partial _t\rho _m+\boldsymbol{u}_m\cdot \nabla \rho _m=0 \quad \mathrm{in}\;\Omega \times (0,T),
\\
\left( \rho _m\partial _t\boldsymbol{u}_m,\boldsymbol{w} \right)
+\left( \rho _m\left( \boldsymbol{u}_m\cdot \nabla \right) \boldsymbol{u}_m,\boldsymbol{w} \right)
 +\left( \eta \left( \chi _m \right) \mathbb{D} \boldsymbol{u}_m,\nabla \boldsymbol{w} \right)
=\left( \nabla \chi _m\otimes \nabla \chi _m,\nabla \boldsymbol{w} \right), \forall\boldsymbol{w}\in \mathbf{V}_m,
\\  
\rho _m\partial _t\chi _m+\rho _m\boldsymbol{u}_m\cdot \nabla \chi _m =-m(\chi _m)\mu _m \quad \mathrm{in}\;\Omega \times (0,T),
\\ 
\rho _m\mu _m=-\Delta \chi _m+\rho _m F^{\prime}\left( \chi _m \right)\quad \mathrm{in}\;\Omega \times (0,T),
\\
\boldsymbol{u}_m=0,\quad \partial _{\mathbf{n}}\chi _m=0	\quad	\mathrm{on}\;\partial \Omega \times (0,T),
\\
\rho _m(\cdot ,0)=\rho _{0m}, ~\boldsymbol{u} _m(\cdot ,0)=\boldsymbol{u} _{0m}, ~\chi _m(\cdot ,0)=\chi _0	\quad	\mathrm{in}\;\Omega.
\end{cases}
\end{equation}

\noindent
{\it\bfseries Step 2. Existence of approximate solutions.}\quad
The existence of approximate solutions $( \rho _m,\boldsymbol{u}_m,\chi _m)$ to (\ref{A.1}) satisfying
\begin{align*}   
&\rho _m\in \mathcal{C} ^1( \overline{Q}_T ), ~~ \quad
\boldsymbol{u}_m\in \mathcal{C}( [ 0,T ] ;\mathbf{V}_m )\cap W^{1, p}(0, T; \mathbf{V}_m), \\
&\chi _m\in \mathcal{C} ([ 0,T] ;{H^1( \Omega )} )\cap{ L^2( 0,T;H^2\left( \Omega \right)) }
\end{align*}
for any $p>1$ and $T>0$, where $Q_T = \Omega \times (0,T)$, is proved by a fixed point argument in Appendix \ref{A-A}.
\vskip1mm
\noindent
{\it\bfseries Step 3. Estimates uniform in $m$.}\quad {
First, the transport equation gives $\rho_* \le \rho_m \le \rho^*$ in $Q_T$ by the standard characteristics (see Appendix \ref{A-A}),
and the maximum principle applied to the approximate Allen--Cahn equation gives $|\chi_m| \le 1$ in $Q_T$ (see Appendix \ref{A-A}), that is }
\begin{equation}   \label{4.17}   
\rho_{*} \leq \rho_{m}(x, t) \leq \rho^{*}, \qquad { |\chi_m|\le1 \quad \text{in}~Q_{T} }.
\end{equation}
Taking $\boldsymbol{w}=\boldsymbol{u}_m$ in the momentum equation (\ref{A.1})$_2$, multiplying the Allen-Cahn equation (\ref{A.1})$_3$ by $\mu_m$, and integrating over $Q_t$ yields the energy identity
\begin{equation}   \label{4.19}   
E_0\left( \rho _m(t),\boldsymbol{u}_m(t),\chi _m(t) \right)
+\iint_{Q_{t}}{{\left(\eta \left( \chi _m \right) \left| \mathbb{D} \boldsymbol{u}_m \right|^2+{m\left( \chi _m \right) \left|  \mu _m \right|^2}\right)}}\mathrm{d}x\mathrm{d}\tau
=E_0\left( \rho _{0m},\boldsymbol{u}_{0m},\chi _0 \right)
\end{equation}
for all $t \in[0, T]$.
{Recalling $F\geq 0$, and passing to the limit in the initial energy, we deduce from \eqref{4.19} that}
\begin{align} \label{4.20}
 \left\|\boldsymbol{u}_{m}\right\|_{L^{\infty}\left(0, T ; \mathbf{H}_{\sigma}\right)}
+\left\|\boldsymbol{u}_{m}\right\|_{L^{2}\left(0, T ; \mathbf{V}_{\sigma}\right)}
+\left\|\nabla \chi_{m}\right\|_{L^{\infty}\left(0, T ; L^{2}(\Omega)\right)}
+\left\|\mu_{m}\right\|_{L^{2}\left(0, T ; L^{2}(\Omega)\right)}
\leq C(E_{0}),
\end{align}
where $E_{0}=E_{0}\left(\rho_{0}, \boldsymbol{u}_{0}, \chi_{0}\right)$.
Here and below, we denote by $C$ a constant whose value may vary from line to line but is independent of $m$.
{ Since $|\chi_m| \le 1$ by \eqref{4.17} and $\nabla \chi_m$ is bounded in $L^\infty(0, T; L^2(\Omega))$ by \eqref{4.20}, there holds}
\begin{align}   \label{4.24}   
\left\|\chi_{m}\right\|_{L^{\infty}\left(0, T ; H^{1}(\Omega)\right)} \leq { C\left(E_{0}\right)}.
\end{align}
Multiplying $(\ref{A.1})_4$ by $-\Delta \chi _m$ and integrating over $\Omega$, { since $|F'(\chi_m)|\le C$ by \eqref{4.17}, we get
$\|\Delta\chi_m\|_{L^2}\le\rho^*(\|\mu_m\|_{L^2}+C)$; hence by Neumann elliptic regularity}
\begin{align}   \label{4.27}
\left\| \chi _m \right\| _{L^{{2}}\left( 0,T;H^2(\Omega ) \right)}\le { C\left(E_0,T \right)}.
\end{align}
Moreover, by the Gagliardo-Nirenberg inequality, { $(\ref{4.20})$ and (\ref{4.27}), we deduce that
\begin{align}\label{u-phi-43}
\|\boldsymbol{u}_m\|_{L^4(0,T;L^3(\Omega ))}
+\|\nabla\chi _m\|_{L^4(0,T;L^3(\Omega))}
\le C\left(E_0,T \right).
\end{align} }
Moreover, it follows directly from the Allen-Cahn equation (\ref{A.1})$_3$, (\ref{4.17}), (\ref{u-phi-43}), and the bound above that
\begin{align}
\label{phi-t}
\| \partial_t\chi _m \| _{L^{2}(0,T;L^{\frac32}(\Omega))}
\le& C\|\boldsymbol{u}_{m}\cdot\nabla\chi_m \| _{L^{2}(0,T;L^{\frac32}(\Omega))}+C\|\mu_m\|_{L^{2}(0,T;L^{\frac32}(\Omega))}
\\
\le& C\|\boldsymbol{u}_{m}\| _{L^{4}(0,T;L^{3}(\Omega))}\|\nabla\chi_m \| _{L^{4}(0,T;L^{3}(\Omega))}+C\|\mu_m\|_{L^{2}(0,T;L^{2}(\Omega))}
\le C. \nonumber
\end{align}
Let $v \in H_{0}^{1}(\Omega)$ with $\|v\|_{H_{0}^{1}(\Omega)} \leq 1$. 
Multiplying the transport equation $(\ref{A.1})_1$ by $v$ and integrating over $\Omega$, we get
$\left. \langle \partial _t\rho _m,v \right. \rangle =\int_{\Omega}{\rho _m}\boldsymbol{u}_m\nabla v\mathrm{d}x$;
using $(\ref{4.17})$ and $(\ref{4.20})$, we arrive at
\begin{align}   \label{4.28}   
\left\| \partial _t\rho _m \right\| _{L^{\infty}\left( 0,T;H^{-1}(\Omega ) \right)}\le C\left( E_0 \right) .
\end{align}

Since the available estimates do not provide a uniform bound on $\partial_t\boldsymbol{u}_{m}$ suitable for compactness, we resort to the translation differences $\boldsymbol{u}_{m}(t+h)-\boldsymbol{u}_{m}(t)$, as in \cite{S-90, A-K-M, L-78, T-95}.
By using $(\ref{A.1})_1$, integrating $(\ref{A.1})_2$ on the time interval $(t, t+h)$, where $t \in[0, T-h]$,
and taking $\boldsymbol{w}=\boldsymbol{u}_{m}(t+h)-\boldsymbol{u}_{m}(t)$, we have
\begin{align}  \label{4.31}   
&\int_{\Omega}{\rho _m}(t+h)\left| \boldsymbol{u}_m(t+h)-\boldsymbol{u}_m(t) \right|^2\mathrm{d}x \nonumber
\\
&=-\int_{\Omega}\left( \rho _m(t+h)-\rho _m(t) \right) \boldsymbol{u}_m(t)\cdot \left( \boldsymbol{u}_m(t+h)-\boldsymbol{u}_m(t) \right) \mathrm{d}x \nonumber
\\
&~~~-\int_t^{t+h}{\int_{\Omega}}\mathrm{div}\left( \rho _m(\tau )\boldsymbol{u}_m(\tau )\otimes \boldsymbol{u}_m(\tau ) \right) \cdot \left( \boldsymbol{u}_m(t+h)-\boldsymbol{u}_m(t) \right) \mathrm{d}x\mathrm{d}\tau
\\
&~~~-\int_t^{t+h}{\int_{\Omega}}\eta \left( \chi _m(\tau ) \right) \mathbb{D} \boldsymbol{u}_m(\tau ):\nabla \left( \boldsymbol{u}_m(t+h)-\boldsymbol{u}_m(t) \right) \mathrm{d}x\mathrm{d}\tau \nonumber
\\
&~~~+\int_t^{t+h}{\int_{\Omega}{{\nabla \chi _m(\tau )\otimes \nabla \chi _m(\tau )}:\nabla \left( \boldsymbol{u}_m(t+h)-\boldsymbol{u}_m(t) \right)}}\mathrm{d}x\mathrm{d}\tau  \nonumber
=\sum_{i=1}^{4}{M}_i(t). \nonumber
\end{align}
{
The terms $M_1$, $M_2$, $M_3$ are estimated exactly as $I_1$, $I_2$, $I_3$ in
\cite[proof of Theorem 4.1, estimates (4.34)--(4.38)]{G-T},
using the transport equation, the density bounds
\eqref{4.17}, the energy bounds \eqref{4.20},
and the uniform upper bound $\eta(\cdot)\le\eta^*$;
this yields
\begin{align} \label{4.32}
\int_0^{T-h} |M_1(t)| \, dt \le C(E_0,T)h^{3/4}, \qquad \int_0^{T-h} \big( |M_2(t)| + |M_3(t)| \big) \, dt \le C(E_0,T)h^{1/2}.
\end{align}
It remains to estimate the capillary term $M_4$, which replaces the terms $I_4$, $I_5$ of \cite{G-T} and is treated as a whole,
since only $\mu_m \in L^2(0, T; L^2(\Omega))$ is available here. By the Gagliardo--Nirenberg inequality, \eqref{4.17} and $(\ref{4.20})$, }
\begin{align*}      
\left| {{M}_4(t)} \right|&\le \int_t^{t+h}{\left\| \nabla \chi _m(\tau ) \right\| _{L^4(\Omega )}^{2}}\mathrm{d}\tau \left( \left\| \nabla \boldsymbol{u}_m(t+h) \right\| _{L^2(\Omega )}+\left\| \nabla \boldsymbol{u}_m(t) \right\| _{L^2(\Omega )} \right)
\\
&\le C\left( E_0 \right) \int_t^{t+h}{\left\| \Delta \chi _m(\tau ) \right\| _{L^2(\Omega )}^{\frac{3}{2}}}\mathrm{d}\tau \left( \left\| \nabla \boldsymbol{u}_m(t+h) \right\| _{L^2(\Omega )}+\left\| \nabla \boldsymbol{u}_m(t) \right\| _{L^2(\Omega )} \right).
\end{align*}
{Integrating in $t$, using Fubini's theorem, and then the Cauchy-Schwarz inequality over the inner integral of length at most $h$, we obtain
\begin{align}  \label{4.35}
\int_0^{T-h}{\left|{M}_4(t) \right|}\mathrm{d}t
\leq C(E_0) \, h^{\frac{1}{2}} \|\nabla \boldsymbol{u}_m\|_{L^2(Q_T)} \int_0^T \|\Delta \chi_m(\tau)\|_{L^2(\Omega)}^{\frac{3}{2}} \,\mathrm{d}\tau \leq C(E_0, T) \, h^{\frac{1}{2}},
\end{align}
where the last line follows from \eqref{4.20}, \eqref{4.27} and H\"{o}lder's inequality. }
Recalling $(\ref{4.17})$, integrating $(\ref{4.31})$ on $(0, T-h)$, and using $(\ref{4.32})$, $(\ref{4.35})$, we infer that
\begin{align*}    
\int_0^{T-h}{\left\| \boldsymbol{u}_m(t+h)-\boldsymbol{u}_m(t) \right\| _{L^2(\Omega )}^{2}}\mathrm{d}t\le C\left(E_0,T \right) h^{\frac{1}{2}}.
\end{align*}
This implies that
\begin{align}    \label{4.45}   
\left\| \boldsymbol{u}_m \right\| _{B_{2,\infty}^{\frac{1}{4}}\left( 0,T;\mathbf{H}_{\sigma} \right)}
\le C(E_0,T),
\end{align}
where $C$ is a positive constant independent of $m$ and $h$.
\vskip1.8mm
\noindent
{\it\bfseries Step 4.  Passage to the limit.}\quad
Thanks to $(\ref{4.17})$, $(\ref{4.20})$, $(\ref{4.24})$, $(\ref{4.27})$, $(\ref{phi-t})$ and $(\ref{4.28})$,
we infer that, up to a subsequence, the corresponding weak convergences hold.
Using $(\ref{4.20})$, $(\ref{4.27})$, $(\ref{phi-t})$, $(\ref{4.45})$, Lemma $\ref{L3.1}$ and the Aubin--Lions lemma, we have
\begin{align}         
\boldsymbol{u}_{m} \rightarrow \boldsymbol{u} & \quad \text { strongly in } L^{2}\left(0, T ; L^{q}(\Omega)\right), \forall q \in[2,6),\label{4.47}\\
\chi_{m} \rightarrow \chi & \quad \text { strongly in } L^{{2}}(0, T ; W^{1, q}(\Omega)), \forall q \in[2,6),\label{4.48} \\
\chi_{m} \rightarrow \chi & \quad \text { strongly in } \mathcal{C}\left([0, T] ; L^{q}(\Omega)\right), \forall q \in[2,6). \label{4.49}
\end{align}
The strong convergence of the density follows from Lemma \ref{Lions}
\begin{align}      \label{4.50}   
\rho_{m} \rightarrow \rho \quad \text { strongly in } \mathcal{C}\left([0, T] ; L^{r}(\Omega)\right), \forall r \in[1, \infty).
\end{align}
Thanks to the above estimates and convergence results, we obtain
\begin{alignat}{2}      \label{4.51}   
&\rho_{m} \boldsymbol{u}_{m} \otimes \boldsymbol{u}_{m} \rightharpoonup \rho \boldsymbol{u} \otimes \boldsymbol{u} && \text { weakly in } L^{{\frac{4}{3}}}(0, T ; L^{2}(\Omega)),  \nonumber
 \\
&\eta \left(\chi_{m}\right) \mathbb{D} \boldsymbol{u}_{m} \rightharpoonup \eta(\chi) \mathbb{D} \boldsymbol{u} && \text { weakly in } L^{2}(Q_T),
\\
&{ {\nabla \chi_{m} \otimes \nabla \chi_{m} \rightarrow \nabla \chi \otimes \nabla \chi}} 
~\quad && { \text { strongly in } L^{1}(Q_T)  \text { and weakly in } L^{{\frac{4}{3}}}(0, T ; L^{2}(\Omega)) }, \nonumber
 \\
& {\rho_{m} F^{\prime}\left(\chi_{m}\right) \rightarrow \rho F^{\prime}(\chi)}
 &&{ \text { strongly in } L^{2}(Q_T)}.  \nonumber
\end{alignat}
{The weak convergence of the capillary stress follows from its strong $L^1(Q_T)$ convergence and its uniform boundedness in $L^{4/3}(0, T; L^2(\Omega))$, obtained by interpolation from \eqref{4.24} and \eqref{4.27}.
}

With the convergences $(\ref{4.51})$ at hand, we pass to the limit in the Galerkin formulations. The continuity equation follows directly from Lemma \ref{Lions}. For the momentum equation, we fix $l \in \mathbb{N}$, and take $m \geq l$ in
\begin{align}       \label{4.53}   
	\left( \partial _t\left( \rho _m\boldsymbol{u}_m \right) ,\boldsymbol{w} \right) +\left( \mathrm{div}\left( \rho _m\boldsymbol{u}_m\otimes \boldsymbol{u}_m \right) ,\boldsymbol{w} \right) &+\left( \eta \left( \chi _m \right) \mathbb{D} \boldsymbol{u}_m,\nabla \boldsymbol{w} \right)
=\left({ \nabla \chi _m\otimes \nabla \chi _m}, \nabla \boldsymbol{w} \right),
\end{align}
for all $\boldsymbol{w} \in \mathbf{V}_{l}$.
Multiply $(\ref{4.53})$ by a test function $\psi \in \mathcal{C}_{c}^{1}([0, T))$ and integrate in time. Letting $m \rightarrow \infty$ and then $l \rightarrow \infty$ by the density of $\cup_{m}\mathbf{V}_m$ in $\mathbf{V}_\sigma$ yields the weak formulation of the momentum equation.
Finally, by \eqref{4.17}, \eqref{u-phi-43}--\eqref{phi-t}, and \eqref{4.47}--\eqref{4.50},
\[
 \rho_m\partial_t\chi_m \rightharpoonup \rho\partial_t\chi,
 \qquad
 \rho_m u_m\cdot\nabla\chi_m
 \rightharpoonup \rho u\cdot\nabla\chi
 \qquad { \text{in }L^2(0,T;L^{3/2}(\Omega))}.
\]
{ Here the limits are identified using the strong convergences
and the uniform boundedness of $\rho_m$.
Since $m(\chi_m)\mu_m\rightharpoonup m(\chi)\mu$ in $L^2(Q_T)$,
we can pass directly to the limit in \eqref{A.1}$_3$,
obtaining $\rho\partial_t\chi+\rho\boldsymbol{u}\cdot\nabla\chi=-m(\chi)\mu$.
Similarly, $\rho_m\mu_m \rightharpoonup \rho\mu$ in $L^2(Q_T)$ (by the strong convergence of $\rho_m$ in $\mathcal{C}([0, T]; L^r(\Omega))$ and the weak convergence of $\mu_m$), and $-\Delta \chi_m \rightharpoonup -\Delta \chi$ in $L^2(Q_T)$; together with the last line of \eqref{4.51} this gives $\rho\mu = -\Delta \chi + \rho F'(\chi)$ a.e. in $\Omega \times (0, T)$.
The boundary condition $\partial_{\mathbf n}\chi=0$
is inherited from the weak convergence in $L^2(0,T;H^2(\Omega))$ and the continuity of the normal trace. }

\noindent{\it Time continuity.}\quad 
Since $\rho \in \mathcal{C}\left([0, T] ; L^{r}(\Omega)\right)$ for all $r \in[1, \infty)$, $\chi \in \mathcal{C}\left([0, T] ; L^{q}(\Omega)\right)$ for all $q \in[2,6)$, and $\chi \in L^{\infty}(0, T ; H^{1}(\Omega))$, Lemma $\ref{L3.2}$ gives $\chi \in \mathcal{C}([0, T] ; (H^{1}(\Omega))_{w})$. Next, we rewrite $(\ref{4.53})$ as follows
\begin{align*}   
\left. \langle \partial _t\mathbb{P} _m\left( \rho _m\boldsymbol{u}_m \right) ,\boldsymbol{w} \right. \rangle =\left( \rho _m\boldsymbol{u}_m\otimes \boldsymbol{u}_m,\nabla \mathbb{P} _m\boldsymbol{w} \right) -\left( \eta \left( \chi _m \right) \mathbb{D} \boldsymbol{u}_m,\nabla \mathbb{P} _m\boldsymbol{w} \right) +\left({\nabla \chi _m\otimes \nabla \chi _m},\nabla \mathbb{P} _m\boldsymbol{w} \right),
\end{align*}
for all $\boldsymbol{w} \in \mathbf{V}_{\sigma}$. { By using $(\ref{4.20})$, $(\ref{4.27})$ and the Gagliardo-Nirenberg inequality, we obtain
\begin{align*}
	\left|  \langle \partial _t\mathbb{P} _m\left( \rho _m\boldsymbol{u}_m \right) , \boldsymbol{w}  \rangle \right|
\le &C\left( \left\| \nabla \boldsymbol{u}_m \right\| _{L^2(\Omega )}^{\frac{3}{2}}+\left\| \mathbb{D} \boldsymbol{u}_m \right\| _{L^2(\Omega )}+{\left\| \Delta \chi _m \right\| _{L^2(\Omega )}^{\frac{3}{2}}} \right) \left\| \nabla \boldsymbol{w}\right\| _{L^2(\Omega )},
\end{align*}
so that $\partial_{t} \mathbb{P}_{m}\left(\rho_{m} \boldsymbol{u}_{m}\right)$ is bounded in
$L^{\frac{4}{3}}\left(0, T ; \mathbf{V}_{\sigma}^{\prime} \right)$.
By \eqref{4.17}, \eqref{4.20}, \eqref{4.47}, and
\eqref{4.50}, $\rho_m\boldsymbol{u}_m$ is bounded in
$L^\infty(0,T;L^2(\Omega))$ and converges to
$\rho\boldsymbol{u}$ strongly in $L^2(0,T;L^2(\Omega))$.
Together with the above time-derivative bound and the compact
embedding $\mathbf{H}_\sigma\Subset\mathbf{V}_\sigma'$, this yields,
up to a subsequence,
\[
 \mathbb{P}_m(\rho_m\boldsymbol{u}_m)
 \longrightarrow \mathbb{P}(\rho\boldsymbol{u})
 \quad\text{in }\mathcal{C}([0,T];\mathbf{V}_\sigma').
\]
Hence Lemma~\ref{L3.2} gives
$\mathbb{P}(\rho\boldsymbol{u})
\in\mathcal{C}([0,T];(\mathbf{H}_\sigma)_w)$.
The weak time continuity of $\boldsymbol{u}$ and
$\sqrt{\rho}\boldsymbol{u}$ then follows from the Helmholtz
decomposition argument in
\cite[proof of Theorem 4.1, ``Time continuity"]{G-T}.
}

\vskip1mm
\noindent{\it Energy inequality.} 
The energy inequality $(\ref{4.7})$ follows from multiplying $(\ref{4.19})$ by a nonnegative test function, integrating in time, and passing to the limit: { the kinetic and gradient energies and the dissipative terms are handled by weak lower semicontinuity, while the potential energy passes to the limit directly, since $\rho_m F(\chi_m) \to \rho F(\chi)$ strongly in $L^1(Q_T)$ by \eqref{4.49}, \eqref{4.50}, $|\chi_m| \le 1$ and dominated convergence.
This gives $(\ref{4.7})$ for almost every $t$. The inequality holds in fact for every $t \in [0, T]$: the kinetic and gradient energies are lower semicontinuous in time, by $\sqrt{\rho}\boldsymbol{u} \in \mathcal{C}([0, T]; (L^2(\Omega))_w)$ and $\chi \in \mathcal{C}([0, T]; (H^1(\Omega))_w)$; the potential energy is continuous in time, by $\rho \in \mathcal{C}([0, T]; L^r(\Omega))$, $\chi \in \mathcal{C}([0, T]; L^2(\Omega))$ and $|\chi| \le 1$; and the dissipative integrals are continuous in their upper limit.
}

\vskip1mm
\noindent{\it The initial data.} {
By \eqref{4.50}, \eqref{4.49}, and the approximate initial data,
$\rho(0)=\rho_0$ and $\chi(0)=\chi_0$.
The limiting momentum equation identifies
$\boldsymbol{u}(0)=\boldsymbol{u}_0$ as in
\cite[proof of Theorem 4.1, ``The initial data"]{G-T}.
The weak time continuity established above, the energy
inequality \eqref{4.7}, and the continuity of the potential
energy imply convergence of the sum of the kinetic and
gradient energies to its initial value.
Consequently, weak convergence and convergence of norms give}
\[
 \lim_{t\to0^+}
 \left(
 \left\|\sqrt{\rho(t)}\,\boldsymbol{u}(t)
       -\sqrt{\rho_0}\,\boldsymbol{u}_0\right\|_{L^2(\Omega)}^2
 +\left\|\nabla\chi(t)-\nabla\chi_0\right\|_{L^2(\Omega)}^2
 \right)=0.
\]
{ Since $\rho(t)\ge\rho_*>0$ and
$\rho(t)\to\rho_0$ in $L^r(\Omega)$ for every $r<\infty$,
this yields $\boldsymbol{u}(t)\to\boldsymbol{u}_0$
in $L^2(\Omega)$.
Together with $\chi(t)\to\chi_0$ in $L^2(\Omega)$,
we also obtain $\chi(t)\to\chi_0$ in $H^1(\Omega)$.
}
\hfill$\Box$

\section{Existence of strong solutions}
This section is devoted to the existence of strong solutions to the system $(\ref{1.1})$--$(\ref{1.2})$.
We prove that these solutions are global in the 2D case and local in the 3D case.
\begin{Definition}
For $T>0$, $(\rho, \boldsymbol{u}, p, \chi)$ is called a strong solution to the problem $(\ref{1.1})$--$(\ref{1.2})$, if
\begin{align*}
& \rho \in \mathcal{C}\left([0, T] ; L^{r}(\Omega)\right) \cap L^{\infty}(\Omega \times(0, T)) \cap L^{\infty}(0, T ; H^{-1}(\Omega)),
\\
& \boldsymbol{u} \in \mathcal{C}\left([0, T] ; \mathbf{V}_{\sigma}\right) \cap L^{2}(0, T ; H^{2}(\Omega)) \cap H^{1}\left(0, T ; \mathbf{H}_{\sigma}\right), \\
& p \in L^{2}(0, T ; H^{1}(\Omega)),
\\
& \chi \in \mathcal{C}([0, T]; ((H^{2}(\Omega))_w ) \cap H^{1}(0, T ; H^{1}(\Omega))
\end{align*}
for any $r \in[1, \infty)$, and $(\ref{1.1})_1$ holds in $\mathcal{D}'(Q_T)$, $(\ref{1.1})_{2-5}$ hold a.e.\ in $Q_T$,
$\partial_{\boldsymbol{n}} \chi=0$ holds on $\partial \Omega\times(0, T)$.
\end{Definition}
\begin{Theorem}       
\label{Th6.1}
Let $\Omega$ be a bounded domain with smooth boundary in $\mathbb{R}^{d}$ ($d=2, 3$).
Assume that $(\rho_0, \boldsymbol{u}_0, \chi_0)\in L^\infty(\Omega)\times\mathbf{V}_{\sigma}\times H^2(\Omega)$
is given and such that $0 < \rho_{*} \leq \rho_{0} \leq \rho^{*}$, $|\chi_0|\le1$ a.e. in $\Omega$
and $\partial_{\boldsymbol{n}} \chi_{0}=0$ on $\partial \Omega$. Then, we have
\begin{enumerate}[(i)]  
\item If $d=2$ and $\eta\equiv1$, for any $T>0$ there exists a strong solution $(\rho, \boldsymbol{u}, p, \chi)$
to the system $(\ref{1.1})$--$(\ref{1.2})$ in $Q_T$.
\item If $d=3$, there exist a time $T_{0}>0$, depending on the norms of the initial data, and a strong solution
$(\rho, \boldsymbol{u}, p, \chi)$ to the system $(\ref{1.1})$--$(\ref{1.2})$ in $Q_{T_0}$.
\end{enumerate}
\end{Theorem}   

\begin{Remark}
The assumption $\eta\equiv 1$ in Theorem \ref{Th6.1}(i) is technically necessary. For the NSCH system in \cite{G-T}, the low-order energy yields $\mu_m\in L^4(L^2)\cap L^2(H^1)$, and the high-order energy inequality has $\|\nabla\mu_m\|_{L^2}^2$ as an energy term on the left-hand side, which is crucial for controlling the variable-viscosity terms $I_{15}$ and $K_{4,5}$ in the estimates below.
In contrast, for the NSAC system, we have $\mu_m\in L^2(L^2)$, with $\|\rho_m\partial_t\chi_m\|_{L^2}^2$ replacing $\|\nabla\mu_m\|_{L^2}^2$ as the principal energy term (see (\ref{6.9}) below). Consequently, the techniques for controlling variable viscosity in \cite{G-T} do not directly apply to the present system, and we assume $\eta\equiv 1$ for global solutions in 2D. In 3D, variable viscosity is handled locally via the high-power estimate (\ref{6.40}).
{Related higher-order estimates for a conserved NSAC system with
nonconstant viscosity and a singular potential were obtained by
Grasselli, Hurm, and Poiatti \cite{G-H-P} in their analysis of local
strong solutions in 3D. }
\end{Remark}

\noindent
{\bf Proof of Theorem \ref{Th6.1}.} The proof is divided into three parts.
\vskip1mm
\noindent
{\it\bfseries Step 1. Approximate regular solutions}\quad
We consider the approximate solutions constructed in the proof of Theorem \ref{Th4.1}.
More precisely, there exist
\begin{align*}   
&\rho _m\in \mathcal{C} ^1( \overline Q_T ) , ~~ \boldsymbol{u}_m\in \mathcal{C}( [ 0,T ] ;\mathbf{V}_m )\cap W^{1, p}(0, T; \mathbf{V}_m),  ~~
 \chi _m\in \mathcal{C} ([ 0,T] ;{H^1( \Omega )} )\cap{ L^2( 0,T;H^2\left( \Omega \right)) }
\end{align*}
for any $p>1$.
The approximate densities $\rho_{m}(x, t)=\rho_{0m}\left(\mathbf{X}_{m}(0, t, x)\right)$, where the characteristic $\mathbf{X}_{m}(s, t, x)$ is given by $(\ref{A.4})$ with $\boldsymbol{v}_m$ replaced by $\boldsymbol{u}_m$, solve
\begin{align} \label{6.1}  
\partial _t\rho _m+\boldsymbol{u}_m\cdot \nabla \rho _m=0 \quad \mathrm{in}\;\Omega \times (0,T).
\end{align}
The approximate solutions $\left(\rho_{m}, \boldsymbol{u}_{m}, \chi_{m}\right)$ satisfy
\begin{align}  
&\left( \rho _m\partial _t\boldsymbol{u}_m,\boldsymbol{w} \right) +\left( \rho _m\left( \boldsymbol{u}_m\cdot \nabla \right) \boldsymbol{u}_m,\boldsymbol{w} \right) +\left( \eta \left( \chi _m \right) \mathbb{D} \boldsymbol{u}_m,\nabla \boldsymbol{w} \right)
={\left( \nabla \chi _m\otimes \nabla \chi _m,\nabla \boldsymbol{w} \right)}, \label{6.2}
\\  
&{\rho _{m}^{2}\partial _t\chi _m+\rho _{m}^{2}\boldsymbol{u}_m\cdot \nabla \chi _m = m(\chi_m)\Delta \chi _m-\rho _m m(\chi_m)F ^\prime\left( \chi _m \right)} \quad \mathrm{in}\;\Omega \times (0,T) \label{6.3}  
\end{align}
for all $\boldsymbol{w}\in \mathbf{V}_m$ and $t\in [0, T]$.
The initial conditions are {$\rho_m(\cdot, 0)=\rho_{0m}$, $\boldsymbol{u}_{m}(\cdot, 0)=\mathbb{P}_{m} \boldsymbol{u}_{0}$, $\chi_{m}(\cdot, 0)=\chi_{0}$}. Moreover, recalling that $\rho_{0m} \rightarrow \rho_{0}$ strongly in $L^{r}(\Omega), \forall r \in[1, \infty)$, $\rho_{0m} \rightharpoonup \rho_{0}$ weakly$^*$ in $L^{\infty}(\Omega)$ and $\boldsymbol{u}_{m}(\cdot, 0) \rightarrow \boldsymbol{u}_{0}$ in $\mathbf{V}_{\sigma}$, it follows from $(\ref{4.17})$, $(\ref{4.19})$, $(\ref{4.24})$, and $(\ref{4.27})$ that
\begin{align}      \label{6.4}   
0<\rho_{*} \leq \rho_{m}(x, t) \leq \rho^{*}, \quad |\chi_m(x, t)|\le 1 \quad \forall(x, t) \in \overline Q_{T}
\end{align}
\begin{align}      \label{6.5}   
\left\|\boldsymbol{u}_{m}\right\|_{L^{\infty}\left(0, T ; \mathbf{H}_{\sigma}\right)} \leq C\left(E_{0}\right), \quad\left\|\boldsymbol{u}_{m}\right\|_{L^{2}\left(0, T ; \mathbf{V}_{\sigma}\right)} \leq C\left(E_{0}\right)
\end{align}
and
\begin{align}      \label{6.6}   
\left\|\chi_{m}\right\|_{L^{\infty}\left(0, T ; H^{1}(\Omega)\right)} \leq C\left(E_{0}\right),
\quad \left\|\chi_{m}\right\|_{L^{2}\left(0, T ; H^{2}(\Omega)\right)} \leq C\left(E_{0}\right),
\end{align}
where $C$ is independent of $m$, and $E_{0}=E_{0}\left(\rho_{0}, \boldsymbol{u}_{0}, \chi_{0}\right)$.
\vskip1.8mm
\noindent
{\it\bfseries Step 2. Higher-order energy inequalities}\quad
Differentiating the equation $(\ref{6.3})$ with respect to $t$,
multiplying the resulting equation by $\partial _t\chi _m$ and integrating over $\Omega$, we can deduce that
\begin{align}           \label{6.9phi}   
	&\frac{\mathrm{d}}{\mathrm{d}t}\int\limits_{\Omega}{\frac{1}{2}\rho _{m}^{2}\left| \partial _t\chi _m \right|^2}\mathrm{d}x+\int\limits_{\Omega}{m(\chi _m)\left| \nabla \partial _t\chi _m \right|^2}\mathrm{d}x
\nonumber\\
	=&\frac{1}{2}\int\limits_{\Omega}{\left| \partial _t\chi _m \right|^2\boldsymbol{u}_m\cdot {\nabla \rho _{m}^{2}}}\mathrm{d}x
+\int\limits_{\Omega}{\partial _t\chi _m\left( \boldsymbol{u}_m\cdot {\nabla \rho _{m}^{2}} \right) \left( \boldsymbol{u}_m\cdot \nabla \chi _m \right)}\mathrm{d}x
-\int\limits_{\Omega}{\rho _{m}^{2}\partial _t\chi _m\partial _t\boldsymbol{u}_m\cdot \nabla \chi _m}\mathrm{d}x\nonumber\\
	&-\int\limits_{\Omega}{\rho _{m}^{2}\partial _t\chi _m\boldsymbol{u}_m\cdot \nabla \partial _t\chi _m}\mathrm{d}x+\int\limits_{\Omega}{m^{\prime}(\chi _m)\left| \partial _t\chi _m \right|^2\Delta \chi _m}\mathrm{d}x-\int\limits_{\Omega}{m^{\prime}(\chi _m)\partial _t\chi _m\nabla \partial _t\chi _m\cdot\nabla \chi _m}\mathrm{d}x\nonumber\\
	&-\int\limits_{\Omega}{\rho _mm^{\prime}(\chi _m)F^{\prime}(\chi _m)\partial _t\chi _m\boldsymbol{u}_m\cdot\nabla \chi _m}\mathrm{d}x-\int\limits_{\Omega}{\rho _mm(\chi _m)F^{\prime\prime}(\chi _m)\partial _t\chi _m\boldsymbol{u}_m\cdot\nabla \chi _m}\mathrm{d}x\\
    &-\int\limits_{\Omega}{\rho _mm(\chi _m)F^{\prime}(\chi _m)\boldsymbol{u}_m\cdot\nabla \partial _t\chi _m}\mathrm{d}x-\int\limits_{\Omega}{\rho _mm^{\prime}(\chi _m)\left| \partial _t\chi _m \right|^2F ^{\prime}(\chi _m)}\mathrm{d}x \nonumber\\
    &-\int\limits_{\Omega}{\rho _mm(\chi _m)F^{\prime\prime}(\chi _m)\left| \partial _t\chi _m \right|^2}\mathrm{d}x. \nonumber
\end{align}
Next, taking $\boldsymbol{w}=\partial _t\boldsymbol{u}_m$ in $(\ref{6.2})$, we have
\begin{align}      \label{6.23}   
\frac{\mathrm{d}}{\mathrm{d}t}& \int_{\Omega}{\frac{1}{2}\eta \left( \chi _m \right) \left| \mathbb{D} \boldsymbol{u}_m \right|^2}\mathrm{d}x
+\int_{\Omega}{\rho _m\left| \partial _t\boldsymbol{u}_m \right|^2}\mathrm{d}x \nonumber
=\int_{\Omega}{-\mathrm{div}\left( \nabla \chi _m \otimes \nabla \chi _m \right)\cdot \partial _t\boldsymbol{u}_m}\mathrm{d}x
\\
&+\int_{\Omega}{-\rho _m\left( \boldsymbol{u}_m\cdot \nabla \right) \boldsymbol{u}_m\cdot \partial _t\boldsymbol{u}_m}\mathrm{d}x
+\int_{\Omega}{\frac{1}{2}\eta ^{\prime}\left( \chi _m \right) \partial _t\chi _m\left| \mathbb{D} \boldsymbol{u}_m \right|^2}\mathrm{d}x.
\end{align}
Using integration by parts on $(\ref{6.9phi})$, and adding $(\ref{6.23})$, we obtain
\begin{align}      \label{6.9}   
	\frac{\mathrm{d}}{\mathrm{d}t} &\left\{\int_{\Omega}{\frac{1}{2}\rho _{m}^{2}\left| \partial _t\chi _m \right|^2}\mathrm{d}x+\int_{\Omega}{\frac{1}{2}\eta \left( \chi _m \right) \left| \mathbb{D} \boldsymbol{u}_m \right|^2}\mathrm{d}x\right\}+m_*\int_{\Omega}{\left| \nabla \partial _t\chi _m \right|^2}\mathrm{d}x+\rho _*\int_{\Omega}{\left| \partial _t\boldsymbol{u}_m \right|^2}\mathrm{d}x \nonumber
\\
\le & {\int_{\Omega}{-{\rho _{m}^{2}} \partial _t\chi _m \boldsymbol{u}_m \cdot {\nabla} \partial _t\chi _m}\mathrm{d}x}
+{\int_{\Omega}{-{\rho _{m}^{2}}\left( \boldsymbol{u}_m\cdot {\nabla} \partial _t\chi _m \right) \left( \boldsymbol{u}_m\cdot \nabla \chi _m \right)}\mathrm{d}x} \nonumber
\\
&+{\int_{\Omega}{-{\rho _{m}^{2}}\partial _{\mathrm{t}}\chi _{\mathrm{m}}\boldsymbol{u}_m\cdot {\nabla} \left( \boldsymbol{u}_m\cdot \nabla \chi _m \right)}\mathrm{d}x}
+{\int_{\Omega}{-\rho _{m}^{2}\partial _t\chi _m\partial _t\boldsymbol{u}_m\cdot \nabla \chi _m}\mathrm{d}x}\nonumber
    \\
&+{\int_{\Omega}{-\rho _{m}^{2}\partial _t\chi _m\boldsymbol{u}_m \cdot\nabla \partial _t\chi _m}\mathrm{d}x}
+{\int_{\Omega}{m^{\prime}(\chi _m)\left| \partial _t\chi _m \right|^2\Delta \chi _m}\mathrm{d}x}
\\
&+{\int_{\Omega}{-m^{\prime}(\chi _m)\partial _t\chi _m\nabla \partial _t\chi _m\cdot \nabla \chi _m}\mathrm{d}x}
+{\int_{\Omega}{-\rho _mm^{\prime}(\chi _m)F^{\prime}(\chi _m)\partial _t\chi _m\boldsymbol{u}_m\cdot\nabla \chi _m}\mathrm{d}x} \nonumber
\\
&+{\int_{\Omega}{-\rho _mm(\chi _m)F^{\prime\prime}(\chi _m)\partial _t\chi _m\boldsymbol{u}_m\cdot\nabla \chi _m}\mathrm{d}x}
+{\int_{\Omega}{-\rho _mm(\chi _m)F^{\prime}(\chi _m)\boldsymbol{u}_m\cdot\nabla \partial _t\chi _m}\mathrm{d}x} \nonumber
\\
&+{\int_{\Omega}{-\rho _mm^{\prime}(\chi _m)\left| \partial _t\chi _m \right|^2F^{\prime}(\chi _m)}\mathrm{d}x}
+{\int_{\Omega}{-\rho _mm(\chi _m)F^{\prime\prime}(\chi _m)\left| \partial _t\chi _m \right|^2}\mathrm{d}x} \nonumber
\\
& +{\int_{\Omega}{-\Delta\chi _m\nabla \chi _m \cdot \partial _t\boldsymbol{u}_m}\mathrm{d}x} +{\int_{\Omega}{-\rho _m\left( \boldsymbol{u}_m\cdot \nabla \right) \boldsymbol{u}_m\cdot \partial _t\boldsymbol{u}_m}\mathrm{d}x} \nonumber
\\
&+{\int_{\Omega}{\frac{1}{2}\eta ^{\prime}\left( \chi _m \right) \partial _t\chi _m\left| \mathbb{D} \boldsymbol{u}_m \right|^2}\mathrm{d}x}
=\sum_{i=1}^{15}I_i. \nonumber
\end{align}
Here we have used
$\operatorname{div}(\nabla\chi_m\otimes\nabla\chi_m)
=\Delta\chi_m\nabla\chi_m+\frac12\nabla|\nabla\chi_m|^2$,
the divergence-free condition and the zero boundary trace of
$\partial_t\boldsymbol{u}_m$.

\noindent
{\it\bfseries Remark on the structure of \eqref{6.9}.}
Unlike the Cahn-Hilliard system in \cite{G-T}, where the test function $w=\partial_t\mu_m$ is admissible and yields a high-order estimate involving $\|\nabla\mu_m\|_{L^2}$ as a key energy term, the second-order Allen-Cahn equation does not admit such a test function. To obtain a dissipative estimate for $\|\nabla\partial_t\chi_m\|_{L^2}$, we are forced to differentiate (\ref{6.3}) directly in time and multiply the result by $\partial_t\chi_m$. This generates $\partial_t\rho_m$ in multiple terms, which are replaced by $-{\rm div}(\rho_m\boldsymbol{u}_m)$ so that the spatial derivative can be transferred to $\mathbf{u}_m$ and $\partial_t\chi_m$ via integration by parts, producing numerous high-order nonlinear terms on the right-hand side of (\ref{6.9}). Owing to the weaker regularity of the second-order Allen-Cahn equation (compared to the fourth-order Cahn-Hilliard system in \cite{G-T}), several of these terms are borderline and must be handled carefully.

\vskip1.8mm
\noindent
{\it\bfseries Three-dimensional case ($d=3$).}\quad
{ We first estimate $I_1,\ldots,I_{12}$. Throughout these estimates,
the coefficients and polynomial factors are bounded by
\eqref{6.4}. We apply H\"older's inequality, the embedding
$H^1(\Omega)\hookrightarrow L^6(\Omega)$, and
Gagliardo--Nirenberg interpolation for the $L^3$- and $L^4$-norms.
The Neumann elliptic estimate controls the spatial derivatives of
$\chi_m$, while \eqref{6.5}--\eqref{6.6} control the lower-order
factors. Young's inequality then gives the displayed small
multiples of the dissipation terms. }
\begin{flalign*}      
\qquad\quad	\left| I_1 \right|&\le C\left\| \partial _t\chi _m \right\| _{L^3\left( \Omega \right)}\left\| \boldsymbol{u}_m\right\| _{L^6\left( \Omega \right)}\left\| \nabla \partial _t\chi _m \right\| _{L^2\left( \Omega \right)} \nonumber
\\
	&\le C\left( \left\| \nabla \partial _t\chi _m \right\| _{L^2\left( \Omega \right)}^{\frac{1}{2}}\left\| \partial _t\chi _m \right\| _{L^2\left( \Omega \right)}^{\frac{1}{2}}+\left\| \partial _t\chi _m \right\| _{L^2\left( \Omega \right)} \right) \left\| \nabla \boldsymbol{u}_m \right\| _{L^2\left( \Omega \right)}\left\| \nabla \partial _t\chi _m \right\| _{L^2\left( \Omega \right)}
\\
	&=C\left\| \nabla \partial _t\chi _m \right\| _{L^2\left( \Omega \right)}^{\frac{3}{2}}\left\| \partial _t\chi _m \right\| _{L^2\left( \Omega \right)}^{\frac{1}{2}}\left\| \nabla \boldsymbol{u}_m\right\| _{L^2\left( \Omega \right)}
\\
&\quad +C\left\| \partial _t\chi _m \right\| _{L^2\left( \Omega \right)}\left\| \nabla \boldsymbol{u}_m\right\| _{L^2\left( \Omega \right)}\left\| \nabla \partial _t\chi _m \right\| _{L^2\left( \Omega \right)} \nonumber
\\
	&\le \frac{m_*}{16}\left\| \nabla \partial _t\chi _m \right\| _{L^2\left( \Omega \right)}^{2}+C\left\| \partial _t\chi _m \right\| _{L^2\left( \Omega \right)}^{2}\left\| \nabla \boldsymbol{u}_m\right\| _{L^2\left( \Omega \right)}^{4}+C\left\| \partial _t\chi _m \right\| _{L^2\left( \Omega \right)}^{2} \left\| \nabla \boldsymbol{u}_m\right\| _{L^2\left( \Omega \right)}^{2},& \nonumber
\end{flalign*}
\begin{flalign*}    
\qquad\quad	\left| I_2 \right|&\le C\left\| \boldsymbol{u}_m \right\| _{L^6\left( \Omega \right)}^{2}\left\| \nabla \partial _t\chi _m \right\| _{L^2\left( \Omega \right)}\left\| \nabla \chi _m \right\| _{L^6\left( \Omega \right)} \nonumber\\
	&\le C\left\| \nabla \boldsymbol{u}_m \right\| _{L^2\left( \Omega \right)}^{2}\left\| \nabla \partial _t\chi _m \right\| _{L^2\left( \Omega \right)}\left\| \Delta \chi _m \right\| _{L^2\left( \Omega \right)} \\
	&\le \frac{m_*}{16}\left\| \nabla \partial _t\chi _m \right\| _{L^2\left( \Omega \right)}^{2}+C\left\| \Delta \chi _m \right\| _{L^2\left( \Omega \right)}^{2}\left\| \nabla \boldsymbol{u}_m \right\| _{L^2\left( \Omega \right)}^{4},& \nonumber
\end{flalign*}
\begin{flalign*}     
\qquad\quad	\left| I_3 \right|&\le C\left\| \partial _t\chi _m \right\| _{L^6\left( \Omega \right)}\left\| \boldsymbol{u}_m \right\| _{L^6\left( \Omega \right)}\left( \left\| \nabla \boldsymbol{u}_m \right\| _{L^2\left( \Omega \right)}\left\| \nabla \chi _m \right\| _{L^6\left( \Omega \right)}+\left\| \boldsymbol{u}_m \right\| _{L^6\left( \Omega \right)}\left\| \Delta \chi _m \right\| _{L^2\left( \Omega \right)} \right) \nonumber\\
	&\le C\left( \left\| \nabla \partial _t\chi _m \right\| _{L^2\left( \Omega \right)}+\left\| \partial _t\chi _m \right\| _{L^2\left( \Omega \right)} \right) \left\| \nabla \boldsymbol{u}_m \right\| _{L^2\left( \Omega \right)}^{2}\left\| \Delta \chi _m \right\| _{L^2\left( \Omega \right)}\\
	&\le \frac{m_*}{16}\left\| \nabla \partial _t\chi _m \right\| _{L^2\left( \Omega \right)}^{2}+C\left\| \partial _t\chi _m \right\| _{L^2\left( \Omega \right)}^{2}+C\left\| \Delta \chi _m \right\| _{L^2\left( \Omega \right)}^{2}\left\| \nabla \boldsymbol{u}_m \right\| _{L^2\left( \Omega \right)}^{4},& \nonumber
\end{flalign*}
\begin{flalign*}   
\qquad\quad	\left| I_4 \right|&\le C\left\| \partial _t\chi _m \right\| _{L^3\left( \Omega \right)}\left\| \partial _t\boldsymbol{u}_m \right\| _{L^2\left( \Omega \right)}\left\| \nabla \chi _m \right\| _{L^6\left( \Omega \right)} \nonumber\\
	&\le C\left( \left\| \nabla \partial _t\chi _m \right\| _{L^2\left( \Omega \right)}^{\frac{1}{2}}\left\| \partial _t\chi _m \right\| _{L^2\left( \Omega \right)}^{\frac{1}{2}}+\left\| \partial _t\chi _m \right\| _{L^2\left( \Omega \right)} \right) \left\| \partial _t\boldsymbol{u}_m \right\| _{L^2\left( \Omega \right)}\left\| \Delta \chi _m \right\| _{L^2\left( \Omega \right)}\\
	&\le \frac{m_*}{16}\left\| \nabla \partial _t\chi _m \right\| _{L^2\left( \Omega \right)}^{2}+\frac{\rho _*}{8}\left\| \partial _t\boldsymbol{u}_m \right\| _{L^2\left( \Omega \right)}^{2}
\\
&\quad +C\left( \left\| \Delta \chi _m \right\| _{L^2\left( \Omega \right)}^{4}+\left\| \Delta \chi _m \right\| _{L^2\left( \Omega \right)}^{2} \right) \left\| \partial _t\chi _m \right\| _{L^2\left( \Omega \right)}^{2},& \nonumber
\end{flalign*}
\begin{flalign*}     
\qquad\quad	\left| I_5 \right|&\le C\left\| \partial _t\chi _m \right\| _{L^3\left( \Omega \right)}\left\| \boldsymbol{u}_m \right\| _{L^6\left( \Omega \right)}\left\| \nabla \partial _t\chi _m \right\| _{L^2\left( \Omega \right)} \nonumber\\
	&\le C\left( \left\| \nabla \partial _t\chi _m \right\| _{L^2\left( \Omega \right)}^{\frac{1}{2}}\left\| \partial _t\chi _m \right\| _{L^2\left( \Omega \right)}^{\frac{1}{2}}+\left\| \partial _t\chi _m \right\| _{L^2\left( \Omega \right)} \right) \left\| \nabla \boldsymbol{u}_m \right\| _{L^2\left( \Omega \right)}\left\| \nabla \partial _t\chi _m \right\| _{L^2\left( \Omega \right)}\\
	&\le \frac{m_*}{16}\left\| \nabla \partial _t\chi _m \right\| _{L^2\left( \Omega \right)}^{2}+C\left\| \nabla \boldsymbol{u}_m \right\| _{L^2\left( \Omega \right)}^{4}\left\| \partial _t\chi _m \right\| _{L^2\left( \Omega \right)}^{2}+C\left\| \nabla \boldsymbol{u}_m \right\| _{L^2\left( \Omega \right)}^{2}\left\| \partial _t\chi _m \right\| _{L^2\left( \Omega \right)}^{2},& \nonumber
\end{flalign*}
\begin{flalign*}    
\qquad\quad	\left| I_6 \right|&\le C\left\| \partial _t\chi _m \right\| _{L^4\left( \Omega \right)}^{2}\left\| \Delta \chi _m \right\| _{L^2\left( \Omega \right)}\nonumber\\
	&\le C\left( \left\| \nabla \partial _t\chi _m \right\| _{L^2\left( \Omega \right)}^{\frac{3}{2}}\left\| \partial _t\chi _m \right\| _{L^2\left( \Omega \right)}^{\frac{1}{2}}+\left\| \partial _t\chi _m \right\| _{L^2\left( \Omega \right)}^{2} \right) \left\| \Delta \chi _m \right\| _{L^2\left( \Omega \right)}\\
	&\le \frac{m_*}{16}\left\| \nabla \partial _t\chi _m \right\| _{L^2\left( \Omega \right)}^{2}+C\left\| \Delta \chi _m \right\| _{L^2\left( \Omega \right)}^{4}\left\| \partial _t\chi _m \right\| _{L^2\left( \Omega \right)}^{2}+C\left\| \Delta \chi _m \right\| _{L^2\left( \Omega \right)}\left\| \partial _t\chi _m \right\| _{L^2\left( \Omega \right)}^{2},& \nonumber
\end{flalign*}
\begin{flalign*}      
\qquad\quad	\left| I_7 \right|&\le C\left\| \partial _t\chi _m \right\| _{L^3\left( \Omega \right)}\left\| \nabla \partial _t\chi _m \right\| _{L^2\left( \Omega \right)}\left\| \nabla \chi _m \right\| _{L^6\left( \Omega \right)}\nonumber\\
	&\le C\left( \left\| \nabla \partial _t\chi _m \right\| _{L^2\left( \Omega \right)}^{\frac{1}{2}}\left\| \partial _t\chi _m \right\| _{L^2\left( \Omega \right)}^{\frac{1}{2}}+\left\| \partial _t\chi _m \right\| _{L^2\left( \Omega \right)} \right)\left\| \nabla \partial _t\chi _m \right\| _{L^2\left( \Omega \right)}\left\| \Delta \chi _m \right\| _{L^2\left( \Omega \right)}
\\
	&\le \frac{m_*}{16}\left\| \nabla \partial _t\chi _m \right\| _{L^2\left( \Omega \right)}^{2}+C\left\| \Delta \chi _m \right\| _{L^2\left( \Omega \right)}^{4}\left\| \partial _t\chi _m \right\| _{L^2\left( \Omega \right)}^{2}+C\left\| \Delta \chi _m \right\| _{L^2\left( \Omega \right)}^{2}\left\| \partial _t\chi _m \right\| _{L^2\left( \Omega \right)}^{2},& \nonumber
\end{flalign*}
\begin{flalign*}      
\qquad\quad	\left| I_8+I_9 \right|
&\le C\left\| \partial _t\chi _m \right\| _{L^2\left( \Omega \right)}\left\| \boldsymbol{u}_m \right\| _{L^6\left( \Omega \right)}\left\| \nabla \chi _m \right\| _{L^3\left( \Omega \right)}\nonumber
\\
&\le C\left( E_0 \right) \left\| \partial _t\chi _m \right\| _{L^2\left( \Omega \right)}
\left\|\nabla \boldsymbol{u}_m \right\| _{L^2\left( \Omega \right)}
\left\| \Delta \chi _m \right\| _{L^2\left( \Omega \right)}
\\
&\le C \left\| \partial _t\chi _m \right\| _{L^2\left( \Omega \right)}^{2}
+C\left\| \Delta \chi _m \right\| _{L^2\left( \Omega \right)}^{2}\left\| \nabla \boldsymbol{u}_m \right\| _{L^2\left( \Omega \right)}^{2},& \nonumber
\end{flalign*}
\begin{flalign*}      
\qquad\quad	\left| I_{10}+ I_{11}+I_{12} \right|
&\le C\left\| \boldsymbol{u}_m \right\| _{L^2\left( \Omega \right)}\left\| \nabla \partial _t\chi _m \right\| _{L^2\left( \Omega \right)}
+C\left\| \partial _t\chi _m \right\| _{L^2\left( \Omega \right)}^{2}\nonumber
\\
	&\le C\left( E_0 \right)\left\| \nabla \partial _t\chi _m \right\| _{L^2\left( \Omega \right)}
+C\left\| \partial _t\chi _m \right\| _{L^2\left( \Omega \right)}^{2}
\\
	&\le \frac{m_*}{16}\left\| \nabla \partial _t\chi _m \right\| _{L^2\left( \Omega \right)}^{2}
+C\left\| \partial _t\chi _m \right\| _{L^2\left( \Omega \right)}^{2}
+C\left( E_0 \right).& \nonumber
\end{flalign*}
Noticing that
\begin{align}  \label{I-13}
\left|  I_{13}  \right|\le \frac{\rho _*}{4}\left\|\partial _t \boldsymbol{u}_m \right\| _{L^2\left( \Omega \right)}^{2}
+C \left\|{\Delta \chi _m \nabla \chi _m }\right\| _{L^2\left( \Omega \right)}^{2},
\end{align}
we need to estimate $\left\|{\Delta \chi _m \nabla \chi _m }\right\| _{L^2\left( \Omega \right)}$
\footnote{
From the basic energy estimate, the available regularity is only $\chi\in L^2(0,T;H^2(\Omega))$. This is insufficient to control the product $\|{\Delta \chi _m \nabla \chi _m }\|_{L^2}$ by separating the two factors through standard Sobolev embeddings. We therefore treat the product $\Delta\chi\nabla\chi$ as a whole.
}first.
Multiplying the equation $(\ref{6.3})$ by $ \nabla \chi_m$ yields
\begin{align}\label{phi-phi}
m\left( \chi _m \right) {\Delta \chi _m \nabla \chi _m}=\rho _{m}^{2}\partial _t\chi _m\nabla \chi _m+\rho _{m}^{2}\left(\boldsymbol{u}_m \cdot \nabla \chi _m \right)\nabla \chi _m+\rho _mm\left( \chi _m \right) F^{\prime}\left( \chi _m \right) \nabla \chi _m.
\end{align}
{ Applying H\"older's inequality to \eqref{phi-phi}, we estimate
$\partial_t\chi_m\nabla\chi_m$ in $L^2(\Omega)$ using the
$L^3$- and $L^6$-norms, and
$(\boldsymbol{u}_m\cdot\nabla\chi_m)\nabla\chi_m$ using three
$L^6$-norms. Interpolation, the Neumann elliptic estimate, and
\eqref{6.4}--\eqref{6.6}, followed by Young's inequality, yield }
\begin{align}      \label{6.24}   
&m_*\left\| {\Delta \chi _m \nabla \chi _m} \right\| _{L^2\left( \Omega \right)}^{2}  \nonumber
\\
&\le C\left\| \partial _t\chi _m\nabla \chi _m \right\| _{L^2\left( \Omega \right)}^{2}+C\left\| \left(\boldsymbol{u}_m \cdot \nabla \chi _m \right)\nabla \chi _m \right\| _{L^2\left( \Omega \right)}^{2}+C\left\| F^{\prime}\left( \chi _m \right) \nabla \chi _m \right\| _{L^2\left( \Omega \right)}^{2}\nonumber
\\
	&\le C\left\| \partial _t\chi _m \right\| _{L^3\left( \Omega \right)}^{2}\left\| \nabla \chi _m \right\| _{L^6\left( \Omega \right)}^{2}+C\left\| \boldsymbol{u}_m \right\| _{L^6\left( \Omega \right)}^{2}\left\| \nabla \chi _m \right\| _{L^6\left( \Omega \right)}^{4}
+C\left\| \nabla \chi _m \right\| _{L^2\left( \Omega \right)}^{2}\nonumber
    \\
	&\le C\left( E_0 \right) \left( \left\| \nabla \partial _t\chi _m \right\| _{L^2\left( \Omega \right)}\left\| \partial _t\chi _m \right\| _{L^2\left( \Omega \right)}+\left\| \partial _t\chi _m \right\| _{L^2\left( \Omega \right)}^{2} \right) \left\| \Delta \chi _m \right\| _{L^2\left( \Omega \right)}^{2}\nonumber
\\
	&\quad +C\left( E_0 \right) \left\| \nabla \boldsymbol{u}_m \right\| _{L^2\left( \Omega \right)}^{2}\left\| \Delta \chi _m \right\| _{L^2\left( \Omega \right)}^{4}+C\left( E_0 \right)\nonumber
\\
	&\le \varepsilon \left\| \nabla \partial _t\chi _m \right\| _{L^2\left( \Omega \right)}^{2}+C\left( E_0 \right) \left( \left\| \Delta \chi _m \right\| _{L^2\left( \Omega \right)}^{4}+\left\| \Delta \chi _m \right\| _{L^2\left( \Omega \right)}^{2} \right) \left\| \partial _t\chi _m \right\| _{L^2\left( \Omega \right)}^{2}\nonumber
\\
    &\quad +C\left( E_0 \right) \left\| \Delta \chi _m \right\| _{L^2\left( \Omega \right)}^{4}\left\| \nabla \boldsymbol{u}_m \right\| _{L^2\left( \Omega \right)}^{2}+C\left( E_0 \right).
\end{align}
Putting (\ref{6.24}) into (\ref{I-13}) and taking $\varepsilon=\frac{m _*}{8}$, we have
\begin{flalign}      \label{6.25}   
\qquad\qquad \left|  I_{13}  \right|&\le \frac{\rho _*}{4}\left\|\partial _t \boldsymbol{u}_m \right\| _{L^2\left( \Omega \right)}^{2}+\frac{m _*}{8} \left\| \nabla \partial _t\chi _m \right\| _{L^2\left( \Omega \right)}^{2}+C\left( E_0 \right) \left\| \Delta \chi _m \right\| _{L^2\left( \Omega \right)}^{4}\left\| \nabla \boldsymbol{u}_m \right\| _{L^2\left( \Omega \right)}^{2} \nonumber
\\
	&\quad +C\left( E_0 \right) \left( \left\| \Delta \chi _m \right\| _{L^2\left( \Omega \right)}^{4}+\left\| \Delta \chi _m \right\| _{L^2\left( \Omega \right)}^{2} \right) \left\| \partial _t\chi _m \right\| _{L^2\left( \Omega \right)}^{2}+C\left( E_0 \right).&
\end{flalign}
{ For $I_{14}$, we use Agmon's inequality
\[
\|\boldsymbol{u}_m\|_{L^\infty(\Omega)}^2
\le C\|\boldsymbol{u}_m\|_{H^1(\Omega)}
      \|\boldsymbol{u}_m\|_{H^2(\Omega)},
\]
together with elliptic regularity and Young's inequality.
For $I_{15}$, we interpolate
$\|\nabla\boldsymbol{u}_m\|_{L^4(\Omega)}$ between its
$L^2(\Omega)$- and $H^1(\Omega)$-norms and then apply Young's
inequality with conjugate exponents $\frac43$ and $4$.
Thus, }
\begin{flalign}     \label{6.26}   
\qquad\qquad \left|  I_{14}  \right|&\le \rho^*
\left\| \boldsymbol{u}_m \right\| _{L^{\infty}\left( \Omega \right)}
\left\| \nabla \boldsymbol{u}_m \right\| _{L^2\left( \Omega \right)}
\left\| \partial _t\boldsymbol{u}_m \right\| _{L^2\left( \Omega \right)}\nonumber
\\
	&\le \frac{\rho _*}{4}\left\| \partial _t\boldsymbol{u}_m \right\| _{L^2\left( \Omega \right)}^{2}
+C\left\| \boldsymbol{u}_m \right\| _{L^{\infty}\left( \Omega \right)}^{2}\left\| \nabla \boldsymbol{u}_m \right\| _{L^2\left( \Omega \right)}^{2}  \nonumber
\\
	&\le \frac{\rho _*}{4}\left\| \partial _t\boldsymbol{u}_m \right\| _{L^2\left( \Omega \right)}^{2}
+C\left\| \boldsymbol{u}_m \right\| _{H^1\left( \Omega \right)}\left\| \boldsymbol{u}_m \right\| _{H^2\left( \Omega \right)}
\left\| \nabla \boldsymbol{u}_m \right\| _{L^2\left( \Omega \right)}^{2}
\\
	&\le \frac{\rho _*}{4}\left\| \partial _t\boldsymbol{u}_m \right\| _{L^2\left( \Omega \right)}^{2}+C\left( E_0 \right) \left\| \Delta \boldsymbol{u}_m \right\| _{L^2\left( \Omega \right)}\left\| \nabla \boldsymbol{u}_m \right\| _{L^2\left( \Omega \right)}^{3}+C\left( E_0 \right) \left\| \nabla \boldsymbol{u}_m \right\| _{L^2\left( \Omega \right)}^{4} \nonumber
\\
	&\le \frac{\rho _*}{4}\left\| \partial _t\boldsymbol{u}_m \right\| _{L^2\left( \Omega \right)}^{2}
+\frac{\varpi}{2}\left\| {\mathbf{A}\boldsymbol{u}_m }\right\| _{L^2\left( \Omega \right)}^{2}
+C\left( E_0 \right) \left\| \nabla \boldsymbol{u}_m \right\| _{L^2\left( \Omega \right)}^{6}
+C\left( E_0 \right),& \nonumber
\\
\label{6.27}   
\qquad\qquad \left|  I_{15}  \right|&\le C\left\| \partial _t\chi _m \right\| _{L^2\left( \Omega \right)}\left\| \nabla \boldsymbol{u}_m \right\| _{L^4\left( \Omega \right)}^{2} \nonumber
\\
	&\le C\left\| \partial _t\chi _m \right\| _{L^2\left( \Omega \right)}
\left\| \nabla \boldsymbol{u}_m \right\| _{H^1\left( \Omega \right)}^{\frac{3}{2}}
\left\| \nabla \boldsymbol{u}_m \right\| _{L^2\left( \Omega \right)}^{\frac{1}{2}}
\\
	&\le C\left\| \partial _t\chi _m \right\| _{L^2\left( \Omega \right)}^{4}
\left\| \nabla \boldsymbol{u}_m \right\| _{L^2\left( \Omega \right)}^{2}+\frac{\varpi}{2}\left\| {\mathbf{A}\boldsymbol{u}_m }\right\| _{L^2\left( \Omega \right)}^{2},& \nonumber
\end{flalign}
where $\varpi$ is a positive (small) constant which will be determined later, {and $\mathbf{A}$ is the Stokes operator.
Here Stokes regularity is used before the final application of
Young's inequality, so that the small higher-order velocity
terms are expressed in terms of
$\|\mathbf{A}\boldsymbol{u}_m\|_{L^2(\Omega)}^2$.}
Putting the estimates on $I_i~ (i=1, \cdots, {15} )$ into $(\ref{6.9})$, we arrive at
\begin{align}      \label{6.28}   
	&\frac{\mathrm{d}}{\mathrm{d}t} \left\{\int_{\Omega}{\frac{1}{2}\rho _{m}^{2}\left| \partial _t\chi _m \right|^2}\mathrm{d}x+\int_{\Omega}{\frac{1}{2}\eta \left( \chi _m \right) \left| \mathbb{D} \boldsymbol{u}_m \right|^2}\mathrm{d}x\right\}+\frac{m _*}{2}\left\| \nabla \partial _t\chi _m  \right\| _{L^2\left( \Omega \right)}^{2}+\frac{\rho _*}{2}\left\| \partial _t\boldsymbol{u}_m \right\| _{L^2\left( \Omega \right)}^{2} \nonumber
\\
	&\le \varpi \left\| \mathbf{A}\boldsymbol{u}_m \right\| _{L^2\left( \Omega \right)}^{2}+\frac{\rho _*}{8}\left\| \partial _t\boldsymbol{u}_m \right\| _{L^2\left( \Omega \right)}^{2}+\frac{m_*}{8}\left\| \nabla \partial _t\chi _m \right\| _{L^2\left( \Omega \right)}^{2}+C\left( E_0 \right) \\
    &\quad +C\left( E_0 \right) \left( \left\| \Delta \chi _m \right\| _{L^2\left( \Omega \right)}^{4}
+\left\| \partial_t \chi _m \right\| _{L^2\left( \Omega \right)}^{4}
+\left\| \nabla \boldsymbol{u}_m \right\| _{L^2\left( \Omega \right)}^{4}+1 \right)
\left( \left\| \partial _t\chi _m \right\| _{L^2\left( \Omega \right)}^{2}
+\left\| \nabla \boldsymbol{u}_m \right\| _{L^2\left( \Omega \right)}^{2} \right). \nonumber
\end{align}
Here the constant $C$ depends on $\varpi$. Next, taking $\boldsymbol{w}=\mathbf{A} \boldsymbol{u}_{m}$ in $(\ref{6.2})$, we obtain
\begin{align}      \label{6.29}   
	\int_{\Omega}{-\mathrm{div}\left( \eta \left( \chi _m \right) \mathbb{D} \boldsymbol{u}_m \right) \cdot \mathbf{A}\boldsymbol{u}_m}\mathrm{d}x=&-\left( \rho _m\partial _t\boldsymbol{u}_m,\mathbf{A}\boldsymbol{u}_m \right) -\left( \rho _m\left( \boldsymbol{u}_m\cdot \nabla \right) \boldsymbol{u}_m,\mathbf{A}\boldsymbol{u}_m \right) \nonumber\\
	&-\left( \mathrm{div}\left( \nabla \chi _m\otimes \nabla \chi _m \right) ,\mathbf{A}\boldsymbol{u}_m \right).
\end{align}
{From \cite[Lemma 3.1, Remark 3.2]{G3-P}}, there exists $\pi_{m} \in L^{2}(0, T ; H^{1}(\Omega))$ such that $-\Delta \boldsymbol{u}_{m}+\nabla \pi_{m}=\mathbf{A} \boldsymbol{u}_{m}$ almost everywhere in $\Omega \times(0, T)$ and such that
\begin{align}      \label{6.30}   
{\left\|\pi_{m}\right\|_{L^{3}(\Omega)}} \leq C\left\|\nabla \boldsymbol{u}_{m}\right\|_{L^{2}(\Omega)}^{\frac{1}{2}}\left\|\mathbf{A} \boldsymbol{u}_{m}\right\|_{L^{2}(\Omega)}^{\frac{1}{2}}, \quad\left\|\pi_{m}\right\|_{H^{1}(\Omega)} \leq C\left\|\mathbf{A} \boldsymbol{u}_{m}\right\|_{L^{2}(\Omega)},
\end{align}
where $C$ is independent of $m$. Then, we rewrite $(\ref{6.29})$ as follows
\begin{align}      \label{6.31}   
	\int_{\Omega}{\frac{\eta \left( \chi _m \right)}{2}\left| \mathbf{A}\boldsymbol{u}_m \right|^2}\mathrm{d}x
&=\left( -\rho _m\partial _t\boldsymbol{u}_m,\mathbf{A}\boldsymbol{u}_m \right)
+\left( -\rho _m\left( \boldsymbol{u}_m\cdot \nabla \right) \boldsymbol{u}_m,\mathbf{A}\boldsymbol{u}_m \right) \nonumber
\\
&\quad +\left( -\Delta \chi _m\nabla \chi _m ,\mathbf{A}\boldsymbol{u}_m \right)
+\left( \eta ^{\prime}\left( \chi _m \right) \mathbb{D} u_m\nabla \chi _m,\mathbf{A}\boldsymbol{u}_m \right)
\\
&\quad +( -\frac{1}{2}\eta ^{\prime}\left( \chi _m \right) \pi _m\nabla \chi _m,\mathbf{A}\boldsymbol{u}_m )
=\sum_{i=1}^5K_i. \nonumber
\end{align}
{ For $K_1$, we apply Young's inequality directly. For $K_2$ and
$K_4$, we use the Sobolev embedding and
\[
\|\nabla\boldsymbol{u}_m\|_{L^3(\Omega)}
\le C\|\nabla\boldsymbol{u}_m\|_{L^2(\Omega)}^{1/2}
      \|\mathbf{A}\boldsymbol{u}_m\|_{L^2(\Omega)}^{1/2},
\]
which follows from interpolation and Stokes regularity.
The term $K_3$ is controlled by \eqref{6.24}, whereas $K_5$
uses the $L^3$-estimate for $\pi_m$ in \eqref{6.30}.
Young's inequality then yields }
\begin{flalign}      \label{6.32}   
\qquad\qquad\quad \left| K_1 \right|&\le \rho ^*\left\| \partial _t\boldsymbol{u}_m \right\| _{L^2\left( \Omega \right)}\left\| \mathbf{A}\boldsymbol{u}_m \right\| _{L^2\left( \Omega \right)} \le \frac{\eta _*}{20}\left\| \mathbf{A}\boldsymbol{u}_m \right\| _{L^2\left( \Omega \right)}^{2}+\frac{5\left( \rho ^* \right) ^2}{\eta _*}\left\| \partial _t\boldsymbol{u}_m \right\| _{L^2\left( \Omega \right)}^{2},&
\\     \label{6.33}   
\qquad\qquad\quad \left| K_2 \right| &\le \rho ^*\left\| \boldsymbol{u}_m \right\| _{L^6\left( \Omega \right)}\left\| \nabla \boldsymbol{u}_m \right\| _{L^3\left( \Omega \right)}\left\| \mathbf{A}\boldsymbol{u}_m \right\| _{L^2\left( \Omega \right)} \nonumber
\\
    &\le C\left( E_0 \right) \left\| \nabla \boldsymbol{u}_m \right\| _{L^2\left( \Omega \right)}^{\frac{3}{2}}\left\| \mathbf{A}\boldsymbol{u}_m \right\| _{L^2\left( \Omega \right)}^{\frac{3}{2}}
    \\
    &\le \frac{\eta _*}{20}\left\| \mathbf{A}\boldsymbol{u}_m \right\| _{L^2\left( \Omega \right)}^{2}+C\left( E_0 \right) \left\| \nabla \boldsymbol{u}_m \right\| _{L^2\left( \Omega \right)}^{6},& \nonumber
\\    \label{6.34}   
\qquad\qquad\quad \left| K_3 \right|&\le \frac{\eta _*}{20}\left\| \mathbf{A}\boldsymbol{u}_m \right\| _{L^2\left( \Omega \right)}^{2}+C \left\|{\Delta \chi _m \nabla \chi _m }\right\| _{L^2\left( \Omega \right)}^{2}
\\
    &\le \frac{\eta _*}{20}\left\| \mathbf{A}\boldsymbol{u}_m \right\| _{L^2\left( \Omega \right)}^{2}+\varepsilon \left\| \nabla \partial _t\chi _m \right\| _{L^2\left( \Omega \right)}^{2}+C\left( E_0 \right) \left\| \Delta \chi _m \right\| _{L^2\left( \Omega \right)}^{4}\left\| \nabla \boldsymbol{u}_m \right\| _{L^2\left( \Omega \right)}^{2}\nonumber
    \\
	&\quad +C\left( E_0 \right) \left( \left\| \Delta \chi _m \right\| _{L^2\left( \Omega \right)}^{4}+\left\| \Delta \chi _m \right\| _{L^2\left( \Omega \right)}^{2} \right) \left\| \partial _t\chi _m \right\| _{L^2\left( \Omega \right)}^{2} +C\left( E_0 \right) ,& \nonumber
\\     \label{6.35}   
\qquad\qquad\quad	\left| K_4 \right|&\le C\left\| \nabla \boldsymbol{u}_m \right\| _{L^3\left( \Omega \right)}\left\| \nabla \chi _m \right\| _{L^6\left( \Omega \right)}\left\| \mathbf{A}\boldsymbol{u}_m \right\| _{L^2\left( \Omega \right)} \nonumber
\\
	&\le C\left\| \nabla \boldsymbol{u}_m \right\| _{L^2\left( \Omega \right)}^{\frac{1}{2}}\left\| \mathbf{A}\boldsymbol{u}_m \right\| _{L^2\left( \Omega \right)}^{\frac{3}{2}}\left\| \Delta \chi _m \right\| _{L^2\left( \Omega \right)}\\
	&\le \frac{\eta _*}{20}\left\| \mathbf{A}\boldsymbol{u}_m \right\| _{L^2\left( \Omega \right)}^{2}+C\left\| \Delta \chi _m \right\| _{L^2\left( \Omega \right)}^{4}\left\| \nabla \boldsymbol{u}_m \right\| _{L^2\left( \Omega \right)}^{2},& \nonumber
\\    \label{6.36}   
\qquad\qquad\quad	
|K_{5}|&\leq C\|\pi_{m}\|_{L^{3}(\Omega)}\|\nabla\chi_{m}\|_{L^{6}(\Omega)}
\|\mathbf{A}\mathbf{u}_{m}\|_{L^{2}(\Omega)} \nonumber
\\
&\leq C{\|\nabla\mathbf{u}_{m}\|_{L^{2}(\Omega)}^{\frac{1}{2}}
\|\mathbf{A}\mathbf{u}_{m}\|_{L^{2}(\Omega)}^{\frac{3}{2}}}\|\Delta\chi_{m}\|_{L^{2}(\Omega)}
\\
&\leq\frac{\eta_{*}}{20}\|\mathbf{A}\mathbf{u}_{m}\|_{L^{2}(\Omega)}^{2}
+C\|\nabla\mathbf{u}_{m}\|_{L^{2}(\Omega)}^{2}{\|\Delta\chi_{m}\|_{L^{2}(\Omega)}^{4}}. \nonumber
\end{flalign}
Combining the above estimates $(\ref{6.32})$-$(\ref{6.36})$ with $(\ref{6.31})$, we arrive at
\begin{align}      \label{6.37}   
	\frac{\eta _*}{4}\left\| \mathbf{A}\boldsymbol{u}_m \right\| _{L^2\left( \Omega \right)}^{2}\le &\frac{5\left( \rho ^* \right) ^2}{\eta _*}\left\| \partial _t\boldsymbol{u}_m \right\| _{L^2\left( \Omega \right)}^{2}+\varepsilon \left\| \nabla \partial _t\chi _m \right\| _{L^2\left( \Omega \right)}^{2}+C\left( E_0 \right) \left\| \nabla \boldsymbol{u}_m \right\| _{L^2\left( \Omega \right)}^{6} \nonumber
\\
	&+C\left( E_0 \right) \left( \left\| \Delta \chi _m \right\| _{L^2\left( \Omega \right)}^{4}+\left\| \Delta \chi _m \right\| _{L^2\left( \Omega \right)}^{2} \right) \left\| \partial _t\chi _m \right\| _{L^2\left( \Omega \right)}^{2}
\\
    &+C\left( E_0 \right) \left\| \Delta \chi _m \right\| _{L^2\left( \Omega \right)}^{4} \left\| \nabla \boldsymbol{u}_m \right\| _{L^2\left( \Omega \right)}^{2}+C\left( E_0 \right). \nonumber
\end{align}
Choose $\delta =\frac{\rho _*\eta _*}{40\left( \rho ^* \right) ^2}$ and set
$\varpi =\frac{\rho _*\eta _{*}^{2}}{320\left( \rho ^* \right) ^2}$, $\varepsilon =\frac{m_*}{8\delta}$. Multiplying $(\ref{6.37})$ by $\delta$ yields
\begin{align}      \label{6.39}   
	2\varpi \left\| \mathbf{A}\boldsymbol{u}_m \right\| _{L^2\left( \Omega \right)}^{2}\le &\frac{\rho _*}{8}\left\| \partial _t\boldsymbol{u}_m \right\| _{L^2\left( \Omega \right)}^{2}+\frac{m_*}{8}\left\| \nabla \partial _t\chi _m \right\| _{L^2\left( \Omega \right)}^{2}
+C\left( E_0 \right)
\\
	&+C\left( E_0 \right) \left( {\left\| \Delta \chi _m \right\| _{L^2\left( \Omega \right)}^{4}}
+\left\|\nabla \boldsymbol{u}_m \right\| _{L^2\left( \Omega \right)}^{4} +1\right)
\left( \left\| \partial _t\chi _m \right\| _{L^2\left( \Omega \right)}^{2}
+\left\| \nabla \boldsymbol{u}_m \right\| _{L^2\left( \Omega \right)}^{2} \right).  \nonumber
\end{align}
Moreover, from the equation $(\ref{6.3})$ we see that
\begin{align*}
\|\Delta\chi_m\|_{L^2(\Omega)}^2
\le &C\|\partial_t\chi_m\|_{L^2(\Omega)}^2
+C\|\boldsymbol{u}_m\|_{L^6(\Omega)}^2\|\nabla\chi_m\|_{L^3(\Omega)}^2+C
\\
\le&C\|\partial_t\chi_m\|_{L^2(\Omega)}^2
+C\|\nabla\boldsymbol{u}_m\|_{L^2(\Omega)}^2\|\nabla\chi_m\|_{L^2(\Omega)}\|\Delta\chi_m\|_{L^2(\Omega)}+C
\\
\le&\frac12\|\Delta\chi_m\|_{L^2(\Omega)}^2
+C\|\partial_t\chi_m\|_{L^2(\Omega)}^2
+C(E_0)\|\nabla\boldsymbol{u}_m\|_{L^2(\Omega)}^4+C,
\end{align*}
i.e.
\begin{align}\label{phi-2}
\|\Delta\chi_m\|_{L^2(\Omega)}^2
\le C\|\partial_t\chi_m\|_{L^2(\Omega)}^2
+C(E_0)\|\nabla\boldsymbol{u}_m\|_{L^2(\Omega)}^4+C.
\end{align}
Set
\begin{align*}
&H_m=H\left( \rho _m,\boldsymbol{u}_m,\chi _m \right) =\int_{\Omega}{\frac{1}{2}\rho _{m}^{2}\left| \partial _t\chi _m \right|^2}\mathrm{d}x+\int_{\Omega}{\frac{1}{2}\eta \left( \chi _m \right) \left| \mathbb{D} \boldsymbol{u}_m \right|^2}\mathrm{d}x,
\\
&F_m=F\left( \boldsymbol{u}_m,\chi _m \right) =\varpi \left\| \mathbf{A}\boldsymbol{u}_m \right\| _{L^2\left( \Omega \right)}^{2}+\frac{\rho _*}{4}\left\| \partial _t\boldsymbol{u}_m \right\| _{L^2\left( \Omega \right)}^{2}+\frac{m_*}{4}\left\| \nabla \partial _t\chi _m \right\| _{L^2\left( \Omega \right)}^{2}.
\end{align*}
Adding $(\ref{6.28})$ and $(\ref{6.39})$ together, and using (\ref{phi-2}) we eventually obtain
\begin{align}      \label{6.40}   
\frac{\mathrm{d}}{\mathrm{d}t}H_m+F_m\le C\left( E_0 \right) \left( C\left( E_0 \right) +H_m \right) ^{{5}},
\end{align}
where the positive constant $C\left(E_{0}\right)$ is independent of the parameter $m$.
\vskip1.8mm
\noindent
{\it\bfseries Two-dimensional case ($d=2$).}\quad
{ In two dimensions, we use Gagliardo--Nirenberg interpolation
for the $L^4$-, $L^6$-, and $L^8$-norms, together with
\eqref{6.4}--\eqref{6.6} and the Neumann elliptic estimate.
In particular, the interpolation exponents for
$\partial_t\chi_m$ in $L^4(\Omega)$ and $L^8(\Omega)$ are
$1/2$ and $3/4$, respectively. Young's inequality gives the
small dissipation terms in the following estimates. }
\begin{flalign*}      
\qquad\quad	\left| I_1 \right|&\le C\left\| \partial _t\chi _m \right\| _{L^4\left( \Omega \right)}\left\| \boldsymbol{u}_m \right\| _{L^4\left( \Omega \right)}\left\| \nabla \partial _t\chi _m \right\| _{L^2\left( \Omega \right)}\nonumber\\
	&\le C\left( E_0 \right) \left( \left\| \nabla \partial _t\chi _m \right\| _{L^2\left( \Omega \right)}^{\frac{1}{2}}\left\| \partial _t\chi _m \right\| _{L^2\left( \Omega \right)}^{\frac{1}{2}}+\left\| \partial _t\chi _m \right\| _{L^2\left( \Omega \right)} \right) \left\| \nabla \boldsymbol{u}_m \right\| _{L^2\left( \Omega \right)}^{\frac{1}{2}}\left\| \nabla \partial _t\chi _m \right\| _{L^2\left( \Omega \right)}\nonumber\\
	&=C\left( E_0 \right) \left\| \partial _t\chi _m \right\| _{L^2\left( \Omega \right)}^{\frac{1}{2}}\left\| \nabla \boldsymbol{u}_m \right\| _{L^2\left( \Omega \right)}^{\frac{1}{2}}\left\| \nabla \partial _t\chi _m \right\| _{L^2\left( \Omega \right)}^{\frac{3}{2}}\\
    &\quad +C\left( E_0 \right) \left\| \partial _t\chi _m \right\| _{L^2\left( \Omega \right)}\left\| \nabla \boldsymbol{u}_m \right\| _{L^2\left( \Omega \right)}^{\frac{1}{2}}\left\| \nabla \partial _t\chi _m \right\| _{L^2\left( \Omega \right)}\nonumber\\
	&\le \frac{m_*}{16}\left\| \nabla \partial _t\chi _m \right\| _{L^2\left( \Omega \right)}^{2}+C\left( E_0 \right) \left\| \nabla \boldsymbol{u}_m \right\| _{L^2\left( \Omega \right)}^{2}\left\| \partial _t\chi _m \right\| _{L^2\left( \Omega \right)}^{2}+C\left( E_0 \right) \left\| \partial _t\chi _m \right\| _{L^2\left( \Omega \right)}^{2},& \nonumber
\end{flalign*}
\begin{flalign*}
\qquad\quad	\left| I_2 \right|&\le C\left\| \boldsymbol{u}_m \right\| _{L^8\left( \Omega \right)}^{2}\left\| \nabla \partial _t\chi _m \right\| _{L^2\left( \Omega \right)}\left\| \nabla \chi _m \right\| _{L^4\left( \Omega \right)} \nonumber
\\
	&\le C\left( E_0 \right) \left\| \nabla \boldsymbol{u}_m \right\| _{L^2\left( \Omega \right)}^{\frac{3}{2}}\left\| \nabla \partial _t\chi _m \right\| _{L^2\left( \Omega \right)}\left\| \Delta \chi _m \right\| _{L^2\left( \Omega \right)}^{\frac{1}{2}}
\\
	&\le \frac{m_*}{16}\left\| \nabla \partial _t\chi _m \right\| _{L^2\left( \Omega \right)}^{2}+C\left( E_0 \right) \left\| \Delta \chi _m \right\| _{L^2\left( \Omega \right)}\left\| \nabla \boldsymbol{u}_m \right\| _{L^2\left( \Omega \right)}^{3} \nonumber
\\
&\le \frac{m_*}{16}\left\| \nabla \partial _t\chi _m \right\| _{L^2\left( \Omega \right)}^{2}
+C\left( E_0 \right) \left(\left\| \Delta \chi _m \right\| _{L^2\left( \Omega \right)}^2
+\|\nabla\boldsymbol{u}_m\|_{L^2(\Omega)}^2\right)
\left\| \nabla \boldsymbol{u}_m \right\| _{L^2\left( \Omega \right)}^2,& \nonumber
\end{flalign*}
\begin{flalign*}
\qquad\quad	\left| I_3 \right|
&\le C\left\| \partial _t\chi _m \right\| _{L^8\left( \Omega \right)}\left\| \boldsymbol{u}_m \right\| _{L^4\left( \Omega \right)}\left( \left\| \nabla \boldsymbol{u}_m \right\| _{L^2\left( \Omega \right)}\left\| \nabla \chi _m \right\| _{L^8\left( \Omega \right)}+\left\| \boldsymbol{u}_m \right\| _{L^8\left( \Omega \right)}\left\| \Delta \chi _m \right\| _{L^2\left( \Omega \right)} \right) \nonumber
\\
&\le C\left( E_0 \right) \left( \left\| \nabla \partial _t\chi _m \right\| _{L^2\left( \Omega \right)}^{\frac{3}{4}}\left\| \partial _t\chi _m \right\| _{L^2\left( \Omega \right)}^{\frac{1}{4}}+\left\| \partial _t\chi _m \right\| _{L^2\left( \Omega \right)} \right) \left\| \nabla \boldsymbol{u}_m \right\| _{L^2\left( \Omega \right)}^{\frac{3}{2}}\left\| \Delta \chi _m \right\| _{L^2\left( \Omega \right)}^{\frac{3}{4}}
\\
	&\quad +C\left( E_0 \right) \left( \left\| \nabla \partial _t\chi _m \right\| _{L^2\left( \Omega \right)}^{\frac{3}{4}}\left\| \partial _t\chi _m \right\| _{L^2\left( \Omega \right)}^{\frac{1}{4}}+\left\| \partial _t\chi _m \right\| _{L^2\left( \Omega \right)} \right) \left\| \nabla \boldsymbol{u}_m \right\| _{L^2\left( \Omega \right)}^{\frac{5}{4}}\left\| \Delta \chi _m \right\| _{L^2\left( \Omega \right)}\nonumber
\\
&\le \frac{m_*}{16}\left\| \nabla \partial _t\chi _m \right\| _{L^2\left( \Omega \right)}^{2}
     +C\left( E_0 \right) \left\| \Delta \chi _m \right\| _{L^2\left( \Omega \right)}^2+C\left( E_0 \right)\nonumber
\\
&\quad
+C\left( E_0 \right) \left( \left\| \nabla \boldsymbol{u}_m \right\| _{L^2\left( \Omega \right)}^{2}+\left\| \Delta \chi _m \right\| _{L^2\left( \Omega \right)}^{2} \right) \left( \left\| \partial _t\chi _m \right\| _{L^2\left( \Omega \right)}^{2}+\left\| \nabla \boldsymbol{u}_m \right\| _{L^2\left( \Omega \right)}^{2} \right),& \nonumber
\end{flalign*}
\begin{flalign*}
\qquad\quad	\left| I_4 \right|&\le C\left\| \partial _t\chi _m \right\| _{L^4\left( \Omega \right)}\left\| \partial _t\boldsymbol{u}_m \right\| _{L^2\left( \Omega \right)}\left\| \nabla \chi _m \right\| _{L^4\left( \Omega \right)} \nonumber\\
	&\le C\left( E_0 \right) \left( \left\| \nabla \partial _t\chi _m \right\| _{L^2\left( \Omega \right)}^{\frac{1}{2}}\left\| \partial _t\chi _m \right\| _{L^2\left( \Omega \right)}^{\frac{1}{2}}+\left\| \partial _t\chi _m \right\| _{L^2\left( \Omega \right)} \right) \left\| \partial _t\boldsymbol{u}_m \right\| _{L^2\left( \Omega \right)}\left\| \Delta \chi _m \right\| _{L^2\left( \Omega \right)}^{\frac{1}{2}}\\
	&\le \frac{m_*}{16}\left\| \nabla \partial _t\chi _m \right\| _{L^2\left( \Omega \right)}^{2}+\frac{\rho _*}{8}\left\| \partial _t\boldsymbol{u}_m \right\| _{L^2\left( \Omega \right)}^{2}+C\left( E_0 \right)\left( \left\| \Delta \chi _m \right\| _{L^2\left( \Omega \right)}^{2} +1\right) \left\| \partial _t\chi _m \right\| _{L^2\left( \Omega \right)}^{2},& \nonumber
\end{flalign*}
\begin{flalign*}
\qquad\quad	\left| I_5 \right|&\le C\left\| \partial _t\chi _m \right\| _{L^4\left( \Omega \right)}\left\| \boldsymbol{u}_m \right\| _{L^4\left( \Omega \right)}\left\| \nabla \partial _t\chi _m \right\| _{L^2\left( \Omega \right)} \nonumber\\
	&\le C\left( E_0 \right) \left( \left\| \nabla \partial _t\chi _m \right\| _{L^2\left( \Omega \right)}^{\frac{1}{2}}\left\| \partial _t\chi _m \right\| _{L^2\left( \Omega \right)}^{\frac{1}{2}}+\left\| \partial _t\chi _m \right\| _{L^2\left( \Omega \right)} \right) \left\| \nabla \boldsymbol{u}_m \right\| _{L^2\left( \Omega \right)}^{\frac{1}{2}}\left\| \nabla \partial _t\chi _m \right\| _{L^2\left( \Omega \right)}\\
	&\le \frac{m_*}{16}\left\| \nabla \partial _t\chi _m \right\| _{L^2\left( \Omega \right)}^{2}+C\left( E_0 \right) \left( \left\| \nabla \boldsymbol{u}_m \right\| _{L^2\left( \Omega \right)}^{2}+1 \right) \left\| \partial _t\chi _m \right\| _{L^2\left( \Omega \right)}^{2},& \nonumber
\end{flalign*}
\begin{flalign*}
\qquad\quad	\left| I_6 \right|&\le C\left\| \partial _t\chi _m \right\| _{L^4\left( \Omega \right)}^{2}\left\| \Delta \chi _m \right\| _{L^2\left( \Omega \right)}\nonumber\\
	&\le C\left( \left\| \nabla \partial _t\chi _m \right\| _{L^2\left( \Omega \right)}\left\| \partial _t\chi _m \right\| _{L^2\left( \Omega \right)}+\left\| \partial _t\chi _m \right\| _{L^2\left( \Omega \right)}^{2} \right) \left\| \Delta \chi _m \right\| _{L^2\left( \Omega \right)}\\
	&\le \frac{m_*}{16}\left\| \nabla \partial _t\chi _m \right\| _{L^2\left( \Omega \right)}^{2}+C\left( \left\| \Delta \chi _m \right\| _{L^2\left( \Omega \right)}^{2}+ 1 \right) \left\| \partial _t\chi _m \right\| _{L^2\left( \Omega \right)}^{2},& \nonumber
\end{flalign*}
\begin{flalign*}
\qquad\quad \left| I_7 \right| & \le C\left\| \partial _t\chi _m \right\| _{L^4\left( \Omega \right)}\left\| \nabla \partial _t\chi _m \right\| _{L^2\left( \Omega \right)}\left\| \nabla \chi _m \right\| _{L^4\left( \Omega \right)} \nonumber\\
	&\le C\left( E_0 \right) \left( \left\| \nabla \partial _t\chi _m \right\| _{L^2\left( \Omega \right)}^{\frac{1}{2}}\left\| \partial _t\chi _m \right\| _{L^2\left( \Omega \right)}^{\frac{1}{2}}+\left\| \partial _t\chi _m \right\| _{L^2\left( \Omega \right)} \right) \left\| \nabla \partial _t\chi _m \right\| _{L^2\left( \Omega \right)}\left\| \Delta \chi _m \right\| _{L^2\left( \Omega \right)}^{\frac{1}{2}}\\
	&\le \frac{m_*}{16}\left\| \nabla \partial _t\chi _m \right\| _{L^2\left( \Omega \right)}^{2}+C\left( E_0 \right) \left( \left\| \Delta \chi _m \right\| _{L^2\left( \Omega \right)}^{2}+ 1 \right) \left\| \partial _t\chi _m \right\| _{L^2\left( \Omega \right)}^{2},& \nonumber
\end{flalign*}

Recalling the equation $(\ref{6.3})$ we get
\begin{align*}
\|\Delta\chi_m\|_{L^2(\Omega)}^2
\le &C\|\partial_t\chi_m\|_{L^2(\Omega)}^2
+C\|\boldsymbol{u}_m\|_{L^4(\Omega)}^2\|\nabla\chi_m\|_{L^4(\Omega)}^2+C
\\
\le&C\|\partial_t\chi_m\|_{L^2(\Omega)}^2
+C(E_0)\|\nabla\boldsymbol{u}_m\|_{L^2(\Omega)}\|\Delta\chi_m\|_{L^2(\Omega)}+C
\\
\le&\frac12\|\Delta\chi_m\|_{L^2(\Omega)}^2
+C\|\partial_t\chi_m\|_{L^2(\Omega)}^2
+C(E_0)\|\nabla\boldsymbol{u}_m\|_{L^2(\Omega)}^2+C,
\end{align*}
i.e.
\begin{align}\label{phi-2-2}
\|\Delta\chi_m\|_{L^2(\Omega)}^2
\le C\|\partial_t\chi_m\|_{L^2(\Omega)}^2
+C(E_0)\|\nabla\boldsymbol{u}_m\|_{L^2(\Omega)}^2+C.
\end{align}
Recalling the equation (\ref{phi-phi}) and using $(\ref{6.6})$, (\ref{phi-2-2}), we see that
\begin{align}      \label{6.55}   
&m_*\left\|{\Delta \chi _m \nabla \chi _m} \right\| _{L^2\left( \Omega \right)}^{2} \nonumber
\\
&\le C\left\| \partial _t\chi _m\nabla \chi _m \right\| _{L^2\left( \Omega \right)}^{2}+C\left\| \left(\boldsymbol{u}_m \cdot \nabla \chi _m \right)\nabla \chi _m \right\| _{L^2\left( \Omega \right)}^{2}+C\left\| F ^{\prime}\left( \chi _m \right) \nabla \chi _m \right\| _{L^2\left( \Omega \right)}^{2}\nonumber
\\
	&\le C\left\| \partial _t\chi _m \right\| _{L^4\left( \Omega \right)}^{2}\left\| \nabla \chi _m \right\| _{L^4\left( \Omega \right)}^{2}+C\left\| \boldsymbol{u}_m \right\| _{L^6\left( \Omega \right)}^{2}\left\| \nabla \chi _m \right\| _{L^6\left( \Omega \right)}^{4}\nonumber
+C\left\| \nabla \chi _m \right\| _{L^2\left( \Omega \right)}^{2}\nonumber
    \\
	&\le C\left( E_0 \right) \left( \left\| \nabla \partial _t\chi _m \right\| _{L^2\left( \Omega \right)}\left\| \partial _t\chi _m \right\| _{L^2\left( \Omega \right)}+\left\| \partial _t\chi _m \right\| _{L^2\left( \Omega \right)}^{2} \right) \left\| \Delta \chi _m \right\| _{L^2\left( \Omega \right)} \nonumber
\\
	&\quad +C\left( E_0 \right) \left\| \nabla \boldsymbol{u}_m \right\| _{L^2\left( \Omega \right)}^\frac{4}{3}\left\| \Delta \chi _m \right\| _{L^2\left( \Omega \right)}^\frac{8}{3}
+C\left( E_0 \right)
\\
	&\le\varepsilon \left\| \nabla \partial _t\chi _m \right\| _{L^2\left( \Omega \right)}^{2}+C\left( E_0 \right) \left( \left\| \Delta \chi _m \right\| _{L^2\left( \Omega \right)}^{2}+\left\| \Delta \chi _m \right\| _{L^2\left( \Omega \right)} \right) \left\| \partial _t\chi _m \right\| _{L^2\left( \Omega \right)}^{2} \nonumber
\\
    &\quad +C\left( E_0 \right) \left\| \Delta \chi _m \right\| _{L^2\left( \Omega \right)}^{2}
    \left\| \nabla \boldsymbol{u}_m \right\| _{L^2\left( \Omega \right)}^2
    +C\left( E_0 \right) \left\| \Delta \chi _m \right\| _{L^2\left( \Omega \right)}^4+C\left( E_0 \right) \nonumber
\\
	&\le\varepsilon \left\| \nabla \partial _t\chi _m \right\| _{L^2\left( \Omega \right)}^{2}
+C\left( E_0 \right) \left( \left\| \Delta \chi _m \right\| _{L^2\left( \Omega \right)}^{2}+1\right)
\left(\left\| \partial _t\chi _m \right\| _{L^2\left( \Omega \right)}^{2}
+ \left\| \nabla \boldsymbol{u}_m \right\| _{L^2\left( \Omega \right)}^2\right)
+C\left( E_0 \right). \nonumber
\end{align}
Taking $\varepsilon=\frac{m _*}{8}$ here, we have
\begin{align}      \label{6.56}   
\left|  I_{13}  \right|
&\le \frac{\rho _*}{4}\left\|\partial _t \boldsymbol{u}_m \right\| _{L^2\left( \Omega \right)}^{2}
+\frac{m _*}{8} \left\| \nabla \partial _t\chi _m \right\| _{L^2\left( \Omega \right)}^{2}
+C\left( E_0 \right)  \nonumber
\\
	&\quad +C\left( E_0 \right) \left( \left\| \Delta \chi _m \right\| _{L^2\left( \Omega \right)}^{2}+1\right)
\left(\left\| \partial _t\chi _m \right\| _{L^2\left( \Omega \right)}^{2}
+ \left\| \nabla \boldsymbol{u}_m \right\| _{L^2\left( \Omega \right)}^2\right).
\end{align}
{For $u\in H_0^1(\Omega)\cap H^2(\Omega)$}, by the 2-D Gagliardo-Nirenberg inequality
\begin{align*}
\left\| u \right\| _{L^{\infty}(\Omega )}\le C\left\| \Delta u \right\| _{L^2\left( \Omega \right)}^{\frac{1}{3}}\left\| u \right\| _{L^4(\Omega )}^{\frac{2}{3}}
\end{align*}
and Ladyzhenskaya's inequality in 2D {and Poincar\'{e}'s inequality},
\begin{align*}
\left\| u \right\| _{L^4(\Omega )}^{2}\le C\left\| \nabla u \right\| _{L^2\left( \Omega \right)}\left\| u \right\| _{L^2\left( \Omega \right)},
\end{align*}
we infer that
\begin{align}      \label{6.57}   
\left\| u \right\| _{L^{\infty}(\Omega )}\le C\left\| \Delta u \right\| _{L^2\left( \Omega \right)}^{\frac{1}{3}}\left\| \nabla u \right\| _{L^2\left( \Omega \right)}^{\frac{1}{3}}\left\| u \right\| _{L^2\left( \Omega \right)}^{\frac{1}{3}}.
\end{align}
By exploiting $(\ref{6.4})$, $(\ref{6.5})$ and $(\ref{6.57})$, we have
\begin{align}      \label{6.58}   
	\left|  I_{14}  \right|&\le \rho^*
\left\| \boldsymbol{u}_m \right\| _{L^{\infty}\left( \Omega \right)}
\left\| \nabla \boldsymbol{u}_m \right\| _{L^2\left( \Omega \right)}
\left\| \partial _t\boldsymbol{u}_m \right\| _{L^2\left( \Omega \right)}  \nonumber\\
	&\le \frac{\rho _*}{4}\left\| \partial _t\boldsymbol{u}_m \right\| _{L^2\left( \Omega \right)}^{2}+C\left\| \boldsymbol{u}_m \right\| _{L^{\infty}\left( \Omega \right)}^{2}\left\| \nabla \boldsymbol{u}_m \right\| _{L^2\left( \Omega \right)}^{2}\\
	&\le \frac{\rho _*}{4}\left\| \partial _t\boldsymbol{u}_m \right\| _{L^2\left( \Omega \right)}^{2}+C\left( E_0 \right) \left\| \Delta \boldsymbol{u}_m \right\| _{L^2\left( \Omega \right)}^{\frac{2}{3}}\left\| \nabla \boldsymbol{u}_m \right\| _{L^2\left( \Omega \right)}^{\frac{8}{3}} \nonumber\\
	&\le \frac{\rho _*}{4}\left\| \partial _t\boldsymbol{u}_m \right\| _{L^2\left( \Omega \right)}^{2}
+\varpi_1 \left\| \mathbf{A}\boldsymbol{u}_m \right\| _{L^2\left( \Omega \right)}^{2}
+C\left( E_0 \right) \left\| \nabla \boldsymbol{u}_m \right\| _{L^2\left( \Omega \right)}^{4},\nonumber
\end{align}
where $\varpi_1$ is a positive (small) constant which will be determined later.

Noticing that ${\eta\equiv 1} $ in $(\ref{6.2})$, putting the above estimates on $I_i ~(i=1, \cdots,  14 )$ into (\ref{6.9}) gives
\begin{align}      \label{6.53}   
	&\frac{\mathrm{d}}{\mathrm{d}t} \left\{\int_{\Omega}{\frac{1}{2}\rho _{m}^{2}\left| \partial _t\chi _m \right|^2}\mathrm{d}x+\int_{\Omega}{\frac{1}{2}\eta \left( \chi _m \right) \left| \mathbb{D} \boldsymbol{u}_m \right|^2}\mathrm{d}x\right\}+\frac{m _*}{2}\left\| \nabla \partial _t\chi _m  \right\| _{L^2\left( \Omega \right)}^{2}+\frac{\rho _*}{2}\left\| \partial _t\boldsymbol{u}_m \right\| _{L^2\left( \Omega \right)}^{2} \nonumber
\\
	&\le \varpi_1 \left\| \mathbf{A}\boldsymbol{u}_m \right\| _{L^2\left( \Omega \right)}^{2}+\frac{\rho _*}{8}\left\| \partial _t\boldsymbol{u}_m \right\| _{L^2\left( \Omega \right)}^{2}+\frac{m_*}{8}\left\| \nabla \partial _t\chi _m \right\| _{L^2\left( \Omega \right)}^{2}+C\left( E_0 \right)\left\| \Delta \chi _m \right\| _{L^2\left( \Omega \right)}^{2} \\
	&\quad +C\left( E_0 \right) \left( \left\| \Delta \chi _m \right\| _{L^2\left( \Omega \right)}^{2}
+\left\| \nabla \boldsymbol{u}_m \right\| _{L^2\left( \Omega \right)}^2+1\right)
\left(\left\| \partial _t\chi _m \right\| _{L^2\left( \Omega \right)}^{2}
+ \left\| \nabla \boldsymbol{u}_m \right\| _{L^2\left( \Omega \right)}^2\right)+C\left( E_0 \right). \nonumber
\end{align}
Here the constant $C$ depends on $\varpi_1$. Since $\eta\equiv1$, we rewrite $(\ref{6.31})$ as follows
\begin{align}      \label{6.60}   
	\int_{\Omega}{\frac{1}{2}\left| \mathbf{A}\boldsymbol{u}_m \right|^2}\mathrm{d}x
=&\left( -\rho _m\partial _t\boldsymbol{u}_m,\mathbf{A}\boldsymbol{u}_m \right)
+\left( -\rho _m\left( \boldsymbol{u}_m\cdot \nabla \right) \boldsymbol{u}_m,\mathbf{A}\boldsymbol{u}_m \right)  \nonumber
\\
&+\left( -\Delta\chi _m\nabla \chi _m, \mathbf{A}\boldsymbol{u}_m \right)
=\sum_{i=1}^3 K_i.
\end{align}
Using $(\ref{6.4})$-$(\ref{6.6})$, $(\ref{6.30})$ and $(\ref{6.55})$, we estimate $K_{1}$, $K_{2}$, $K_{3}$ as follows
\begin{flalign}      \label{6.61}   
\qquad\qquad \left| K_1 \right|&\le \rho ^*\left\| \partial _t\boldsymbol{u}_m \right\| _{L^2\left( \Omega \right)}\left\| \mathbf{A}\boldsymbol{u}_m \right\| _{L^2\left( \Omega \right)} \le \frac{1}{12}\left\| \mathbf{A}\boldsymbol{u}_m \right\| _{L^2\left( \Omega \right)}^{2}+3\left( \rho ^* \right) ^2 \left\| \partial _t\boldsymbol{u}_m \right\| _{L^2\left( \Omega \right)}^{2},&
\\      \label{6.62}   
\qquad\qquad \left| K_2 \right|&\le \rho ^*\left\| \boldsymbol{u}_m \right\| _{L^4\left( \Omega \right)}\left\| \nabla \boldsymbol{u}_m \right\| _{L^4\left( \Omega \right)}\left\| \mathbf{A}\boldsymbol{u}_m \right\| _{L^2\left( \Omega \right)} \nonumber\\
	&\le C\left( E_0 \right) \left\| \nabla \boldsymbol{u}_m \right\| _{L^2\left( \Omega \right)}\left\| \mathbf{A}\boldsymbol{u}_m \right\| _{L^2\left( \Omega \right)}^{\frac{3}{2}}\\
	&\le \frac{1}{12}\left\| \mathbf{A}\boldsymbol{u}_m \right\| _{L^2\left( \Omega \right)}^{2}+C\left( E_0 \right) \left\| \nabla \boldsymbol{u}_m \right\| _{L^2\left( \Omega \right)}^{4},& \nonumber
\\     \label{6.63}   
\qquad\qquad \left| K_3 \right|&\le \frac{1}{12}\left\| \mathbf{A}\boldsymbol{u}_m \right\| _{L^2\left( \Omega \right)}^{2}+C \left\|{\Delta \chi _m \nabla \chi _m }\right\| _{L^2\left( \Omega \right)}^{2}
\\
    &\le \frac{1}{12}\left\| \mathbf{A}\boldsymbol{u}_m \right\| _{L^2\left( \Omega \right)}^{2}
    +\varepsilon \left\| \nabla \partial _t\chi _m \right\| _{L^2\left( \Omega \right)}^{2}
    +C\left( E_0 \right) \nonumber
    \\
	&\quad +C\left( E_0 \right) \left( \left\| \Delta \chi _m \right\| _{L^2\left( \Omega \right)}^{2}+1\right)
\left(\left\| \partial _t\chi _m \right\| _{L^2\left( \Omega \right)}^{2}
+ \left\| \nabla \boldsymbol{u}_m \right\| _{L^2\left( \Omega \right)}^2\right).& \nonumber
\end{flalign}
Combining the above estimates $(\ref{6.61})$-$(\ref{6.63})$ with $(\ref{6.60})$, we arrive at
\begin{align}      \label{6.64}   
	\frac{1}{4}\left\| \mathbf{A}\boldsymbol{u}_m \right\| _{L^2\left( \Omega \right)}^{2}&\le 3\left( \rho ^* \right) ^2\left\| \partial _t\boldsymbol{u}_m \right\| _{L^2\left( \Omega \right)}^{2}+\varepsilon \left\| \nabla \partial _t\chi _m \right\| _{L^2\left( \Omega \right)}^{2}+C\left( E_0 \right)
\\
	&\quad +C\left( E_0 \right) \left( \left\| \Delta \chi _m \right\| _{L^2\left( \Omega \right)}^{2}
+\left\| \nabla \boldsymbol{u}_m \right\| _{L^2\left( \Omega \right)}^2+1\right)
\left(\left\| \partial _t\chi _m \right\| _{L^2\left( \Omega \right)}^{2}
+ \left\| \nabla \boldsymbol{u}_m \right\| _{L^2\left( \Omega \right)}^2\right). \nonumber
\end{align}
Taking $\delta_1 =\frac{\rho _*}{24\left( \rho ^* \right) ^2}$ and setting $\varpi_1 =\frac{\delta_1}{8}$,
$\varepsilon =\frac{m_*}{8\delta_1}$, {then} multiplying the inequality $(\ref{6.64})$ by $\delta_1$, we have
\begin{align}      \label{6.66}   
	2\varpi_1 \left\| \mathbf{A}\boldsymbol{u}_m \right\| _{L^2\left( \Omega \right)}^{2} & \le \frac{\rho _*}{8}\left\| \partial _t\boldsymbol{u}_m \right\| _{L^2\left( \Omega \right)}^{2}
+\frac{m_*}{8}\left\| \nabla \partial _t\chi _m \right\| _{L^2\left( \Omega \right)}^{2}
+C\left( E_0 \right)
\\
	&\quad +C\left( E_0 \right) \left( \left\| \Delta \chi _m \right\| _{L^2\left( \Omega \right)}^{2}
+\left\| \nabla \boldsymbol{u}_m \right\| _{L^2\left( \Omega \right)}^2+1\right)
\left(\left\| \partial _t\chi _m \right\| _{L^2\left( \Omega \right)}^{2}
+ \left\| \nabla \boldsymbol{u}_m \right\| _{L^2\left( \Omega \right)}^2\right). \nonumber
\end{align}
Set
\begin{align*}
&\tilde{H}_m=\tilde{H}\left( \rho _m,\boldsymbol{u}_m,\chi _m \right) =\int_{\Omega}{\frac{1}{2}\rho _{m}^{2}\left| \partial _t\chi _m \right|^2}\mathrm{d}x+\int_{\Omega}{{\frac{1}{2}}\left| \mathbb{D} \boldsymbol{u}_m \right|^2}\mathrm{d}x,
\\
&\tilde{F}_m=\tilde{F}\left( \boldsymbol{u}_m,\chi _m \right) =\varpi_1 \left\| \mathbf{A}\boldsymbol{u}_m \right\| _{L^2\left( \Omega \right)}^{2}+\frac{\rho _*}{4}\left\| \partial _t\boldsymbol{u}_m \right\| _{L^2\left( \Omega \right)}^{2}+\frac{m_*}{4}\left\| \nabla \partial _t\chi _m \right\| _{L^2\left( \Omega \right)}^{2}.
\end{align*}
Putting $(\ref{6.53})$ and $(\ref{6.66})$ together, we obtain
\begin{align}      \label{6.67}   
\frac{\mathrm{d}}{\mathrm{d}t}\tilde{H}_m+\tilde{F}_m\le C\left( E_0 \right) \left( \left\| \nabla \boldsymbol{u}_m \right\| _{L^2(\Omega )}^{2}+\left\| \Delta \chi _m \right\| _{L^2(\Omega )}^{2}+1 \right) \tilde{H}_m+C\left( E_0 \right) \left\| \Delta \chi _m \right\| _{L^2\left( \Omega \right)}^{2}+C\left( E_0 \right).
\end{align}
\vskip1.8mm
\noindent
{\it\bfseries Step 3. Uniform estimates and passage to the limit}\quad
If $d=3$,  {Bihari's inequality applied to $(\ref{6.40})$ gives
\[
C(E_0)+H_m(t)
\le
\frac{C(E_0)+H_m(0)}
{\left[
1-4C(E_0)\bigl(C(E_0)+H_m(0)\bigr)^4t
\right]^{1/4}},
\]
as long as
\[
4C(E_0)\bigl(C(E_0)+H_m(0)\bigr)^4t<1.
\] }
Since
\begin{align*}
H_m(0)\le \mathrm{C}\left\| \partial _t\chi(0) \right\| _{L^2(\Omega )}^{2}+\mathrm{C}\left\| \nabla \boldsymbol{u}_m(0) \right\| _{L^2\left( \Omega \right)}^{2}
\le C\left(\left\|\boldsymbol{u}_0\right\|_{\mathbf{v}_{\sigma}}^2+\left\|\chi_0\right\|_{H^2(\Omega)}^2\right),
\end{align*}
we infer that there exist $T_{0}>0$ and $C_{1}>0$, depending only on the norms $\left\|\boldsymbol{u}_{0}\right\|_{\mathbf{v}_{\sigma}}$,
$\left\|\chi_0\right\|_{H^2(\Omega)}$, $\rho_{*}$, $\rho^{*}$, $\eta_{*}$ and $\eta^{*}$, such that
\begin{align}      \label{6.68}   
\mathop {\mathrm{sup}}_{t\in \left[ 0,T_0 \right]}\left( \left\| \partial _t\chi _m(t) \right\| _{L^2(\Omega )}^{2}+\left\| \nabla \boldsymbol{u}_m(t) \right\| _{L^2(\Omega )}^{2} \right) \le C_1\quad \,\,\mathrm{if} ~d=3,
\end{align}
where $T_{0}$ and $C_{1}$ are independent of $m$.

If $d=2$, {by Gr\"{o}nwall's inequality applied to} $(\ref{6.67})$ we deduce that
\begin{align*}
\tilde{H}_m(t)
\le& \left( \tilde{H}_m(0)+\int_0^T{\mathrm{C}\left( E_0 \right) \left( \left\| \Delta \chi _m \right\| _{L^2\left( \Omega \right)}^{2}+1 \right)}\mathrm{d}t \right)
\\
&\cdot\exp\left\{C(E_0) \int_0^T{\left( \left\| \nabla \boldsymbol{u}_m \right\| _{L^2\left( \Omega \right)}^{2}+\left\| \Delta \chi _m \right\| _{L^2\left( \Omega \right)}^{2}+1 \right)}\mathrm{d}t\right\},
\quad\forall t\in \left[ 0,T \right].
\end{align*}
Then, recalling $(\ref{6.5})$ and $(\ref{6.6})$, for any $T>0$ there exists $C_{2}=C_{2}(T)$, also depending on the norms
$\left\|\boldsymbol{u}_{0}\right\|_{\mathbf{v}_{\sigma}}$, $\rho_{*}$, $\rho^{*}$, $\eta_{*}$ and $\eta^{*}$, such that
\begin{align*}     
\mathop {\mathrm{sup}}_{t\in [0,T]}\left( \left\| \partial _t\chi _m(t) \right\| _{L^2(\Omega )}^{2}+\left\| \nabla \boldsymbol{u}_m(t) \right\| _{L^2(\Omega )}^{2} \right) \le C_2\quad \mathrm{if}~d=2.
\end{align*}
We will now use the notations $\bar{C}_{1}$ and $\bar{C}_{2}$ to denote generic constants depending on $C_{1}$ and $C_{2}$, but independent of $m$. We also infer from $(\ref{6.40})$ that
\begin{align}      \label{6.70}   
\int_0^{T_0}{\left( \left\| \mathbf{A}\boldsymbol{u}_m(\tau ) \right\| _{L^2(\Omega )}^{2}+\left\| \partial _t\boldsymbol{u}_m(\tau ) \right\| _{L^2(\Omega )}^{2}+\left\| \nabla \partial _t\chi _m(\tau ) \right\| _{L^2(\Omega )}^{2} \right)}\mathrm{d}\tau \le \bar{C}_1\quad \mathrm{if}~ d=3,
\end{align}
and from $(\ref{6.67})$ that
\begin{align}      \label{6.71}   
\int_0^{T}{\left( \left\| \mathbf{A}\boldsymbol{u}_m(\tau ) \right\| _{L^2(\Omega )}^{2}+\left\| \partial _t\boldsymbol{u}_m(\tau ) \right\| _{L^2(\Omega )}^{2}+\left\| \nabla \partial _t\chi _m(\tau ) \right\| _{L^2(\Omega )}^{2} \right)}\mathrm{d}\tau \le \bar{C}_2\quad \mathrm{if}~ d=2.
\end{align}
More precisely, thanks to Lemma \ref{Lions}, $(\ref{6.4})$-$(\ref{6.6})$, (\ref{phi-2}) and the above estimates $(\ref{6.68})$-$(\ref{6.71})$, we deduce that (up to a subsequence)
\begin{alignat}{2}     \label{6.72}   
&\rho_{m} \rightharpoonup \rho && \text { weakly$^*$ in } L^{\infty}(\Omega \times(0, T_0)), \nonumber
\\
&\rho_{m} \rightarrow \rho &&\text { strongly in } \mathcal{C}\left([0, T_0] ; L^{r}(\Omega)\right), \forall r \in[1, \infty), \nonumber
\\
&\boldsymbol{u}_{m} \rightharpoonup \boldsymbol{u} && \text { weakly$^*$ in } L^{\infty}\left(0, T_0 ; \mathbf{V}_{\sigma}\right), \nonumber
\\
&\boldsymbol{u}_{m} \rightharpoonup \boldsymbol{u} && \text { weakly in } L^{2}(0, T_0 ; H^{2}(\Omega))\cap H^{1}\left(0, T_0 ; \mathbf{H}_{\sigma}\right),
\\
&\chi_{m} \rightharpoonup \chi && \text { weakly$^*$ in } L^{\infty}(0, T_0; {H^{2}(\Omega)}),   \nonumber
\\
&\partial_t\chi_{m} \rightharpoonup \partial_t\chi && \text { weakly in }  L^{2}(0, T_0 ; H^{1}(\Omega)), \nonumber
\\
&{\partial_t\chi_{m} \rightharpoonup \partial_t\chi} && \text { weakly$^*$ in }
{L^{\infty}(0, T_0; L^{2}(\Omega))}. \nonumber
\end{alignat}
Thanks to $(\ref{6.72})$, we infer from the Aubin--Lions lemma that
\begin{align}        
\boldsymbol{u}_{m} \rightarrow \boldsymbol{u} & \quad
\text { strongly in } { L^2\left([0, T_0] ; \mathbf{V}_{\sigma}\right)\cap \mathcal{C}([0,T_0]; H^s), ~0\le s<1}, \label{6.73}
\\
\chi_{m} \rightarrow \chi & \quad \text { strongly in } \mathcal{C}([0, T_0] ; W^{1, q}(\Omega)), \forall q \in[2,6),\label{6.74}
\end{align}
We notice that, if $d=2$, $(\ref{6.72})$--$(\ref{6.74})$ hold by replacing $\left[0, T_{0}\right]$ with $[0, T]$, and $(\ref{6.72})$--$(\ref{6.74})$ imply the convergence of the nonlinear terms. Thus, in a standard manner, we can pass to the limit as $m \rightarrow \infty$ in $(\ref{6.2})$--$(\ref{6.3})$ and obtain a limit solution $(\rho, \boldsymbol{u}, p, \chi)$ satisfying $(\ref{1.1})$--$(\ref{1.2})$ as stated in Theorem $\ref{Th6.1}$.
\hfill$\Box$

\section{Global strong solution in three dimensions}
While \cite{G-T} established global weak, 2D global and 3D local strong solutions
for the NSCH system, and \cite{F-N-G} subsequently proved global 3D strong solutions
under the smallness condition \eqref{eq:JDE-smallness}, the global 3D strong solutions and
exponential decay rates here are the first such results for the fully
coupled NSAC system with density-dependent Allen-Cahn dynamics.
In this section, we prove that, under smallness assumptions on the initial data, the 3D local strong solution obtained in Theorem $\ref{Th6.1}$ can be extended to a global one on $[0, \infty)$.

We first define the following energy functional
\begin{align}      \label{8.1}
\mathcal{E}(t)
=\left\| \boldsymbol{u}\right\| _{H^1(\Omega )}^{2}
+\left\| \chi ^2-1 \right\| _{L^2(\Omega )}^{2}
+\left\| \nabla \chi \right\| _{H^1(\Omega )}^{2}
+\left\| \partial _t\chi \right\| _{L^2(\Omega )}^{2},
\end{align}
and dissipative functional
\begin{align}      \label{8.2}
\mathcal{D}(t)
=&\left\| \partial _t\boldsymbol{u} \right\| _{L^2(\Omega )}^{2}
+\left\|  \boldsymbol{u}\right\| _{H^1(\Omega )}^{2}
+\left\| \partial _t\chi \right\| _{H^1(\Omega )}^{2}
+\left\| \chi ^2-1 \right\| _{L^2(\Omega )}^{2}
+\left\|  \nabla \chi \right\| _{H^1(\Omega )}^{2}
+\left\|  \mu \right\| _{L^2(\Omega )}^{2}.
\end{align}
\begin{Theorem}       
\label{Th8.1}
Let $d=3$ and $\eta \equiv 1$.
We assume that $\left. (\rho, {\boldsymbol{u}},\chi ) \right|_{t=0}=\left( \rho _0, {\boldsymbol{u}}_0,\chi _0 \right) \in L^{\infty}(\Omega )\times \mathbf{V}_{\sigma}(\Omega )\times H^2(\Omega)$ such that $0<\rho_{*} \le \rho_{0} \le \rho^{*}$, $\left|\chi_{0}\right| \le 1$ a.e. in $\Omega$, $\partial_n\chi_0=0$ on $\partial\Omega$. There exists a small positive constant $\varepsilon_{0}$ such that if
\begin{align}      \label{8.3}
\left\|{\boldsymbol{u}}_{0}\right\|_{\mathbf{V}_{\sigma}(\Omega)}+\left\|\chi_{0}^{2}-1 \right\|_{L^{2}(\Omega)}+\left\|\nabla \chi_{0}\right\|_{H^{1}(\Omega)} \le \varepsilon_{0},
\end{align}
then the problem $(\ref{1.1})$--$(\ref{1.2})$ has a global-in-time strong solution $(\boldsymbol{u}, \chi)$ in three dimensions with
\begin{align}      \label{8.4}
\mathcal{E}(t)+\int_0^t{\mathcal{D}}(s)\mathrm{d}s\le C\varepsilon _0,
\end{align}
for all $t>0$. In addition, there exists a positive constant $\sigma$ such that
\begin{align*}
\left\|{\boldsymbol{u}}(\cdot,t)\right\|_{L^2(\Omega)}^2
+\left\|\chi^2(\cdot,t)-1\right\|_{L^2(\Omega)}^2
+\left\|\nabla\chi(\cdot,t)\right\|_{L^2(\Omega)}^2
\le C{\mathcal E}(0)e^{-\sigma t},
\quad t\ge0,
\end{align*}
where $\mathcal{E} (0)=\left\| \boldsymbol{u}_0 \right\| _{H^1(\Omega )}^{2}+\left\| \chi _{0}^{2}-1 \right\| _{L^2(\Omega )}^{2}+\left\| \nabla \chi _0 \right\| _{H^1(\Omega )}^{2}+\left\| \partial _t\chi (0) \right\| _{L^2(\Omega )}^{2}\le C_0 \varepsilon _0$, with $\partial_t\chi(0)$ determined by the phase equation.
\end{Theorem}   

\begin{Remark}[{\bf Comparison with the smallness condition in \cite{F-N-G}}]\label{R-small}
{The smallness condition \eqref{8.3} places the initial
velocity near zero and the initial phase field near one
of the spatially homogeneous pure phases $\chi=\pm1$.
Thus, the global existence and decay asserted in
Theorem~\ref{Th8.1} are consistent with the expected
stability of these equilibria.}
The smallness condition in \cite{F-N-G} for the NSCH system requires
\begin{equation}\label{JDE-smallness}
\|\nabla u_0\|_{L^2(\Omega)}+\|\nabla\mu_0\|_{L^2(\Omega)}+\|\rho_0\|_{L^\infty(\Omega)}\leq\varepsilon_0.
\end{equation}
The fourth-order Cahn-Hilliard structure provides $\nabla\mu\in L^2(0,T;L^2(\Omega))$, which yields sufficient regularity to control the capillary stress.
The smallness condition in \cite{F-N-G} involves the
chemical-potential gradient and the magnitude of the
initial density.
In contrast, \eqref{8.3} is expressed directly in terms
of the initial velocity and phase field, while the
density is only required to remain between fixed
positive bounds.

In the present NSAC system, the energy estimate yields $\mu\in L^2(0,T;L^2(\Omega))$ without gradient information, and the techniques for controlling smallness in \cite{F-N-G} do not directly apply. This reflects the different structural roles of the chemical potential in the second-order Allen-Cahn and fourth-order Cahn-Hilliard dynamics.
\end{Remark}

We first give estimates for the 3D local-in-time strong solution obtained in Theorem $\ref{Th6.1}$, which will be used in the proof of
Theorem \ref{Th8.1}.
\begin{Proposition}   \label{R8.1}
Let $d=3$. Under the assumptions of Theorem \ref{Th8.1}, if there exists a universal positive constant $\delta_0$ such that $\mathcal{E}(0) \leq \delta_0$, then the solution $(\boldsymbol{u}, \chi)$ to the problem $(\ref{1.1})$--$(\ref{1.2})$
obtained in Theorem $\ref{Th6.1}$ on the interval $[0, T_0]$ satisfies
$$
\sup _{0 \leq t \leq T} \mathcal{E}(t)+\int_0^T \mathcal{D}(t) \mathrm{d} t \leq C_1 \mathcal{E}(0)
$$
for any $T\in(0, T_0]$.
\end{Proposition}

The following proposition is crucial for proving Theorem $\ref{Th8.1}$.
\begin{Proposition}  \label{pr8.1}         
There exists some positive constant $\delta$ such that if $(\boldsymbol{u}, \chi)$ is a solution of the problem $(\ref{1.1})$--$(\ref{1.2})$ on $\Omega \times[0, T]$, with $T >0$, satisfying
\begin{align}      \label{8.6}
\mathop {\mathrm{sup}} \limits_{0\le t\le T}\mathcal{E} (t)\le \delta,
\end{align}
then
\begin{align}      \label{8.7}
\mathop {\mathrm{sup}} \limits_{0\le t\le T}\mathcal{E} (t)+\int_0^T{\mathcal{D}}(t)\mathrm{d}t\le C_2\left(\varepsilon _0+\mathcal{E}(0)\right),
\end{align}
and there exists a positive constant $\sigma>0$ such that
\begin{align*}
\left\|{\boldsymbol{u}}(\cdot,t)\right\|_{L^2(\Omega)}^2
+\left\|\chi^2(\cdot,t)-1\right\|_{L^2(\Omega)}^2
+\left\|\nabla\chi(\cdot,t)\right\|_{L^2(\Omega)}^2
\le C{\mathcal E}(0)e^{-\sigma t},
\quad 0\le t\le T.
\end{align*}
\end{Proposition}
\noindent
{\bf Proof.} The proof consists of four steps.

\noindent
{\it\bfseries Step 1. Preliminary energy estimates and high-order energy inequality}\quad
Similar to $(\ref{4.19})$ and using $(\ref{8.3})$, it holds that
\begin{align}      \label{8.9}   
	\mathop {\mathrm{sup}} \limits_{0\le t\le T}\left( \left\| \boldsymbol{u} \right\| _{L^2(\Omega )}^{2}+\left\|  \nabla \chi \right\| _{L^2(\Omega )}^{2}+\left\| \chi ^2-1 \right\| _{L^2(\Omega )}^{2} \right) +\int_0^T{\left( \left\| \mathbb{D} \boldsymbol{u} \right\| _{L^2(\Omega )}^{2}+\left\|  \mu \right\| _{L^2(\Omega )}^{2} \right)}\mathrm{d}t \le C_0\varepsilon _0.
\end{align}
We can also deduce an inequality analogous to $(\ref{6.9})$, with the subscript $m$ omitted
\begin{align}     \label{8.10}   
\frac{\mathrm{d}}{\mathrm{d}t}\int_{\Omega}{\left( \frac{1}{2}\rho ^2\left| \partial _t\chi \right|^2+\frac{1}{2}|\mathbb{D} \boldsymbol{u}|^2 \right)}\mathrm{d}x+m_*\left\| \nabla \partial _t\chi \right\| _{L^2(\Omega )}^{2}+\rho _*\left\| \partial _t\boldsymbol{u} \right\| _{L^2(\Omega )}^{2}\le \sum_{i=1}^{15}{I_i}.
\end{align}
By the product rule, $|\chi|\leq1$, and the embedding
$H^1(\Omega)\hookrightarrow L^4(\Omega)$, we have
\begin{align*}
\|\chi^2-1\|_{H^2(\Omega)}
\leq C\left(
\|\chi^2-1\|_{L^2(\Omega)}
+\|\nabla\chi\|_{H^1(\Omega)}
+\|\nabla\chi\|_{L^4(\Omega)}^2
\right)
\leq C\left(\sqrt{\mathcal E(t)}+\mathcal E(t)\right).
\end{align*}
Since $H^2(\Omega)\hookrightarrow L^\infty(\Omega)$ in 3D, we may choose $\delta>0$ sufficiently small in \eqref{8.6} so that
$$
\|\chi^2-1\|_{L^\infty(\Omega)}
\leq C(\sqrt{\delta}+\delta)\leq\frac12.
$$
Consequently, $\chi^2\geq1/2$ and $F''(\chi)=3\chi^2-1\geq1/2$.
Hence, we have
$$
I_{12}=\int_{\Omega}{-}\rho m(\chi )F^{\prime\prime}(\chi )\left| \partial _t\chi \right|^2\mathrm{d}x\le0.
$$
Noticing that $\eta \equiv 1$, we have $I_{15}=\int_{\Omega}{\frac{1}{2}\eta ^{\prime}\left( \chi\right) \partial _t\chi\left| \mathbb{D} \boldsymbol{u}\right|^2}\mathrm{d}x=0$. Therefore, $(\ref{8.10})$ can be written in the following form
\begin{align}      \label{8.11}   
	\frac{\mathrm{d}}{\mathrm{d}t}&\int_{\Omega}{\left( \frac{1}{2}\rho ^2\left| \partial _t\chi \right|^2+\frac{1}{2}|\mathbb{D} \boldsymbol{u}|^2 \right)}\mathrm{d}x+m_*\left\| \nabla \partial _t\chi \right\| _{L^2(\Omega )}^{2}+\rho _*\left\| \partial _t\boldsymbol{u} \right\| _{L^2(\Omega )}^{2} \nonumber\\
	\le &\int_{\Omega}{-\rho ^2\partial _t\chi \boldsymbol{u}\cdot \nabla \partial _t\chi}\mathrm{d}x+\int_{\Omega}{-\rho ^2\left( \boldsymbol{u}\cdot \nabla \partial _t\chi \right) \left( \boldsymbol{u}\cdot \nabla \chi \right)}\mathrm{d}x +\int_{\Omega}{-\rho ^2\partial _{\mathrm{t}}\chi \boldsymbol{u}\cdot \nabla \left( \boldsymbol{u}\cdot \nabla \chi \right)}\mathrm{d}x \nonumber\\
	&+\int_{\Omega}{-\rho ^2\partial _t\chi \partial _t\boldsymbol{u}\cdot \nabla \chi}\mathrm{d}x+\int_{\Omega}{-\rho ^2\partial _t\chi \boldsymbol{u}\cdot \nabla \partial _t\chi}\mathrm{d}x+\int_{\Omega}{m^{\prime}(\chi )\left| \partial _t\chi \right|^2\Delta \chi}\mathrm{d}x \nonumber
\\
&+\int_{\Omega}{-m^{\prime}(\chi )\partial _t\chi \nabla \partial _t\chi \cdot \nabla \chi}\mathrm{d}x
+\int_{\Omega}{-\rho m^{\prime}(\chi )F^{\prime}(\chi )\partial _t\chi \boldsymbol{u}\cdot \nabla \chi}\mathrm{d}x
\\
&+\int_{\Omega}{-\rho m(\chi )F^{\prime\prime}(\chi )\partial _t\chi \boldsymbol{u}\cdot \nabla \chi}\mathrm{d}x \nonumber
+\int_{\Omega}{-\rho m(\chi )F^{\prime}(\chi )\boldsymbol{u}\cdot \nabla \partial _t\chi}\mathrm{d}x
\\
&+\int_{\Omega}{-\rho m^{\prime}(\chi )\left| \partial _t\chi \right|^2F^{\prime}(\chi )}\mathrm{d}x
+\int_{\Omega}{-\mathrm{div}\left( \nabla \chi \otimes \nabla \chi \right) \cdot \partial _t\boldsymbol{u}}\mathrm{d}x \nonumber
+\int_{\Omega}{-\rho \left( \boldsymbol{u}\cdot \nabla \right) \boldsymbol{u}\cdot \partial _t\boldsymbol{u}}\mathrm{d}x.
\end{align}

From $(\ref{1.1})_{4,5}$, we have
\begin{align}     \label{8.12}   
\rho ^2\partial _t\chi +\rho ^2\boldsymbol{u}\cdot \nabla \chi =m(\chi )\Delta \chi -\rho m(\chi )F^{\prime}(\chi ).
\end{align}
Multiplying the above equality by $2\chi$ yields
\begin{align}     \label{8.13}   
	&\rho ^2\partial _t\left( \chi ^2-1 \right) +\rho ^2\boldsymbol{u}\cdot \nabla \left( \chi ^2-1 \right) \nonumber\\
	&=m(\chi )\Delta \left( \chi ^2-1 \right) -2m(\chi )|\nabla \chi |^2-2\rho m(\chi )\left( \chi ^2-1 \right) ^2-2\rho m(\chi )\left( \chi ^2-1 \right).
\end{align}
Then multiplying $(\ref{8.13})$ by $\chi^2-1$ and integrating the result over $\Omega$, we get
\begin{align}     \label{8.14}   
	&\frac{\mathrm{d}}{\mathrm{d}t}\int_{\Omega}{\frac{1}{2}}\rho ^2\left| \chi ^2-1 \right|^2\mathrm{d}x+\int_{\Omega}{2}\rho m(\chi )\left| \chi ^2-1 \right|^2\mathrm{d}x+\int_{\Omega}{m}(\chi )\left| \nabla \left( \chi ^2-1 \right) \right|^2\mathrm{d}x \nonumber
\\
	=&-\int_{\Omega}{m^{\prime}}(\chi )\left( \chi ^2-1 \right) \nabla \chi \cdot \nabla \left( \chi ^2-1 \right) \mathrm{d}x
-\int_{\Omega}{2}\rho m(\chi )\left( \chi ^2-1 \right) ^3\mathrm{d}x  \nonumber
\\
&-\int_{\Omega}{2}m(\chi )|\nabla \chi |^2\left( \chi ^2-1 \right) \mathrm{d}x.
\end{align}
By adding $(\ref{8.11})$ and $(\ref{8.14})$ together, we obtain
\begin{align}      \label{8.15}   
\frac{\mathrm{d}}{\mathrm{d}t}&\int_{\Omega}{\left( \frac{1}{2}\rho ^2\left| \partial _t\chi \right|^2+\frac{1}{2}|\mathbb{D} \boldsymbol{u}|^2+\frac{1}{2}\rho ^2\left| \chi ^2-1 \right|^2 \right)}\mathrm{d}x+m_*\left\| \nabla \partial _t\chi \right\| _{L^2(\Omega )}^{2}+\rho _*\left\| \partial _t\boldsymbol{u} \right\| _{L^2(\Omega )}^{2}\nonumber
\\
&\quad +2\rho _*m_*\left\| \chi ^2-1 \right\| _{L^2(\Omega )}^{2}+m_*\left\| \nabla \left( \chi ^2-1 \right) \right\| _{L^2(\Omega )}^{2}\nonumber
\\
\le & {\int_{\Omega}{-\rho ^{2} \partial _t\chi  \boldsymbol{u} \cdot \nabla \partial _t\chi }\mathrm{d}x}
+{\int_{\Omega}{-\rho ^{2}\left( \boldsymbol{u}\cdot \nabla \partial _t\chi  \right) \left( \boldsymbol{u}\cdot \nabla \chi  \right)}\mathrm{d}x} \nonumber
\\
&+{\int_{\Omega}{-\rho ^{2}\partial _{\mathrm{t}}\chi \boldsymbol{u}\cdot \nabla \left( \boldsymbol{u}\cdot \nabla \chi  \right)}\mathrm{d}x}
+{\int_{\Omega}{-\rho ^{2}\partial _t\chi \partial _t\boldsymbol{u}\cdot \nabla \chi }\mathrm{d}x}\nonumber
\\
&+{\int_{\Omega}{-\rho ^{2}\partial _t\chi \boldsymbol{u} \cdot\nabla \partial _t\chi }\mathrm{d}x}
+{\int_{\Omega}{m^{\prime}(\chi )\left| \partial _t\chi  \right|^2\Delta \chi }\mathrm{d}x} \nonumber
\\
&+{\int_{\Omega}{-m^{\prime}(\chi )\partial _t\chi \nabla \partial _t\chi \cdot \nabla \chi }\mathrm{d}x}
+{\int_{\Omega}{-\rho m^{\prime}(\chi )F^{\prime}(\chi )\partial _t\chi \boldsymbol{u}\cdot\nabla \chi }\mathrm{d}x}
\\
&+{\int_{\Omega}{-\rho m(\chi )F^{\prime\prime}(\chi )\partial _t\chi \boldsymbol{u}\cdot\nabla \chi }\mathrm{d}x}
+{\int_{\Omega}{-\rho m(\chi )F^{\prime}(\chi )\boldsymbol{u}\cdot\nabla \partial _t\chi }\mathrm{d}x} \nonumber
\\
&+{\int_{\Omega}{-\rho m^{\prime}(\chi )\left| \partial _t\chi  \right|^2F^{\prime}(\chi )}\mathrm{d}x}
+{\int_{\Omega}{-\mathrm{div}\left( \nabla \chi  \otimes \nabla \chi  \right)\cdot \partial _t\boldsymbol{u}}\mathrm{d}x} \nonumber
\\
&+{\int_{\Omega}{-\rho \left( \boldsymbol{u}\cdot \nabla \right) \boldsymbol{u}\cdot \partial _t\boldsymbol{u}}\mathrm{d}x}
+{\int_{\Omega}{-m^{\prime}(\chi )\left( \chi ^2-1 \right) \nabla \chi \cdot \nabla \left( \chi ^2-1 \right)}\mathrm{d}x} \nonumber
\\
&+{\int_{\Omega}{-2m(\chi )|\nabla \chi |^2\left( \chi ^2-1 \right)}\mathrm{d}x}
+{\int_{\Omega}{-2\rho m(\chi )\left( \chi ^2-1 \right) ^3}\mathrm{d}x}=\sum_{i=1}^{16}J_i.  \nonumber
\end{align}

\noindent
{\it\bfseries Step 2. Absorption of the nonlinear terms}\quad
{We estimate the terms on the right-hand side of \eqref{8.15} using H\"older's inequality, the Sobolev
embeddings, and the interpolation inequalities stated below.
Throughout this step, $C$ is independent of $T$ and $\delta\leq1$.}

{The homogeneous Dirichlet condition for $\boldsymbol{u}$ and the homogeneous Neumann condition for $\chi$ imply
$$
\|\boldsymbol{u}\|_{L^6(\Omega)}
\leq C\|\nabla\boldsymbol{u}\|_{L^2(\Omega)},
\qquad
\|\nabla\chi\|_{H^1(\Omega)}
\leq C\|\Delta\chi\|_{L^2(\Omega)}.
$$
Here the second estimate follows by applying the Neumann elliptic estimate to $\chi$ minus its spatial mean. We also use
\begin{align*}
\|f\|_{L^3(\Omega)}
&\leq C\|f\|_{L^2(\Omega)}^{1/2}
          \|f\|_{H^1(\Omega)}^{1/2},\\
\|\nabla\chi\|_{L^4(\Omega)}^2
&\leq C\|\nabla\chi\|_{L^2(\Omega)}^{1/2}
          \|\Delta\chi\|_{L^2(\Omega)}^{3/2}.
\end{align*}
By \eqref{8.1} and \eqref{8.2}, the norms entering $\mathcal E(t)$ and $\mathcal D(t)$ are bounded by
$\sqrt{\mathcal E(t)}$ and $\sqrt{\mathcal D(t)}$, respectively. In particular, $\|\Delta\chi\|_{L^2(\Omega)}$
has both bounds, up to a fixed constant.
After extracting a factor $\sqrt{\mathcal E(t)}$, we use $ab\leq(a^2+b^2)/2$ to control the remaining factors by $\mathcal D(t)$.
Additional energy factors are controlled using $\mathcal E(t)\leq\delta\leq1$.}

{Since $\rho$ is bounded above and below, $|\chi|\leq1$, and $m,m',F''$ are bounded on $[-1,1]$, we obtain the following estimates.}
\begin{align}      \label{8.16}   
~~~~~~~~~~~~~~~~~~~ \left| J_1 \right|&\le C\left\| \partial _t\chi \right\| _{L^3(\Omega )}\left\| \boldsymbol{u}\right\| _{L^6(\Omega )}\left\| \nabla \partial _t\chi \right\| _{L^2(\Omega )} \nonumber
\\
&\le C\left\| \nabla \boldsymbol{u}\right\| _{L^2(\Omega )}\left\| \partial _t\chi \right\| _{H^1(\Omega )}^{2}
\le C\sqrt{\mathcal{E} (t)}\mathcal{D} (t),&
\end{align}
\begin{align}      \label{8.17}   
~~~~~~~~~~~~~~~~~~~  \left| J_2 \right|&\le C\left\| \boldsymbol{u}\right\| _{L^6(\Omega )}^{2}\left\| \nabla \partial _t\chi \right\| _{L^2(\Omega )}\left\| \nabla \chi \right\| _{L^6(\Omega )} \nonumber
\\
&\le C\left\| \nabla \boldsymbol{u}\right\| _{L^2(\Omega )}^{2}\left\| \nabla \partial _t\chi \right\| _{L^2\left( \Omega \right)}\left\| \Delta \chi \right\| _{L^2\left( \Omega \right)}
\le C\sqrt{\mathcal{E} (t)}\mathcal{D} (t),&
\end{align}
{
For $J_3$, expanding
\[
\nabla(\boldsymbol{u}\cdot\nabla\chi)
=(\nabla\boldsymbol{u})^T\nabla\chi
+(\boldsymbol{u}\cdot\nabla)\nabla\chi,
\]
and applying H\"older's inequality, we obtain}
\begin{align}\label{8.18}
|J_3|
&\leq
C\|\partial_t\chi\|_{L^6(\Omega)}
 \|\boldsymbol{u}\|_{L^6(\Omega)}
 \|\nabla\boldsymbol{u}\|_{L^2(\Omega)}
 \|\nabla\chi\|_{L^6(\Omega)}
\nonumber\\
&\quad{ +
C\|\partial_t\chi\|_{L^6(\Omega)}
 \|\boldsymbol{u}\|_{L^6(\Omega)}^2
 \|\nabla^2\chi\|_{L^2(\Omega)}}
\nonumber\\
&\leq
C\|\partial_t\chi\|_{H^1(\Omega)}
 { \|\boldsymbol{u}\|_{H^1(\Omega)}^2
 \|\nabla\chi\|_{H^1(\Omega)}}
\nonumber\\
&{\leq C\mathcal E(t)\mathcal D(t)}
\leq C\sqrt{\mathcal E(t)}\mathcal D(t).
\end{align}
{The terms $J_4,\ldots,J_9$ are estimated directly
using the bounds stated at the beginning of this step:}
\begin{align}      \label{8.19}   
~~~~~~~~~~~~~~~~~~~ 	\left| J_4 \right|&\le C\left\| \partial _t\chi \right\| _{L^6(\Omega )}\left\| \partial _t\boldsymbol{u} \right\| _{L^2(\Omega )}\left\|\nabla \chi \right\| _{L^3(\Omega )} \nonumber
\\
	&\le C\left\| \partial _t\chi \right\| _{H^1(\Omega )}\left\| \partial _t\boldsymbol{u} \right\| _{L^2(\Omega )}\left\| \nabla \chi \right\| _{H^1(\Omega )}
\le C\sqrt{\mathcal{E} (t)}\mathcal{D} (t),&
\end{align}
\begin{align}      \label{8.20}   
~~~~~~~~~~~~~~~~~~~ 	\left| J_5 \right|&\le C\left\| \partial _t\chi \right\| _{L^3(\Omega )}
\left\| \boldsymbol{u}\right\| _{L^6(\Omega )}
\left\| \nabla \partial _t\chi \right\| _{L^2(\Omega )} \nonumber
\\
&\le C\left\| \partial _t\chi \right\| _{H^1(\Omega )}
\left\| \nabla \boldsymbol{u}\right\| _{L^2(\Omega )}
\left\| \nabla \partial _t\chi \right\| _{L^2(\Omega )}
\le C\sqrt{\mathcal{E} (t)}\mathcal{D} (t),&
\end{align}
\begin{align}      \label{8.21}   
~~~~~~~~~~~~~~~~~~~  \left| J_6 \right|&\le C\left\| \partial _t\chi \right\| _{L^3(\Omega )}
\left\| \partial _t\chi \right\| _{L^6(\Omega )}
\left\| \Delta \chi \right\| _{L^2(\Omega )} \nonumber
\\
&\le C\left\| \partial _t\chi \right\| _{H^1(\Omega )}^{2}\left\| \Delta \chi \right\| _{L^2(\Omega )}
\le C\sqrt{\mathcal{E} (t)}\mathcal{D} (t),&
\end{align}
\begin{align}      \label{8.22}   
~~~~~~~~~~~~~~~~~~~ 	\left| J_7 \right|&\le C\left\| \partial _t\chi \right\| _{L^3(\Omega )}\left\| \nabla \partial _t\chi \right\| _{L^2(\Omega )}\left\| \nabla \chi \right\| _{L^6(\Omega )} \nonumber
\\
&\le C\left\| \partial _t\chi \right\| _{H^1(\Omega )}
\left\| \nabla \partial _t\chi \right\| _{L^2(\Omega )}
\left\| \Delta \chi \right\| _{L^2(\Omega )}
\le C\sqrt{\mathcal{E} (t)}\mathcal{D} (t),&
\end{align}
\begin{align}      \label{8.23}   
~~~~~~~~~~~~		\left| J_8+J_9 \right|
&\le C\left\| \partial _t\chi \right\| _{L^6(\Omega )}
\left\| \boldsymbol{u} \right\| _{L^3(\Omega )}
\left\| \nabla \chi \right\| _{L^2(\Omega )} \nonumber
\\
&\le C\left\| \partial _t\chi \right\| _{H^1(\Omega )}
\left\| \boldsymbol{u} \right\| _{H^1(\Omega )}
\left\| \nabla \chi \right\| _{L^2(\Omega )}
\le C\sqrt{\mathcal{E} (t)}\mathcal{D} (t),&
\end{align}
{For $J_{10}$, we use
$|F'(\chi)|\leq|\chi^2-1|$ and apply the interpolation
inequality above to $\chi^2-1$, }
\begin{align}      \label{8.24}   
~~~~~~~~~~~~~~~~~~~  \left| J_{10} \right|
&\le C\left\| \chi ^2-1 \right\| _{L^3(\Omega )}
\left\| \boldsymbol{u} \right\| _{L^6(\Omega )}
\left\| \nabla \partial _t\chi \right\| _{L^2(\Omega )} \nonumber
\\
&\le C\left( \left\| 2\chi \nabla \chi \right\| _{L^2(\Omega )}^{\frac{1}{2}}\left\| \chi ^2-1 \right\| _{L^2(\Omega )}^{\frac{1}{2}}+\left\| \chi ^2-1 \right\| _{L^2(\Omega )} \right) \left\| \nabla \boldsymbol{u} \right\| _{L^2(\Omega )}\left\| \nabla \partial _t\chi \right\| _{L^2(\Omega )} \nonumber\\
	&\le C\left( \left\| \nabla \chi \right\| _{L^2(\Omega )}+\left\| \chi ^2-1 \right\| _{L^2(\Omega )} \right) \left\| \nabla \boldsymbol{u} \right\| _{L^2(\Omega )}\left\| \nabla \partial _t\chi \right\| _{L^2(\Omega )}\\
	&\le C\sqrt{\mathcal{E} (t)}\mathcal{D} (t),& \nonumber
\end{align}
\begin{align}      \label{8.25}   
~~~~~~~~~~~~~~~~~~~  	\left| J_{11} \right|
&\le C\left\| \chi ^2-1 \right\| _{L^2(\Omega )}
\left\| \partial _t\chi \right\| _{L^3(\Omega )}
\left\| \partial _t\chi \right\| _{L^6(\Omega )} \nonumber
\\
&\le C\left\| \chi ^2-1 \right\| _{L^2(\Omega )}
\left\| \partial _t\chi \right\| _{H^1(\Omega )}^{2}
\le C\sqrt{\mathcal{E} (t)}\mathcal{D} (t).&
\end{align}
From the equation $(\ref{8.12})$, we see that
\begin{align}      \label{8.26}   
	\left\| \Delta \chi \nabla \chi \right\| _{L^2(\Omega )}\le &C\left( \left\| \partial _t\chi \nabla \chi \right\| _{L^2(\Omega )}+\left\| \boldsymbol{u}\cdot \nabla \chi \nabla \chi \right\| _{L^2(\Omega )}+\left\| \chi \left( \chi ^2-1 \right) \nabla \chi \right\| _{L^2(\Omega )} \right) \nonumber\\
	\le &C\left( \left\| \partial _t\chi \right\| _{L^6(\Omega )}\left\| \nabla \chi \right\| _{L^3(\Omega )}+\left\| \boldsymbol{u} \right\| _{L^6(\Omega )}\left\| \nabla \chi \right\| _{L^6\left( \Omega \right)}^{2}+\left\| \chi ^2-1 \right\| _{L^6\left( \Omega \right)}\left\| \nabla \chi \right\| _{L^3(\Omega )} \right) \nonumber\\
	\le &C\left\| \partial _t\chi \right\| _{H^1(\Omega )}\left\| \nabla \chi \right\| _{H^1(\Omega )}+C\left\| \nabla \boldsymbol{u} \right\| _{L^2\left( \Omega \right)}\left\| \Delta \chi \right\| _{L^2\left( \Omega \right)}^{2} \nonumber\\
	&+C\left( \left\| \nabla \chi \right\| _{L^2\left( \Omega \right)}+\left\| \chi ^2-1 \right\| _{L^2(\Omega )} \right) \left\| \nabla \chi \right\| _{H^1(\Omega )} \nonumber\\
    \le& C\sqrt{\mathcal{E} (t)}\sqrt{\mathcal{D} (t)}.
\end{align}
{For $J_{12}$, we use
\[
\operatorname{div}(\nabla\chi\otimes\nabla\chi)
=\Delta\chi\nabla\chi+\nabla\left(\frac12|\nabla\chi|^2\right).
\]
The gradient term has zero pairing with
$\partial_t\boldsymbol{u}$, which belongs to the solenoidal
$L^2$ space and is orthogonal to gradients.
This cancellation is justified first at the smooth
approximation level and then passed to the limit.
Consequently, \eqref{8.26} yields}
\begin{align}\label{8.27}
|J_{12}|
&{=\left|\int_\Omega
\Delta\chi\nabla\chi\cdot\partial_t\boldsymbol{u}
\,\mathrm{d}x\right|}
\leq
\|\partial_t\boldsymbol{u}\|_{L^2(\Omega)}
\|\Delta\chi\nabla\chi\|_{L^2(\Omega)}
\leq C\sqrt{\mathcal E(t)}\mathcal D(t).
\end{align}  
{ To estimate $J_{13}$, we first derive an $H^2$ estimate
for the velocity.
Since $\eta\equiv1$, the momentum equation becomes
\[
-\frac12\Delta\boldsymbol{u}+\nabla\widetilde p
=-\rho\partial_t\boldsymbol{u}
-\rho\boldsymbol{u}\cdot\nabla\boldsymbol{u}
-\Delta\chi\nabla\chi,
\qquad
\widetilde p=p+\frac12|\nabla\chi|^2.
\]
We normalize $\widetilde p$ to have zero mean.
The elliptic regularity theory for the stationary Stokes problem \cite{T-77}
and H\"older's inequality give}
\begin{align*}
\|\boldsymbol{u}\|_{H^2(\Omega)}
+\|\widetilde p\|_{H^1(\Omega)}
&\leq C\left(
\|\partial_t\boldsymbol{u}\|_{L^2(\Omega)}
+\|\boldsymbol{u}\|_{L^6(\Omega)}
 \|\nabla\boldsymbol{u}\|_{L^3(\Omega)}
+\|\Delta\chi\nabla\chi\|_{L^2(\Omega)}
\right)\\
&\leq C\left(
\|\partial_t\boldsymbol{u}\|_{L^2(\Omega)}
+\|\nabla\boldsymbol{u}\|_{L^2(\Omega)}^{3/2}
 \|\boldsymbol{u}\|_{H^2(\Omega)}^{1/2}
+\|\Delta\chi\nabla\chi\|_{L^2(\Omega)}
\right).
\end{align*}
{ Here we used the interpolation inequality
\[
\|\nabla\boldsymbol{u}\|_{L^3(\Omega)}
\leq C\|\nabla\boldsymbol{u}\|_{L^2(\Omega)}^{1/2}
          \|\boldsymbol{u}\|_{H^2(\Omega)}^{1/2}.
\]
Young's inequality allows us to absorb one half of
$\|\boldsymbol{u}\|_{H^2(\Omega)}$, yielding
\begin{align*}
\|\boldsymbol{u}\|_{H^2(\Omega)}
&\leq C\left(
\|\partial_t\boldsymbol{u}\|_{L^2(\Omega)}
+\|\nabla\boldsymbol{u}\|_{L^2(\Omega)}^3
+\|\Delta\chi\nabla\chi\|_{L^2(\Omega)}
\right)\\
&\leq C\left(
1+\mathcal E(t)+\sqrt{\mathcal E(t)}
\right)\sqrt{\mathcal D(t)}
\leq C\sqrt{\mathcal D(t)},
\end{align*}
where we used \eqref{8.26} and $\mathcal E(t)\leq1$.
Consequently, }
\begin{align}      \label{8.28}   
~~~~~~~~~~~~~~~~~~~~~~~~~~~	\left| J_{13} \right|
&\le C\left\| \boldsymbol{u} \right\| _{L^6(\Omega )}\left\| \nabla \boldsymbol{u} \right\| _{L^3\left( \Omega \right)}\left\| \partial _t\boldsymbol{u} \right\| _{L^2\left( \Omega \right)}\nonumber
\\
&\le C\left\| \nabla \boldsymbol{u} \right\| _{L^2\left( \Omega \right)}
\left\| \nabla \boldsymbol{u} \right\| _{H^1\left( \Omega \right)}
\left\| \partial _t\boldsymbol{u} \right\| _{L^2(\Omega )}
\le C\sqrt{\mathcal{E} (t)}\mathcal{D} (t),&
\end{align}
{For $J_{14}+J_{15}$, we substitute
$\nabla(\chi^2-1)=2\chi\nabla\chi$ and use the
$L^4$ interpolation inequality above.
For $J_{16}$, we use
$H^2(\Omega)\hookrightarrow L^\infty(\Omega)$
and the product estimate for $\chi^2-1$ established in Step 1.}
\begin{align}      \label{8.29}   
~~~~~~~~~~~~~~~~~~	\left| J_{14}+J_{15} \right|
&\le C\left\| \nabla \chi \right\| _{L^4\left( \Omega \right)}^{2}\left\| \chi ^2-1 \right\| _{L^2(\Omega )} \nonumber
\\
&\le C\left\| \nabla \chi \right\| _{L^2\left( \Omega \right)}^{\frac{1}{2}}
\left\| \Delta \chi \right\| _{L^2\left( \Omega \right)}^{\frac{3}{2}}
\left\| \chi ^2-1 \right\| _{L^2(\Omega )}
\le C\sqrt{\mathcal{E} (t)}\mathcal{D} (t),&
\end{align}
\begin{align}      \label{8.30}   
~~~~~~~~~~~~~~~~~~~~~~~~~~~	\left| J_{16} \right|&\le C\left\| \chi ^2-1 \right\| _{L^{\infty}(\Omega )}\left\| \chi ^2-1 \right\| _{L^2\left( \Omega \right)}^{2} \nonumber\\
	&\le C\left\| \chi ^2-1 \right\| _{H^2\left( \Omega \right)}\left\| \chi ^2-1 \right\| _{L^2\left( \Omega \right)}^{2}\\
	&\le C\left( \left\| \chi ^2-1 \right\| _{L^2(\Omega )}+\left\| \nabla \chi \right\| _{H^1\left( \Omega \right)} \right) \left\| \chi ^2-1 \right\| _{L^2\left( \Omega \right)}^{2} \nonumber\\
	&\le C\sqrt{\mathcal{E} (t)}\mathcal{D} (t).& \nonumber
\end{align}

\vskip2mm
\noindent
{\it\bfseries Step 3. Closing the energy estimates}\quad
Combining $(\ref{8.15})$-$(\ref{8.30})$, we arrive at
\begin{align}      \label{8.31}   
  \frac{\mathrm{d}}{\mathrm{d}t}&\int_{\Omega}{\left( \frac{1}{2}\rho ^2\left| \partial _t\chi \right|^2+\frac{1}{2}|\mathbb{D} \boldsymbol{u}|^2+\frac{1}{2}\rho ^2\left| \chi ^2-1 \right|^2 \right)}\mathrm{d}x+m_*\left\| \nabla \partial _t\chi \right\| _{L^2(\Omega )}^{2}+\rho _*\left\| \partial _t\boldsymbol{u} \right\| _{L^2(\Omega )}^{2} \nonumber\\
  &\quad +2\rho _*m_*\left\| \chi ^2-1 \right\| _{L^2(\Omega )}^{2}+m_*\left\| \nabla \left( \chi ^2-1 \right) \right\| _{L^2(\Omega )}^{2}
  \le C\sqrt{\mathcal{E} (t)}\mathcal{D} (t).
\end{align}
{Integrating \eqref{8.31} over $(0,t)$, retaining the
initial higher-order energy, and using the density bounds
and Korn's inequality, we obtain bounds for
$\partial_t\chi$ and $\mathbb D\boldsymbol{u}$.
Combining these bounds with the lower-order estimate
\eqref{8.9}, and then taking the supremum over
$0\leq t\leq T$, gives}
\begin{align}      \label{8.32}   
	&\mathop {\mathrm{sup}} \limits_{0\le t\le T}\left( \left\| \boldsymbol{u} \right\| _{L^2\left( \Omega \right)}^{2}+\left\| \nabla \chi \right\| _{L^2\left( \Omega \right)}^{2}+\left\| \chi ^2-1 \right\| _{L^2\left( \Omega \right)}^{2}+\left\| \partial _t\chi \right\| _{L^2\left( \Omega \right)}^{2}+\left\| \mathbb{D} \boldsymbol{u} \right\| _{L^2\left( \Omega \right)}^{2} \right) \nonumber
\\
	&~+\int_0^T\Big( \left\| \mathbb{D} \boldsymbol{u} \right\| _{L^2\left( \Omega \right)}^{2}
+\left\| \mu \right\| _{L^2\left( \Omega \right)}^{2}
+\left\| \partial _t\boldsymbol{u} \right\| _{L^2\left( \Omega \right)}^{2}
+\left\| \nabla \partial _t\chi \right\| _{L^2\left( \Omega \right)}^{2}  \nonumber
\\
&\hskip17mm+\left\| \chi ^2-1 \right\| _{L^2\left( \Omega \right)}^{2}
+\left\| \nabla \left( \chi ^2-1 \right) \right\| _{L^2\left( \Omega \right)}^{2} \Big)\mathrm{d}t \nonumber
\\
	&\le C\left(\varepsilon _0+\mathcal{E}(0)\right)+C\int_0^T{\sqrt{\mathcal{E} (t)}\mathcal{D} (t)}\mathrm{d}t.
\end{align}
{ To recover the $H^1$ norm of $\nabla\chi$ in
$\mathcal E(t)$, we solve \eqref{8.12} for $\Delta\chi$.
The convection term is estimated using
\[
\|\boldsymbol{u}\cdot\nabla\chi\|_{L^2(\Omega)}^2
\leq C\|\nabla\boldsymbol{u}\|_{L^2(\Omega)}^2
       \|\nabla\chi\|_{L^2(\Omega)}
       \|\Delta\chi\|_{L^2(\Omega)}.
\]
Using this bound in \eqref{8.12} and then applying
Young's inequality, we obtain }
\begin{align*}
	\left\| \Delta \chi \right\| _{L^2\left( \Omega \right)}^{2}\le &C\left( \left\| \partial _t\chi \right\| _{L^2\left( \Omega \right)}^{2}+\left\| \boldsymbol{u}\cdot \nabla \chi \right\| _{L^2\left( \Omega \right)}^{2}+\left\| \chi ^3-\chi \right\| _{L^2\left( \Omega \right)}^{2} \right)\\
	\le &C\left( \left\| \partial _t\chi \right\| _{L^2\left( \Omega \right)}^{2}+\left\| \boldsymbol{u} \right\| _{L^6\left( \Omega \right)}^{2}\left\| \nabla \chi \right\| _{L^3\left( \Omega \right)}^{2}+\left\| \chi ^2-1 \right\| _{L^2\left( \Omega \right)}^{2} \right)\\
	\le &C\left\| \partial _t\chi \right\| _{L^2\left( \Omega \right)}^{2}+C\left\| \nabla \boldsymbol{u} \right\| _{L^2\left( \Omega \right)}^{2}\left\| \nabla \chi \right\| _{L^2\left( \Omega \right)}\left\| \Delta \chi \right\| _{L^2\left( \Omega \right)}+C\left\| \chi ^2-1 \right\| _{L^2\left( \Omega \right)}^{2}\\
	\le & \frac{1}{2}\left\| \Delta \chi \right\| _{L^2\left( \Omega \right)}^{2}+C\left\| \partial _t\chi \right\| _{L^2\left( \Omega \right)}^{2}+C\left\| \nabla \boldsymbol{u} \right\| _{L^2\left( \Omega \right)}^{4}\left\| \nabla \chi \right\| _{L^2\left( \Omega \right)}^{2}+C\left\| \chi ^2-1 \right\| _{L^2\left( \Omega \right)}^{2},
\end{align*}
i.e.
\begin{align}      \label{8.33}   
	\left\| \Delta \chi \right\| _{L^2\left( \Omega \right)}^{2}\le &C\left( \left\| \partial _t\chi \right\| _{L^2\left( \Omega \right)}^{2}+\left\| \nabla \boldsymbol{u} \right\| _{L^2\left( \Omega \right)}^{4}\left\| \nabla \chi \right\| _{L^2\left( \Omega \right)}^{2}+\left\| \chi ^2-1 \right\| _{L^2\left( \Omega \right)}^{2} \right) \nonumber\\
	\le &C\left( \left\| \partial _t\chi \right\| _{L^2\left( \Omega \right)}^{2}+\left\| \nabla \chi \right\| _{L^2\left( \Omega \right)}^{2}+\left\| \chi ^2-1 \right\| _{L^2\left( \Omega \right)}^{2} \right)\\
	\le &C { \left(\varepsilon_0+\mathcal E(0)\right) } +C\int_0^T{\sqrt{\mathcal{E} (t)}\mathcal{D} (t)}\mathrm{d}t. \nonumber
\end{align}
{
For the time-integrated estimates, we use the
chemical-potential relation and the phase equation:
\begin{align*}
\|\Delta\chi\|_{L^2(\Omega)}^2
&\leq C\left(
\|\mu\|_{L^2(\Omega)}^2
+\|\chi^2-1\|_{L^2(\Omega)}^2
\right),\\
\|\partial_t\chi\|_{L^2(\Omega)}^2
&\leq C\|\mu\|_{L^2(\Omega)}^2
+C\|\boldsymbol{u}\cdot\nabla\chi\|_{L^2(\Omega)}^2\\
&\leq C\|\mu\|_{L^2(\Omega)}^2
+C\|\nabla\boldsymbol{u}\|_{L^2(\Omega)}^2
  \|\nabla\chi\|_{L^2(\Omega)}
  \|\Delta\chi\|_{L^2(\Omega)}\\
&\leq C\|\mu\|_{L^2(\Omega)}^2
+C\delta\|\nabla\boldsymbol{u}\|_{L^2(\Omega)}^2.
\end{align*}
The last inequality follows from \eqref{8.6}.
Integrating these estimates and using \eqref{8.9}
and \eqref{8.32}, we obtain}
\begin{align}      \label{8.34}   
	&\int_0^T{\left( \left\| \Delta \chi \right\| _{L^2\left( \Omega \right)}^{2}+\left\| \partial _t\chi \right\| _{L^2\left( \Omega \right)}^{2} \right)}\mathrm{d}t \nonumber\\
	&\le C\int_0^T{\left( \left\| \mu \right\| _{L^2\left( \Omega \right)}^{2}+\left\| \chi ^2-1 \right\| _{L^2\left( \Omega \right)}^{2}+\left\| \nabla \boldsymbol{u} \right\| _{L^2\left( \Omega \right)}^{2} \right)}\mathrm{d}t \nonumber\\
	&\le C{ \left(\varepsilon_0+\mathcal E(0)\right) }+C\int_0^T{\sqrt{\mathcal{E} (t)}\mathcal{D} (t)}\mathrm{d}t.
\end{align}
{ The Neumann elliptic estimate controls
$\|\nabla\chi\|_{H^1(\Omega)}$ by
$\|\Delta\chi\|_{L^2(\Omega)}$, while Poincar\'e's
and Korn's inequalities control
$\|\boldsymbol{u}\|_{H^1(\Omega)}$ by
$\|\mathbb D\boldsymbol{u}\|_{L^2(\Omega)}$.
Thus \eqref{8.32}--\eqref{8.34} control every term
in $\mathcal E(t)$ and $\mathcal D(t)$.
Using $\mathcal E(t)\leq\delta$, we obtain }
\begin{align}      \label{8.35}   
\mathop {\mathrm{sup}} \limits_{0\le t\le T}\mathcal{E} (t)+\int_0^T{\mathcal{D} (t)}\mathrm{d}t
&\le C{ \left(\varepsilon_0+\mathcal E(0)\right) }+C\int_0^T{\sqrt{\mathcal{E} (t)}\mathcal{D} (t)}\mathrm{d}t \nonumber
\\
&\le C_2\left(\varepsilon _0+\mathcal{E}(0)\right)+C_2\delta ^{\frac{1}{2}}\int_0^T{\mathcal{D} (t)}\mathrm{d}t.
\end{align}
{ Choose $\delta$ sufficiently small that $C_2\delta ^{\frac{1}{2}}<\frac{1}{2}$.
Absorbing this term into the left-hand side and enlarging
the final constant $C_2$, we obtain \eqref{8.7}.
The constants and the smallness threshold $\delta$ are independent
of $T$. }

\vskip2mm
\noindent
{\it\bfseries Step 4. Exponential decay to equilibrium}\quad
{ Testing the momentum equation by $\boldsymbol{u}$ and
the phase equation by $\mu$, and using the continuity
equation and the boundary conditions, we obtain}
\begin{equation}\label{energy}
\frac{\mathrm d}{\mathrm dt}\int_\Omega
\left(
\frac12\rho|\boldsymbol{u}|^2
+\frac12|\nabla\chi|^2
+\frac14\rho|\chi^2-1|^2
\right)\,\mathrm dx
+\|\mathbb D\boldsymbol{u}\|_{L^2(\Omega)}^2
+\int_\Omega m(\chi)|\mu|^2\,\mathrm dx=0.
\end{equation}
{Here the transport terms cancel by the continuity
equation, and the capillary term cancels the corresponding
phase-gradient contribution. }

{ Adding \eqref{8.14} to \eqref{energy}, we set}
\[
\mathcal A(t)=\int_\Omega
\left(
\frac12\rho|\boldsymbol{u}|^2
+\frac12|\nabla\chi|^2
+\frac14\rho|\chi^2-1|^2
+\frac12\rho^2|\chi^2-1|^2
\right)\,\mathrm dx.
\]
{By $\nabla(\chi^2-1)=2\chi\nabla\chi$, $|\chi|\leq1$,
and the bounds on $\rho,m,m'$, the three terms on the
right-hand side of \eqref{8.14} satisfy
\begin{align*}
&\left|\int_\Omega
m'(\chi)(\chi^2-1)\nabla\chi\cdot
\nabla(\chi^2-1)\,\mathrm dx\right|
\\
&\quad+
\left|\int_\Omega
2m(\chi)|\nabla\chi|^2(\chi^2-1)\,\mathrm dx\right|
+
\left|\int_\Omega
2\rho m(\chi)(\chi^2-1)^3\,\mathrm dx\right|
\\
&\leq
C\|\chi^2-1\|_{L^\infty(\Omega)}
\left(
\|\nabla\chi\|_{L^2(\Omega)}^2
+\|\chi^2-1\|_{L^2(\Omega)}^2
\right).
\end{align*}
The product estimate in Step 1 gives
$\|\chi^2-1\|_{L^\infty(\Omega)}\leq C\sqrt{\delta}$.
Moreover, the Neumann elliptic estimate and the
chemical-potential relation imply
\[
\|\nabla\chi\|_{L^2(\Omega)}^2
\leq C\|\Delta\chi\|_{L^2(\Omega)}^2
\leq C\left(
\|\mu\|_{L^2(\Omega)}^2
+\|\chi^2-1\|_{L^2(\Omega)}^2
\right).
\]
Consequently, the preceding three terms are bounded by
\[
C\sqrt{\delta}\left(
\|\mu\|_{L^2(\Omega)}^2
+\|\chi^2-1\|_{L^2(\Omega)}^2
\right).
\]
Reducing $\delta$ further so that
$C\sqrt{\delta}\leq\min\{m_*/2,\rho_*m_*\}$,
we obtain
\begin{align}\label{decay}
\mathcal A'(t)
&+\|\mathbb D\boldsymbol{u}\|_{L^2(\Omega)}^2
+\frac{m_*}{2}\|\mu\|_{L^2(\Omega)}^2
+\rho_*m_*\|\chi^2-1\|_{L^2(\Omega)}^2
\nonumber\\
&\quad+
m_*\|\nabla(\chi^2-1)\|_{L^2(\Omega)}^2
\leq0.
\end{align}
The density bounds, the Poincar\'e's and Korn's inequalities,
and the preceding estimate for $\nabla\chi$ give
\[
\mathcal A(t)\leq C\left(
\|\mathbb D\boldsymbol{u}\|_{L^2(\Omega)}^2
+\|\mu\|_{L^2(\Omega)}^2
+\|\chi^2-1\|_{L^2(\Omega)}^2
\right).
\]
Thus \eqref{decay} implies
$\mathcal A'(t)+\sigma\mathcal A(t)\leq0$
for some $\sigma>0$ independent of $T$.
Multiplying by $e^{\sigma t}$ and integrating, we find
\[
\mathcal A(t)\leq\mathcal A(0)e^{-\sigma t}.
\]
Since $\mathcal A(0)\leq C\mathcal E(0)$, this yields}
\[
\|\boldsymbol{u}(t)\|_{L^2(\Omega)}^2
+\|\nabla\chi(t)\|_{L^2(\Omega)}^2
+\|\chi^2(t)-1\|_{L^2(\Omega)}^2
\leq C\mathcal E(0)e^{-\sigma t},
\qquad 0\leq t\leq T.
\]
\hfill$\Box$

\vskip2mm
\noindent
{\it\bfseries  Proof of Theorem \ref{Th8.1}}. Recall the constant $\delta$ and $C_2$ in Proposition \ref{pr8.1} and $\delta_0$ and $C_1$ in Proposition \ref{R8.1}, and, without loss of generality, assume $C_1, C_2 \geq 1$. It follows from the definition of $\mathcal{E}(0)$ and (\ref{8.3}) that there exists a positive constant $C_0 \geq 1$ such that $\mathcal{E}(0) \le C_0 \varepsilon_0$ for

\begin{align}      \label{8.36}   
\varepsilon_0=\min \left\{\frac{\delta_0}{2 C_0 C_1 C_2}, \frac{\delta}{2 C_0 C_1 C_2}\right\}.
\end{align}
Define
\begin{align}      \label{8.37}   
T^*=\mathop {\mathrm{sup}} \limits_{T}\left\{ \begin{array}{c|c}
	T>0&		\begin{array}{l}
	\mathrm{There}~\mathrm{exists}~\mathrm{a}~\mathrm{solution}~\mathrm{to}~$(\ref{1.1})$-$(\ref{1.2})$~\mathrm{on}~[0,T]\\
	\mathrm{satisfying}~\mathop {\mathrm{sup}} \limits_{0\le t \le T}\mathcal{E} (t)\le 2 C_0 C_1 C_2\varepsilon _0.\\
\end{array}\\
\end{array} \right\}
\end{align}
Thus, it follows from Theorem \ref{Th6.1} that $T^* \geq T_0>0$. We are going to show that $T^*=\infty$ by contradiction. Suppose that $T^*<\infty$. Then for any $0<T<T^*$, it follows from $(\ref{8.36})$ and $(\ref{8.37})$ that
\begin{align}      \label{8.38}   
\mathop {\mathrm{sup}} \limits_{0\le t\le T}\mathcal{E} (t)\le 2 C_0 C_1 C_2\varepsilon _0 \le \delta.
\end{align}
Then Proposition \ref{pr8.1} shows
\begin{align}      \label{8.39}   
\sup _{0\le t\le T} \mathcal{E}(t) \leq C_2\left(\varepsilon_0+\mathcal{E}(0)\right) \leq C_2\left(1+C_0\right) \varepsilon_0,
\end{align}
which implies
\begin{align}      \label{8.40}   
\sup _{0\le t\le T} \mathcal{E}(t) \leq \delta_0.
\end{align}
Now taking the above $T$ as the initial time, one can apply $(\ref{8.39})$ and Proposition \ref{R8.1} again to find that there is a solution to $(\ref{1.1})$-$(\ref{1.2})$ on $\left[T, T+T_0\right]$ satisfying
\begin{align}      \label{8.41}   
\sup _{T\le t\le T+T_0} \mathcal{E}(t) \leq C_1 \mathcal{E}(T) \leq C_1 C_2\left(1+C_0\right) \varepsilon_0 \leq 2 C_0 C_1 C_2 \varepsilon_0.
\end{align}
This contradicts the definition of $T^*$, and so $T^*=\infty$. Therefore, the three-dimensional local strong solution in Theorem \ref{Th6.1} can be extended globally.
\hfill$\Box$

\section{Weak-Strong Uniqueness}
In this section, we prove a weak-strong uniqueness theorem for the nonhomogeneous incompressible NSAC system with concentration-dependent mobility.
The proof relies on a relative energy argument; by treating the density difference through the transport equation, we avoid requiring spatial regularity on the weak solution density, thereby requiring only $\rho_0\in L^\infty(\Omega)$.
\begin{Theorem}       
\label{Th7.1}
Let $({\rho},{\boldsymbol{u}},{\chi})$ be the weak solution obtained in Theorem \ref{Th4.1}, with $d=2,3$.
Moreover, we assume that there exists a strong solution $(\bar{\rho},\bar{\boldsymbol{u}},\bar{\chi})$
to the problem $(\ref{1.1})$--$(\ref{1.2})$ satisfying
$\bar{\rho},\bar{\boldsymbol{u}}\in \mathcal{C} (\bar{\Omega}\times [0,T])$, $\left| \bar{\chi}\right|\le 1$ a.e. in $\Omega\times(0, T)$,
$\nabla \bar{\rho}, \partial _t\bar{\boldsymbol{u}}\in L^2( 0,T; L^3(\Omega ))$,
$\nabla \bar{\boldsymbol{u}}, \Delta \bar{\chi}, \bar{\mu} \in L^2\left( 0,T; L^{\infty}(\Omega ) \right)$, and
\begin{align}      \label{7.1}   
\left(\bar{\rho}, \bar{\boldsymbol{u}}, \bar{\chi} \right)(\cdot , 0)=\left(\rho _0,\boldsymbol{u}_0, \chi _0\right) \quad
\mathrm{in} ~\Omega.
\end{align}
Then $\rho =\bar{\rho},\boldsymbol{u}=\bar{\boldsymbol{u}},\chi =\bar{\chi},\mu =\bar{\mu}$ a.e. in $\Omega \times \left( 0,T \right)$.
\end{Theorem}   
\begin{Remark} \label{R-WSU}
{
A relative energy approach has also been used by
Grasselli, Hurm, and Poiatti
\cite[Theorem 3.10]{G-H-P}
to establish conditional weak--strong uniqueness for
the conserved Navier--Stokes--Allen--Cahn system with
unmatched viscosities and a singular potential.
Our proof follows the same general strategy, with
an additional estimate for the density difference
based on the transport equation.}
The uniqueness result for local strong solutions in
\cite{L-H} requires $\rho_0\in W^{2,6}(\Omega)$.
In contrast, our weak--strong uniqueness argument does
not require higher-order spatial regularity of the
weak solution density, whose initial value is only
assumed to belong to $L^\infty(\Omega)$.
The comparison strong solution must, however, satisfy
the additional regularity assumptions of
Theorem~\ref{Th7.1}.
Under these assumptions, any weak solution in the class
specified by that theorem coincides with the strong
solution originating from the same initial data.
\end{Remark}

\noindent
{\bf Proof of Theorem \ref{Th7.1}.} The proof is divided into four steps.
\vskip1mm
\noindent
{\it\bfseries Step 1. Estimate for the velocity difference}\quad
Multiplying $(\ref{1.1})_2$ by $\bar{\boldsymbol{u}}$ and integrating over $\Omega \times \left( 0,t \right)$,
we have
\begin{align}      \label{7.2}   
	&\int_{\Omega}{\rho \boldsymbol{u}\cdot \bar{\boldsymbol{u}}\mathrm{d}x}
+\int_0^t{\int_{\Omega}{\eta (\chi )\mathbb{D} \boldsymbol{u}:\mathbb{D} \bar{\boldsymbol{u}}\mathrm{d}x\mathrm{d}s}}
\\ \nonumber
	&=\int_{\Omega}{\rho _0}\boldsymbol{u}_0\cdot\bar{\boldsymbol{u}}(0)\mathrm{d}x+\int_0^t{\int_{\Omega}{\nabla}\chi \otimes \nabla \chi :\nabla \bar{\boldsymbol{u}}\mathrm{d}x\mathrm{d}s}+\int_0^t{\int_{\Omega}{\rho \boldsymbol{u}\cdot \left\{ \partial _t\bar{\boldsymbol{u}}+\boldsymbol{u}\cdot \nabla \bar{\boldsymbol{u}} \right\} \mathrm{d}x\mathrm{d}s}}.
\end{align}
Next, we rewrite the equation $(\ref{1.1})_2$ satisfied by the strong solution as follows
\begin{align}      \label{7.3}   
	&\rho \partial _t\bar{\boldsymbol{u}}+\rho \boldsymbol{u}\cdot \nabla \bar{\boldsymbol{u}}-\mathrm{div(}\eta (\chi )\mathbb{D} \bar{\boldsymbol{u}})+\nabla \bar{p}\\
	=&-\mathrm{div(}\nabla \bar{\chi}\otimes \nabla \bar{\chi})+(\rho -\bar{\rho})(\partial _t\bar{\boldsymbol{u}}+\bar{\boldsymbol{u}}\cdot \nabla \bar{\boldsymbol{u}})+\rho (\boldsymbol{u}-\bar{\boldsymbol{u}})\cdot \nabla \bar{\boldsymbol{u}}
-{\rm div}((\eta (\chi )-\eta (\bar{\chi}))\mathbb{D} \bar{\boldsymbol{u}}).\nonumber
\end{align}
Multiplying $(\ref{7.3})$ by $\boldsymbol{u}$ and integrating over $\Omega \times \left( 0,t \right)$,  we find
\begin{align}      \label{7.4}   
	&\int_0^t{\int_{\Omega}{\left\{ \rho \partial _t\bar{\boldsymbol{u}}
+\rho \boldsymbol{u}\cdot \nabla \bar{\boldsymbol{u}} \right\}  \cdot \boldsymbol{u} \mathrm{d}x\mathrm{d}s}}
+\int_0^t{\int_{\Omega}{\eta(\chi )\mathbb{D} \bar{\boldsymbol{u}}: \mathbb{D} \boldsymbol{u}\mathrm{d}x\mathrm{d}s}} \nonumber \\
	=&\int_0^t{\int_{\Omega}{\nabla\bar{\chi}\otimes \nabla \bar{\chi}:\nabla \boldsymbol{u}\mathrm{d}x\mathrm{d}s}}+\int_0^t{\int_{\Omega}{(\rho -\bar{\rho})\left( \partial _t\bar{\boldsymbol{u}}+\bar{\boldsymbol{u}}\cdot \nabla \bar{\boldsymbol{u}} \right) \cdot \boldsymbol{u}\mathrm{d}x\mathrm{d}s}}  \\
	&+\int_0^t{\int_{\Omega}{\rho(\boldsymbol{u}-\bar{\boldsymbol{u}})\cdot \nabla \bar{\boldsymbol{u}}\cdot \boldsymbol{u}\mathrm{d}x\mathrm{d}s}}+\int_0^t{\int_{\Omega}{(\eta (\chi )-\eta (\bar{\chi}))\mathbb{D} \bar{\boldsymbol{u}}: \mathbb{D} \boldsymbol{u}\mathrm{d}x\mathrm{d}s}}.\nonumber
\end{align}
By summing $(\ref{7.2})$ and $(\ref{7.4})$, we arrive at
\begin{align}      \label{7.5}   
	&\int_{\Omega}{\rho \boldsymbol{u}\cdot \bar{\boldsymbol{u}}\mathrm{d}x}
+2\int_0^t{\int_{\Omega}{\eta (\chi )\mathbb{D} \bar{\boldsymbol{u}}: \mathbb{D} \boldsymbol{u}\mathrm{d}x\mathrm{d}s}}\nonumber \\
	&=\int_{\Omega}{\rho _0\boldsymbol{u}_0\cdot \bar{\boldsymbol{u}}(0)\mathrm{d}x}+\int_0^t{\int_{\Omega}{(\nabla \bar{\chi}\otimes \nabla \bar{\chi}:\nabla \boldsymbol{u}+\nabla \chi \otimes \nabla \chi :\nabla \bar{\boldsymbol{u}})\mathrm{d}x\mathrm{d}s}}\\
	&\quad +\int_0^t{\int_{\Omega}{(}}\rho -\bar{\rho})\left( \partial _t\bar{\boldsymbol{u}}+\bar{\boldsymbol{u}}\cdot \nabla \bar{\boldsymbol{u}} \right) \cdot \boldsymbol{u}\mathrm{d}x\mathrm{d}s+\int_0^t{\int_{\Omega}{\rho}}(\boldsymbol{u}-\bar{\boldsymbol{u}})\cdot \nabla \bar{\boldsymbol{u}}\cdot \boldsymbol{u}\mathrm{d}x\mathrm{d}s \nonumber\\
	&\quad +\int_0^t{\int_{\Omega}{(\eta (\chi )-\eta (\bar{\chi}))\mathbb{D} \bar{\boldsymbol{u}}: \mathbb{D} \boldsymbol{u}\mathrm{d}x\mathrm{d}s}}.\nonumber
\end{align}
Multiplying $(\ref{7.3})$ by $\bar{\boldsymbol{u}}$ and integrating over $\Omega \times \left( 0,t \right)$,  we get
\begin{align}      \label{7.6}   
	&\frac{1}{2}\int_{\Omega}{\rho |\bar{\boldsymbol{u}}|^2\mathrm{d}x}+\int_0^t{\int_{\Omega}{\eta (\chi )|\mathbb{D} \bar{\boldsymbol{u}}|^2\mathrm{d}x\mathrm{d}s}} \nonumber\\
	&=\frac{1}{2}\int_{\Omega}{\rho _0|\bar{\boldsymbol{u}}\left( 0 \right) |^2\mathrm{d}x}+\int_0^t{\int_{\Omega}{\nabla \bar{\chi}\otimes \nabla \bar{\chi}:\nabla \bar{\boldsymbol{u}}\mathrm{d}x\mathrm{d}s}}+\int_0^t{\int_{\Omega}{(\rho -\bar{\rho})\left( \partial _t\bar{\boldsymbol{u}}+\bar{\boldsymbol{u}}\cdot \nabla \bar{\boldsymbol{u}} \right) \cdot \bar{\boldsymbol{u}}\mathrm{d}x\mathrm{d}s}}  \nonumber
\\
	&\quad+\int_0^t{\int_{\Omega}{\rho (\boldsymbol{u}-\bar{\boldsymbol{u}})\cdot \nabla \bar{\boldsymbol{u}}\cdot \bar{\boldsymbol{u}}\mathrm{d}x\mathrm{d}s}}+\int_0^t{\int_{\Omega}{(\eta (\chi )-\eta (\bar{\chi}))|\mathbb{D} \bar{\boldsymbol{u}}|^2\mathrm{d}x\mathrm{d}s}}.
\end{align}
By exploiting $(\ref{4.7})$, we have
\begin{align}      \label{7.7}   
	&\frac{1}{2}\int_{\Omega}{\rho |\boldsymbol{u}|^2\mathrm{d}x}
+\frac{1}{2}\int_{\Omega}{|\nabla \chi |^2\mathrm{d}x}
+\int_{\Omega}{\rho F(\chi )\mathrm{d}x}+\int_0^t{\int_{\Omega}{\left( \eta (\chi )|\mathbb{D} \boldsymbol{u}|^2+m(\chi) |\mu|^2 \right) \mathrm{d}x\mathrm{d}s}} \nonumber\\
	&\le  \frac{1}{2}\int_{\Omega}{\rho _0}\left| \boldsymbol{u}_0 \right|^2\mathrm{d}x+\frac{1}{2}\int_{\Omega}{\left| \nabla \chi _0 \right|^2}\mathrm{d}x+\int_{\Omega}{\rho _0}F\left( \chi _0 \right) \mathrm{d}x.
\end{align}
Then, we add up $(\ref{7.6})$ and $(\ref{7.7})$ and subtract $(\ref{7.5})$ to obtain
\begin{align}      \label{7.8}   
	&\frac{1}{2}\int_{\Omega}{\rho |\boldsymbol{u}-\bar{\boldsymbol{u}}|^2\mathrm{d}x}
+\frac{1}{2}\int_{\Omega}{|\nabla \chi |^2\mathrm{d}x}
+\int_{\Omega}{\rho F(\chi )\mathrm{d}x}
+\int_0^t{\int_{\Omega}{\left( \eta (\chi )|\mathbb{D} \boldsymbol{u}-\mathbb{D} \bar{\boldsymbol{u}}|^2+ m(\chi) |\mu|^2 \right) \mathrm{d}x\mathrm{d}s}} \nonumber
\\
	&\le \frac{1}{2}\int_{\Omega}{\left| \nabla \chi _0 \right|^2\mathrm{d}x}
+\int_{\Omega}{\rho _0F(\chi _0)\mathrm{d}x}
+\int_0^t{\int_{\Omega}{\left( \nabla \bar{\chi}\otimes \nabla \bar{\chi}:\nabla (\bar{\boldsymbol{u}}-\boldsymbol{u})-\nabla \chi \otimes \nabla \chi :\nabla \bar{\boldsymbol{u}} \right) \mathrm{d}x\mathrm{d}s}} \nonumber\\
	&\quad +\int_0^t{\int_{\Omega}{(\rho -\bar{\rho})(\partial _t\bar{\boldsymbol{u}}+\bar{\boldsymbol{u}}\cdot \nabla \bar{\boldsymbol{u}})\cdot (\bar{\boldsymbol{u}}-\boldsymbol{u})\mathrm{d}x\mathrm{d}s}}+\int_0^t{\int_{\Omega}{\rho (\boldsymbol{u}-\bar{\boldsymbol{u}})\cdot \nabla \bar{\boldsymbol{u}}\cdot (\bar{\boldsymbol{u}}-\boldsymbol{u})\mathrm{d}x\mathrm{d}s}}\nonumber\\
	&\quad +\int_0^t{\int_{\Omega}{(\eta (\chi )-\eta (\bar{\chi}))\mathbb{D} \bar{\boldsymbol{u}} : (\mathbb{D} \bar{\boldsymbol{u}}-\mathbb{D} \boldsymbol{u})\mathrm{d}x\mathrm{d}s}}.
\end{align}

\noindent
{\it\bfseries  Step 2. Phase-field and chemical potential differences} \quad
Multiplying $(\ref{1.1})_4$ by $\bar{\mu}$ and integrating over $\Omega \times \left( 0,t \right)$,  we get
\begin{align}      \label{7.9}   
	&\int_{\Omega}{\Big( \frac{1}{2}|\nabla \bar{\chi}|^2+\bar{\rho}F(\bar{\chi}) \Big)}\mathrm{d}x+\int_0^t{\int_{\Omega}{ m(\chi) |\bar{\mu}|^2}}\mathrm{d}x\mathrm{d}s \nonumber
\\
	=&\int_{\Omega}{\left( \frac{1}{2}\left| \nabla \bar{\chi}_0 \right|^2+\bar{\rho}_0F( \bar{\chi}_0 ) \right)}\mathrm{d}x
-\int_0^t{\int_{\Omega}{\nabla}}\bar{\chi}\otimes \nabla \bar{\chi}:\nabla \bar{\boldsymbol{u}}\mathrm{d}x\mathrm{d}s
\\
& + \int_0^t{\int_{\Omega}{\frac{m\left( \bar{\chi} \right) -m\left( \chi \right)}{\bar{\rho}}\bar{\mu}\Delta \bar{\chi}\mathrm{d}x\mathrm{d}s}}
+\int_0^t{\int_{\Omega}{\left( m\left( \chi \right) -m\left( \bar{\chi} \right) \right) \bar{\mu}F^{\prime}(\bar{\chi}) \mathrm{d}x\mathrm{d}s}}.
\nonumber
\end{align}
Next, we also need to obtain the cross terms of $\chi,\bar{\chi}$ and $\mu,\bar{\mu}$. According to $(\ref{1.1})_4$, we have
\begin{align*}
\partial _t\chi +\boldsymbol{u}\cdot \nabla \chi =-\frac{{m(\chi)} }{\rho}\mu.
\end{align*}
Multiplying the above equation by $\Delta \bar{\chi}$ and integrating over $\Omega \times \left( 0,t \right)$, we get
\begin{align}      \label{7.10}   
\int_0^t{\int_{\Omega}{\partial _t(\nabla \chi )\cdot \nabla \bar{\chi}\mathrm{d}x\mathrm{d}s}}=\int_0^t{\int_{\Omega}{\boldsymbol{u}\cdot \nabla \chi \Delta \bar{\chi}\mathrm{d}x\mathrm{d}s}}+\int_0^t{\int_{\Omega}{\frac{{m(\chi)} }{\rho}\mu \Delta \bar{\chi}\mathrm{d}x\mathrm{d}s}}.
\end{align}
Similarly, it holds that
\begin{align}      \label{7.11}   
\int_0^t{\int_{\Omega}{\partial _t(\nabla \bar{\chi})\cdot \nabla \chi \mathrm{d}x\mathrm{d}s}}=\int_0^t{\int_{\Omega}{\bar{\boldsymbol{u}}\cdot \nabla \bar{\chi}\Delta \chi \mathrm{d}x\mathrm{d}s}}+\int_0^t{\int_{\Omega}{\frac{{m(\bar\chi)} }{\bar{\rho}}\bar{\mu}\Delta \chi \mathrm{d}x\mathrm{d}s}}.
\end{align}
By summing $(\ref{7.10})$ and $(\ref{7.11})$ together, we have
\begin{align}      \label{7.12}   
\int_{\Omega}{\nabla \chi \cdot \nabla \bar{\chi}\mathrm{d}x}=&\int_{\Omega}{\nabla \chi _0\cdot \nabla \bar{\chi}_0\mathrm{d}x}+\int_0^t{\int_{\Omega}{\boldsymbol{u}\cdot \nabla \chi \Delta \bar{\chi}\mathrm{d}x\mathrm{d}s}}+\int_0^t{\int_{\Omega}{\bar{\boldsymbol{u}}\cdot \nabla \bar{\chi}\Delta \chi \mathrm{d}x\mathrm{d}s}} \nonumber \\
&+\int_0^t{\int_{\Omega}{\frac{{m(\chi)} }{\rho}\mu \Delta \bar{\chi}\mathrm{d}x\mathrm{d}s}}
+\int_0^t{\int_{\Omega}{\frac{{m(\bar \chi)} }{\bar{\rho}}\bar{\mu}\Delta \chi \mathrm{d}x\mathrm{d}s}}.
\end{align}
Moreover, there holds
\begin{align}      \label{7.13}   
	&\int_0^t{\int_{\Omega}{2{m(\chi)} \mu \bar{\mu}\mathrm{d}x\mathrm{d}s}}\nonumber
\\
&=\int_0^t{\int_{\Omega}{{m(\chi)}\mu \bar{\mu}\mathrm{d}x\mathrm{d}s}}
+\int_0^t{\int_{\Omega}{{m(\chi)} \mu \bar{\mu}\mathrm{d}x\mathrm{d}s}} \nonumber
\\
	&=\int_0^t{\int_{\Omega}{{m(\chi)} \mu \left( -\frac{1}{\bar{\rho}}\Delta \bar{\chi}+F^{\prime}(\bar{\chi}) \right) \mathrm{d}x\mathrm{d}s}}
+\int_0^t{\int_{\Omega}{{m(\chi)} \left( -\frac{1}{\rho}\Delta \chi +F^{\prime}(\chi ) \right) \bar{\mu}\mathrm{d}x\mathrm{d}s}} \nonumber
\\
	&=\int_0^t{\int_{\Omega}{-\frac{{m(\chi)} }{\bar{\rho}}\mu \Delta \bar{\chi}\mathrm{d}x\mathrm{d}s}}
+\int_0^t{\int_{\Omega}{{m(\chi)} \mu F^{\prime}(\bar{\chi})\mathrm{d}x\mathrm{d}s}}
+\int_0^t{\int_{\Omega}{-\frac{{m(\chi)} }{\rho}\bar{\mu}\Delta \chi \mathrm{d}x\mathrm{d}s}} \nonumber
\\
&\quad +\int_0^t{\int_{\Omega}{{m(\chi)} \bar{\mu}F^{\prime}(\chi )\mathrm{d}x\mathrm{d}s}}.
\end{align}
Next, we add $(\ref{7.8})$ and $(\ref{7.9})$ and subtract $(\ref{7.12})$ and $(\ref{7.13})$ to obtain
\begin{align}          \label{7.14-1}   
	&\frac{1}{2}\int_{\Omega}{\left( \rho |\boldsymbol{u}-\bar{\boldsymbol{u}}|^2+|\nabla \chi -\nabla \bar{\chi}|^2 \right)\mathrm{d}x}
+\int_{\Omega}{\left( \rho F(\chi )+\bar{\rho}F(\bar{\chi}) \right)\mathrm{d}x} \nonumber
\\
	&+\int_0^t{\int_{\Omega}{\left( \eta (\chi )|\mathbb{D}\boldsymbol{u}-\mathbb{D} \bar{\boldsymbol{u}}|^2
+m(\chi)|\mu -\bar{\mu}|^2 \right)\mathrm{d}x\mathrm{d}s}} \nonumber
\\
	&\le 2\int_{\Omega}{\rho _0}F\left( \chi _0 \right) \mathrm{d}x+\underbrace{\int_0^t{\int_{\Omega}{\left(\nabla \bar{\chi}\otimes \nabla \bar{\chi}:\nabla (\bar{\boldsymbol{u}}-\boldsymbol{u})-\nabla \chi \otimes \nabla \chi :\nabla \bar{\boldsymbol{u}}  \right) \mathrm{d}x\mathrm{d}s}} }_{N_1} \nonumber
\\
	&\quad +\underbrace{\int_0^t{\int_{\Omega}{-\nabla \bar{\chi}\otimes \nabla \bar{\chi}:\nabla \bar{\boldsymbol{u}}}}\mathrm{d}x\mathrm{d}s}_{N_2}+\underbrace{\int_0^t{\int_{\Omega}{-\boldsymbol{u}\cdot \nabla \chi \Delta \bar{\chi}}}\mathrm{d}x\mathrm{d}s}_{N_3}+\underbrace{\int_0^t{\int_{\Omega}{-\bar{\boldsymbol{u}}\cdot \nabla \bar{\chi}\Delta \chi}}\mathrm{d}x\mathrm{d}s}_{N_4} \nonumber
\\
	&\quad +\int_0^t{\int_{\Omega}{(\rho -\bar{\rho})(\partial _t\bar{\boldsymbol{u}}+\bar{\boldsymbol{u}}\cdot \nabla \bar{\boldsymbol{u}})\cdot (\bar{\boldsymbol{u}}-\boldsymbol{u})\mathrm{d}x\mathrm{d}s}}+\int_0^t{\int_{\Omega}{\rho (\boldsymbol{u}-\bar{\boldsymbol{u}})\cdot \nabla \bar{\boldsymbol{u}}\cdot (\bar{\boldsymbol{u}}-\boldsymbol{u})\mathrm{d}x\mathrm{d}s}}
\\
	&\quad +\int_0^t{\int_{\Omega}{(\eta (\chi )-\eta (\bar{\chi}))\mathbb{D} \bar{\boldsymbol{u}}: (\mathbb{D} \bar{\boldsymbol{u}}-\mathbb{D} \boldsymbol{u})\mathrm{d}x\mathrm{d}s}}+\underbrace{\int_0^t{\int_{\Omega}{\left( m\left( \chi \right) -m\left( \bar{\chi} \right) \right) \bar{\mu}F^{\prime}(\bar{\chi})}}\mathrm{d}x\mathrm{d}s}_{{N}_5} \nonumber
\\
	&\quad +\underbrace{\int_0^t{\int_{\Omega}{-m\left( \chi \right)}\mu F^{\prime}(\bar{\chi})}\mathrm{d}x\mathrm{d}s}_{{N}_6}+\underbrace{\int_0^t{\int_{\Omega}{-m\left( \chi \right)}\bar{\mu}F^{\prime}(\chi )}\mathrm{d}x\mathrm{d}s}_{{N}_7}
+\underbrace{\int_0^t{\int_{\Omega}{\frac{m\left( \bar{\chi} \right) }{\bar{\rho}}\bar{\mu}\Delta \bar{\chi}}}\mathrm{d}x\mathrm{d}s}_{{N}_8} \nonumber
\\
	&\quad +\underbrace{\int_0^t{\int_{\Omega}{-\frac{m\left( \chi \right)}{\bar{\rho}}\bar{\mu}\Delta \bar{\chi}}}\mathrm{d}x\mathrm{d}s}_{{N}_9}
+\underbrace{\int_0^t{\int_{\Omega}{-\frac{m\left( \chi \right)}{\rho}\mu \Delta \bar{\chi}}}\mathrm{d}x\mathrm{d}s}_{{N}_{10}}
+\underbrace{\int_0^t{\int_{\Omega}{-\frac{m\left( \bar{\chi} \right)}{\bar{\rho}}}\bar{\mu}\Delta \chi}\mathrm{d}x\mathrm{d}s}_{{N}_{11}} \nonumber
\\
	&\quad +\underbrace{\int_0^t{\int_{\Omega}{\frac{m\left( \chi \right)}{\bar{\rho}}\mu \Delta \bar{\chi}}}\mathrm{d}x\mathrm{d}s}_{{N}_{12}}
+\underbrace{\int_0^t{\int_{\Omega}{\frac{m\left( \chi \right)}{\rho}}\bar{\mu}\Delta \chi}\mathrm{d}x\mathrm{d}s}_{{N}_{13}}. \nonumber
\end{align}
We rewrite some terms on the right side of inequality $(\ref{7.14-1})$ as follows
\begin{align*}    
	\sum_{i=1}^4{N_i} =&\int_0^t{\int_{\Omega}{\left( -\nabla \bar{\chi}\otimes \nabla \bar{\chi}:\nabla \boldsymbol{u}-\nabla \chi \otimes \nabla \chi :\nabla \bar{\boldsymbol{u}}-\boldsymbol{u}\cdot \nabla \chi \Delta \bar{\chi}-\bar{\boldsymbol{u}}\cdot \nabla \bar{\chi}\Delta \chi \right) \mathrm{d}x\mathrm{d}s}}
\\
	=&\int_0^t{\int_{\Omega}{\left( \boldsymbol{u}\cdot \nabla \bar{\chi}\Delta \bar{\chi}+\bar{\boldsymbol{u}}\cdot \nabla \chi \Delta \chi -\boldsymbol{u}\cdot \nabla \chi \Delta \bar{\chi}-\bar{\boldsymbol{u}}\cdot \nabla \bar{\chi}\Delta \chi \right) \mathrm{d}x\mathrm{d}s}}
\\
	=&\int_0^t{\int_{\Omega}{\left( \bar{\boldsymbol{u}}\cdot \nabla (\chi -\bar{\chi})\Delta \chi -\boldsymbol{u}\cdot \nabla (\chi -\bar{\chi})\Delta \bar{\chi} \right) \mathrm{d}x\mathrm{d}s}}
\\
	=&\int_0^t{\int_{\Omega}{\left\{ \bar{\boldsymbol{u}}\cdot \nabla (\chi -\bar{\chi})\Delta \chi -\bar{\boldsymbol{u}}\cdot \nabla (\chi -\bar{\chi})\Delta \bar{\chi}+\bar{\boldsymbol{u}}\cdot \nabla (\chi -\bar{\chi})\Delta \bar{\chi}-\boldsymbol{u}\cdot \nabla (\chi -\bar{\chi})\Delta \bar{\chi} \right\} \mathrm{d}x\mathrm{d}s}}
\\
	=&\int_0^t{\int_{\Omega}{\left\{ \bar{\boldsymbol{u}}\cdot \nabla (\chi -\bar{\chi})\Delta (\chi -\bar{\chi})+(\bar{\boldsymbol{u}}-\boldsymbol{u})\cdot \nabla (\chi -\bar{\chi})\Delta \bar{\chi} \right\} \mathrm{d}x\mathrm{d}s}}
\\
	=&\int_0^t{\int_{\Omega}{\left\{{ (\bar{\boldsymbol{u}}-\boldsymbol{u})\cdot \nabla (\chi -\bar{\chi})\Delta \bar{\chi}-\nabla \bar{\boldsymbol{u}}:\nabla (\chi -\bar{\chi})\otimes \nabla (\chi -\bar{\chi})} \right\} \mathrm{d}x\mathrm{d}s}},
\\
	\sum_{i=5}^7{N_i} =&\int_0^t{\int_{\Omega}{m\left( \chi \right) \left( \bar\mu -\mu \right) \left( F^{\prime}(\bar{\chi})-F^{\prime}(\chi ) \right) \mathrm{d}x\mathrm{d}s}}-\int_0^t{\int_{\Omega}{m\left( \chi \right) \left( \bar{\mu}F^{\prime}(\bar\chi )+\mu F^{\prime}(\chi ) \right) \mathrm{d}x\mathrm{d}s}}  \\
	&-\int_0^t{\int_{\Omega}{m\left( \bar{\chi} \right) \bar{\mu}F^{\prime}(\bar{\chi})\mathrm{d}x\mathrm{d}s}}+\int_0^t{\int_{\Omega}{m\left( \chi \right) \bar{\mu}F^{\prime}(\bar{\chi})\mathrm{d}x\mathrm{d}s}} \\
	=&\int_0^t{\int_{\Omega}{m\left( \chi \right) \left(\bar\mu -\mu \right) \left( F^{\prime}(\bar{\chi})-F^{\prime}(\chi ) \right) \mathrm{d}x\mathrm{d}s}}-\int_0^t{\int_{\Omega}{m\left( \chi \right) \mu F^{\prime}(\chi )\mathrm{d}x\mathrm{d}s}} \\
	&-\int_0^t{\int_{\Omega}{m\left( \bar{\chi} \right) \bar{\mu}F^{\prime}(\bar{\chi})\mathrm{d}x\mathrm{d}s}} \\
	=&\int_0^t{\int_{\Omega}{m\left( \chi \right) \left(\bar\mu -\mu \right) \left( \bar{\chi}^3-\chi ^3+\chi -\bar{\chi} \right) \mathrm{d}x\mathrm{d}s}}+\int_0^t{\int_{\Omega}{\left( \rho \partial _t\chi +\rho \boldsymbol{u}\cdot \nabla \chi \right) F^{\prime}(\chi )\mathrm{d}x\mathrm{d}s}} \\
	&+\int_0^t{\int_{\Omega}{\left( \bar{\rho}\partial _t\bar{\chi}+\bar{\rho}\bar{\boldsymbol{u}}\cdot \nabla \bar{\chi} \right) F^{\prime}(\bar{\chi})\mathrm{d}x\mathrm{d}s}} \\
	=&\int_0^t{\int_{\Omega}{{ m( \chi)} \left(\bar\mu -\mu \right) \left( \bar\chi -\chi \right) \left( \bar{\chi}^2+\bar{\chi}\chi +\chi ^2-1 \right) \mathrm{d}x\mathrm{d}s}}+\int_{\Omega}{\left( \rho F(\chi )+\bar{\rho}F(\bar{\chi}) \right)}\mathrm{d}x \\
    &-2\int_{\Omega}{ \rho _0F\left( \chi _0 \right)} \mathrm{d}x.
\end{align*}
%
Furthermore, noticing that
%
%
\begin{align*}
 	\Delta \bar{\chi}-\Delta \chi &=\left( \bar{\rho}F^{\prime}(\bar{\chi})-\bar{\rho}\bar{\mu} \right) -\left( \rho F^{\prime}(\chi )-\rho \mu \right)\\
	&=\bar{\rho}F^{\prime}(\bar{\chi})-\bar{\rho}F^{\prime}(\chi )+\bar{\rho}F^{\prime}(\chi )-\rho F^{\prime}(\chi )+\rho \mu -\rho \bar{\mu}+\rho \bar{\mu}-\bar{\rho}\bar{\mu}\\
	&=\bar{\rho}(\bar{\chi}-\chi )\left( \bar{\chi}^2+\bar{\chi}\chi +\chi ^2-1 \right) +(\bar{\rho}-\rho )\left( \chi ^3-\chi \right)+\rho (\mu -\bar{\mu})+(\rho -\bar{\rho})\bar{\mu},
\end{align*}
there holds
\begin{align*}  
    \sum_{i=8}^{13}{N_i} =&-\int_0^t{\int_{\Omega}{\left( \frac{1}{\bar{\rho}}-\frac{1}{\rho} \right) m\left( \chi \right) \bar{\mu}\Delta \chi}}\mathrm{d}x\mathrm{d}s+\int_0^t{\int_{\Omega}{\left( \frac{1}{\bar{\rho}}-\frac{1}{\rho} \right) m\left( \chi \right) \mu \Delta \bar{\chi}}}\mathrm{d}x\mathrm{d}s
\\
    &+\int_0^t{\int_{\Omega}{\frac{m\left( \bar{\chi} \right)}{\bar{\rho}}}\bar{\mu}\Delta \left( \bar{\chi}-\chi \right)}\mathrm{d}x\mathrm{d}s-\int_0^t{\int_{\Omega}{\frac{m\left( \chi \right)}{\bar{\rho}}}\bar{\mu}\Delta \left( \bar{\chi}-\chi \right)}\mathrm{d}x\mathrm{d}s
\\
    &+\int_0^t{\int_{\Omega}{\left( \frac{1}{\bar{\rho}}-\frac{1}{\rho} \right) m\left( \chi \right) \bar{\mu}\Delta \bar{\chi}}}\mathrm{d}x\mathrm{d}s-\int_0^t{\int_{\Omega}{\left( \frac{1}{\bar{\rho}}-\frac{1}{\rho} \right) m\left( \chi \right) \bar{\mu}\Delta \bar{\chi}}}\mathrm{d}x\mathrm{d}s
\\
    =&\int_0^t{\int_{\Omega}{\frac{\rho -\bar{\rho}}{\bar{\rho}\rho}}m\left( \chi \right) \bar{\mu}\Delta \left( \bar{\chi}-\chi \right)}\mathrm{d}x\mathrm{d}s+\int_0^t{\int_{\Omega}{\frac{\rho -\bar{\rho}}{\bar{\rho}\rho}m\left( \chi \right)}\left( \mu -\bar{\mu} \right) \Delta \bar{\chi}}\mathrm{d}x\mathrm{d}s
\\
    &+\int_0^t{\int_{\Omega}{\frac{m\left( \bar{\chi} \right) -m\left( \chi \right)}{\bar{\rho}}\bar{\mu}\Delta \left( \bar{\chi}-\chi \right) \mathrm{d}x\mathrm{d}s}}
\\
    =&\int_0^t{\int_{\Omega}{\frac{\rho -\bar{\rho}}{\rho}{ m(\chi )}\bar{\mu}(\bar{\chi}-\chi )\left( \bar{\chi}^2+\bar{\chi}\chi +\chi ^2-1 \right) \mathrm{d}x\mathrm{d}s}}
\\
    &+\int_0^t{\int_{\Omega}{-\frac{|\rho -\bar{\rho}|^2}{\rho \bar{\rho}}{m(\chi )}\bar{\mu}\left( \chi ^3-\chi \right)}}\mathrm{d}x\mathrm{d}s+\int_0^t{\int_{\Omega}{\frac{\rho -\bar{\rho}}{\bar{\rho}}{m(\chi )}\bar{\mu}(\mu -\bar{\mu})}}\mathrm{d}x\mathrm{d}s
\\
    &+\int_0^t{\int_{\Omega}{\frac{|\rho -\bar{\rho}|^2}{\rho \bar{\rho}}{m(\chi )}|\bar{\mu}|^2}}\mathrm{d}x\mathrm{d}s+\int_0^t{\int_{\Omega}{\frac{\rho -\bar{\rho}}{\bar{\rho} \rho}{m(\chi )}(\mu -\bar{\mu})\Delta \bar{\chi}}}\mathrm{d}x\mathrm{d}s
\\
    &+\int_0^t{\int_{\Omega}{{\left( m\left( \bar{\chi} \right) -m\left( \chi \right) \right) \bar{\mu}(\bar{\chi}-\chi )\left( \bar{\chi}^2+\bar{\chi}\chi +\chi ^2-1 \right)} \mathrm{d}x\mathrm{d}s}}
\\
    &+\int_0^t{\int_{\Omega}{{\frac{m(\bar{\chi})-m(\chi )}{\bar{\rho}}\bar{\mu}\left( \bar{\rho}-\rho \right) \left( \chi ^3-\chi \right)}}}\mathrm{d}x\mathrm{d}s
\\
    &+\int_0^t{\int_{\Omega}{{\frac{m(\bar{\chi})-m(\chi )}{\bar{\rho}}\rho \bar{\mu}\left( \mu -\bar{\mu} \right)}}}\mathrm{d}x\mathrm{d}s
+\int_0^t{\int_{\Omega}{{\frac{m(\bar{\chi})-m(\chi )}{\bar{\rho}}\left(\rho -\bar{\rho}\right) |\bar{\mu}|^2}}}\mathrm{d}x\mathrm{d}s.
\end{align*}
In summary, the inequality $(\ref{7.14-1})$ can be expressed in the following form
\begin{align}      \label{7.14}   
	&\frac{1}{2}\int_{\Omega}{\left( \rho |\boldsymbol{u}-\bar{\boldsymbol{u}}|^2+|\nabla \chi -\nabla \bar{\chi}|^2 \right)}\mathrm{d}x+\int_0^t{\int_{\Omega}{\left( \eta (\chi )|\mathbb{D} \boldsymbol{u}-\mathbb{D} \bar{\boldsymbol{u}}|^2+{m(\chi)}|\mu -\bar{\mu}|^2 \right)}}\mathrm{d}x\mathrm{d}s \nonumber
\\
	\le& \int_0^t{\int_{\Omega}{\left\{ {(\bar{\boldsymbol{u}}-\boldsymbol{u})\cdot \nabla (\chi -\bar{\chi})\Delta \bar{\chi}-\nabla \bar{\boldsymbol{u}}:\nabla (\chi -\bar{\chi})\otimes \nabla (\chi -\bar{\chi})} \right\} \mathrm{d}x\mathrm{d}s}} \nonumber
\\
	&+\int_0^t{\int_{\Omega}{(\rho -\bar{\rho})(\partial _t\bar{\boldsymbol{u}}+\bar{\boldsymbol{u}}\cdot \nabla \bar{\boldsymbol{u}})\cdot (\bar{\boldsymbol{u}}-\boldsymbol{u})\mathrm{d}x\mathrm{d}s}}+\int_0^t{\int_{\Omega}{\rho (\boldsymbol{u}-\bar{\boldsymbol{u}})\cdot \nabla \bar{\boldsymbol{u}}\cdot (\bar{\boldsymbol{u}}-\boldsymbol{u})\mathrm{d}x\mathrm{d}s}} \nonumber
\\
	&+\int_0^t{\int_{\Omega}{(\eta (\chi )-\eta (\bar{\chi})) \mathbb{D} \bar{\boldsymbol{u}}:  (\mathbb{D} \bar{\boldsymbol{u}}-\mathbb{D} \boldsymbol{u})\mathrm{d}x\mathrm{d}s}} \nonumber
\\
	&+\int_0^t{\int_{\Omega}{{ m(\chi)} \left(\bar\mu -\mu \right) \left( \bar\chi -\chi \right) \left( \bar{\chi}^2+\bar{\chi}\chi +\chi ^2-1 \right) \mathrm{d}x\mathrm{d}s}} \nonumber
\\
	&+\int_0^t{\int_{\Omega}{\frac{\rho -\bar{\rho}}{\rho}{ m(\chi )}\bar{\mu}(\bar{\chi}-\chi )\left( \bar{\chi}^2+\bar{\chi}\chi +\chi ^2-1 \right) \mathrm{d}x\mathrm{d}s}}
\\
    &+\int_0^t{\int_{\Omega}{-\frac{|\rho -\bar{\rho}|^2}{\rho \bar{\rho}}{m(\chi )}\bar{\mu}\left( \chi ^3-\chi \right)}}\mathrm{d}x\mathrm{d}s
    +\int_0^t{\int_{\Omega}{\frac{\rho -\bar{\rho}}{\bar{\rho}}{m(\chi )}\bar{\mu}(\mu -\bar{\mu})}}\mathrm{d}x\mathrm{d}s \nonumber
\\
    &+\int_0^t{\int_{\Omega}{\frac{|\rho -\bar{\rho}|^2}{\rho \bar{\rho}}{m(\chi )}|\bar{\mu}|^2}}\mathrm{d}x\mathrm{d}s+\int_0^t{\int_{\Omega}{\frac{\rho -\bar{\rho}}{\bar{\rho} \rho}{m(\chi )}(\mu -\bar{\mu})\Delta \bar{\chi}}}\mathrm{d}x\mathrm{d}s \nonumber
\\
    &+\int_0^t{\int_{\Omega}{{\left( m\left( \bar{\chi} \right) -m\left( \chi \right) \right) \bar{\mu}(\bar{\chi}-\chi )\left( \bar{\chi}^2+\bar{\chi}\chi +\chi ^2-1 \right)} \mathrm{d}x\mathrm{d}s}} \nonumber
\\
    &+\int_0^t{\int_{\Omega}{{\frac{m(\bar{\chi})-m(\chi )}{\bar{\rho}}\bar{\mu}\left( \bar{\rho}-\rho \right) \left( \chi ^3-\chi \right)}}}\mathrm{d}x\mathrm{d}s
    +\int_0^t{\int_{\Omega}{{\frac{m(\bar{\chi})-m(\chi )}{\bar{\rho}}\rho \bar{\mu}\left( \mu -\bar{\mu} \right)}}}\mathrm{d}x\mathrm{d}s \nonumber
\\
    &+\int_0^t{\int_{\Omega}{{\frac{m(\bar{\chi})-m(\chi )}{\bar{\rho}}\left(\rho -\bar{\rho}\right) |\bar{\mu}|^2}}}\mathrm{d}x\mathrm{d}s. \nonumber
\end{align}

\noindent
{\it\bfseries  Step 3. Closure for density and phase-field differences} \quad
Next, we estimate $\left\| \bar{\chi}-\chi \right\| _{L^2\left( \Omega \right)}$. According to $(\ref{1.1})_4-(\ref{1.1})_5$, we have
\begin{align}  \label{7.15}   
	&\partial _t(\bar{\chi}-\chi )+\bar{\boldsymbol{u}}\cdot \nabla \bar{\chi}-\boldsymbol{u}\cdot \nabla \chi \nonumber
\\
    =&\frac{(\rho -\bar{\rho})(\rho +\bar{\rho})}{\bar{\rho}^2\rho ^2}{m(\bar\chi)}\Delta \bar{\chi}
    +{\frac{m\left( \bar{\chi} \right) -m\left( \chi \right)}{\rho ^2}\Delta \bar{\chi}}
    +\frac{{m(\chi) }}{\rho ^2}(\Delta \bar{\chi}-\Delta \chi )
    -\frac{{m(\bar\chi) }}{\bar{\rho}}\left( \bar{\chi}^3-\bar{\chi} \right) \nonumber
\\
    &+\frac{{m(\chi) }}{\rho}\left( \chi ^3-\chi \right).
\end{align}
Multiplying $(\ref{7.15})$ by $\bar{\chi}-\chi$ and integrating over $\Omega \times \left( 0,t \right)$,  we can derive that
\begin{align}      \label{7.16}   
	&\frac{1}{2}\int_{\Omega}{|\bar{\chi}-\chi |^2\mathrm{d}x} \nonumber
\\
=&-\int_0^t{\int_{\Omega}{(\bar{\boldsymbol{u}}-\boldsymbol{u})\cdot \nabla \bar{\chi}(\bar{\chi}-\chi )\mathrm{d}x\mathrm{d}s}}
+\int_0^t{\int_{\Omega}{\frac{(\rho -\bar{\rho})(\rho +\bar{\rho})}{\bar{\rho}^2\rho ^2}{m(\bar\chi)}\Delta \bar{\chi}(\bar{\chi}-\chi )\mathrm{d}x\mathrm{d}s}} \nonumber
\\
	&+\int_0^t{\int_{\Omega}{\frac{{m(\chi)}}{\rho}(\mu -\bar{\mu})(\bar{\chi}-\chi )\mathrm{d}x\mathrm{d}s}}
+\int_0^t{\int_{\Omega}{\frac{{m(\chi)}\bar{\mu}}{\rho ^2}(\rho -\bar{\rho})(\bar{\chi}-\chi )\mathrm{d}x\mathrm{d}s}}\nonumber
\\
	&-\int_0^t{\int_{\Omega}{\frac{(\rho -\bar{\rho})(\rho +\bar{\rho})}{\bar{\rho}\rho ^2}{m(\chi)}\left( \bar{\chi}^3-\bar{\chi} \right) (\bar{\chi}-\chi )\mathrm{d}x\mathrm{d}s}}
+\int_0^t{\int_{\Omega}{{\frac{m\left( \bar{\chi} \right) -m\left( \chi \right)}{\rho ^2}\Delta \bar{\chi}(\bar{\chi}-\chi ) }\mathrm{d}x\mathrm{d}s}}  \nonumber
\\
&{+\int_0^t{\int_{\Omega}{{\frac{m\left( {\chi} \right) -m\left( \bar\chi \right)}{\bar\rho} (\bar{\chi}^3-\bar\chi)(\bar{\chi}-\chi ) }\mathrm{d}x\mathrm{d}s}}}.
\end{align}

According to (\ref{1.1})$_1$, we have
\begin{align}      \label{7.17}   
\partial _t(\rho -\bar{\rho})+\boldsymbol{u}\cdot\nabla(\rho -\bar{\rho})
=(\bar{\boldsymbol{u}}-\boldsymbol{u})\cdot \nabla \bar{\rho}.
\end{align}
Multiplying $(\ref{7.17})$ by $\rho-\bar{\rho}$ and integrating over $\Omega \times \left( 0,t \right)$, we obtain
\begin{align}      \label{7.18}   
\frac{1}{2}\int_{\Omega}{|\rho -\bar{\rho}|^2\mathrm{d}x}
= \int_0^t{\int_{\Omega}{(\bar{\boldsymbol{u}}-\boldsymbol{u})\cdot \nabla \bar{\rho}(\rho-\bar\rho)\mathrm{d}x\mathrm{d}s}}.
\end{align}

\noindent
{\it\bfseries  Step 4. Gr\"{o}nwall argument}\quad
Finally, by exploiting $(\ref{7.14})$, $(\ref{7.16})$ and $(\ref{7.18})$, we have
\begin{align}      \label{7.19}   
	\frac{1}{2}&\int_{\Omega}{\left\{ \rho |\boldsymbol{u}-\bar{\boldsymbol{u}}|^2+|\nabla \chi -\nabla \bar{\chi}|^2+|\bar{\chi}-\chi |^2+|\rho -\bar{\rho}|^2 \right\} \mathrm{d}x} \nonumber
\\
&\qquad~~+\int_0^t{\int_{\Omega}{\left( \eta(\chi)|\mathbb{D} \boldsymbol{u}-\mathbb{D} \bar{\boldsymbol{u}}|^2+m(\chi)|\mu -\bar{\mu}|^2 \right) \mathrm{d}x\mathrm{d}s}} \nonumber
\\
	\le& \int_0^t{\int_{\Omega}{\left\{ (\bar{\boldsymbol{u}}-\boldsymbol{u})\cdot \nabla (\chi -\bar{\chi})\Delta \bar{\chi}-\nabla \bar{\boldsymbol{u}}:\nabla (\chi -\bar{\chi})\otimes \nabla (\chi -\bar{\chi}) \right\} \mathrm{d}x\mathrm{d}s}}\nonumber
\\
	&+\int_0^t{\int_{\Omega}{(\rho -\bar{\rho})(\partial _t\bar{\boldsymbol{u}}+\bar{\boldsymbol{u}}\cdot \nabla \bar{\boldsymbol{u}})\cdot (\bar{\boldsymbol{u}}-\boldsymbol{u})\mathrm{d}x\mathrm{d}s}}
+\int_0^t{\int_{\Omega}{\rho (\boldsymbol{u}-\bar{\boldsymbol{u}})\cdot \nabla \bar{\boldsymbol{u}}\cdot (\bar{\boldsymbol{u}}-\boldsymbol{u})\mathrm{d}x\mathrm{d}s}}\nonumber
\\
	&+\int_0^t{\int_{\Omega}{(\eta (\chi )-\eta (\bar{\chi})) \mathbb{D} \bar{\boldsymbol{u}}:  (\mathbb{D} \bar{\boldsymbol{u}}-\mathbb{D} \boldsymbol{u})\mathrm{d}x\mathrm{d}s}}
+\int_0^t{\int_{\Omega}{-\frac{|\rho -\bar{\rho}|^2}{\rho \bar{\rho}}{m(\chi )}\bar{\mu}\left( \chi ^3-\chi \right)}}\mathrm{d}x\mathrm{d}s
\nonumber \\
    &+\int_0^t{\int_{\Omega}{{m( \chi )} \left(\bar{\mu} -\mu \right) \left( \bar\chi -\chi \right) \left( \bar{\chi}^2+\bar{\chi}\chi +\chi ^2-1 \right) \mathrm{d}x\mathrm{d}s}}
    +\int_0^t{\int_{\Omega}{\frac{\rho -\bar{\rho}}{\bar{\rho}}{m(\chi )}\bar{\mu}(\mu -\bar{\mu})}}\mathrm{d}x\mathrm{d}s
\nonumber \\
    &+\int_0^t{\int_{\Omega}{\frac{\rho -\bar{\rho}}{\rho}{m(\chi )}\bar{\mu}(\bar{\chi}-\chi )\left( \bar{\chi}^2+\bar{\chi}\chi +\chi ^2-1 \right) \mathrm{d}x\mathrm{d}s}}
    +\int_0^t{\int_{\Omega}{\frac{|\rho -\bar{\rho}|^2}{\rho \bar{\rho}}{m(\chi )}|\bar{\mu}|^2}}\mathrm{d}x\mathrm{d}s
\nonumber\\
    &+\int_0^t{\int_{\Omega}{\frac{\rho -\bar{\rho}}{\bar{\rho}\rho}{m(\chi )}(\mu -\bar{\mu})\Delta \bar{\chi}}}\mathrm{d}x\mathrm{d}s
    +\int_0^t{\int_{\Omega}{{\frac{m(\bar{\chi})-m(\chi )}{\bar{\rho}}\bar{\mu}\left( \bar{\rho}-\rho \right) \left( \chi ^3-\chi \right)}}}\mathrm{d}x\mathrm{d}s
\nonumber\\
    &+\int_0^t{\int_{\Omega}{{\left( m\left( \bar{\chi} \right) -m\left( \chi \right) \right) \bar{\mu}(\bar{\chi}-\chi )\left( \bar{\chi}^2+\bar{\chi}\chi +\chi ^2-1 \right)} \mathrm{d}x\mathrm{d}s}}
\nonumber \\
    &+\int_0^t{\int_{\Omega}{{\frac{m(\bar{\chi})-m(\chi )}{\bar{\rho}}\rho \bar{\mu}\left( \mu -\bar{\mu} \right)}}}\mathrm{d}x\mathrm{d}s
     +\int_0^t{\int_{\Omega}{{\frac{m(\bar{\chi})-m(\chi )}{\bar{\rho}}\left( \rho-\bar{\rho} \right) |\bar{\mu}|^2}}} \mathrm{d}x\mathrm{d}s
\nonumber \\
    & +\int_0^t{\int_{\Omega}{(\boldsymbol{u}-\bar{\boldsymbol{u}})\cdot \nabla \bar{\chi}(\bar{\chi}-\chi )}}\mathrm{d}x\mathrm{d}s
    +\int_0^t{\int_{\Omega}{\frac{(\rho -\bar{\rho})(\rho +\bar{\rho})}{\bar{\rho}^2\rho ^2}{m(\bar\chi)}\Delta \bar{\chi}(\bar{\chi}-\chi )\mathrm{d}x\mathrm{d}s}}
\nonumber\\
	&+\int_0^t{\int_{\Omega}{\frac{{m(\chi)}}{\rho}(\mu -\bar{\mu})(\bar{\chi}-\chi )\mathrm{d}x\mathrm{d}s}}
     +\int_0^t{\int_{\Omega}{\frac{{m(\chi)}\bar{\mu}}{\rho ^2}(\rho -\bar{\rho})(\bar{\chi}-\chi )\mathrm{d}x\mathrm{d}s}}
\nonumber\\
	&+\int_0^t{\int_{\Omega}{\frac{(\bar{\rho}-\rho )(\rho +\bar{\rho})}{\bar{\rho}\rho ^2}{m(\chi)}\left( \bar{\chi}^3-\bar{\chi} \right) (\bar{\chi}-\chi )\mathrm{d}x\mathrm{d}s}}
+\int_0^t{\int_{\Omega}{(\bar{\boldsymbol{u}}-\boldsymbol{u})\cdot \nabla \bar{\rho}(\rho-\bar\rho)\mathrm{d}x\mathrm{d}s}}
\nonumber\\
    &+\int_0^t{\int_{\Omega}{{\left(m\left( \bar\chi \right) -m\left( \chi \right)\right)\left(\frac{\Delta\bar{\chi}}{\rho ^2}{-\frac{\bar\chi^3-\bar\chi}{\bar\rho}}\right) (\bar{\chi}-\chi )}}}\mathrm{d}x\mathrm{d}s
    {=\sum_{i=1}^{21}R_i}.
\end{align}
{We estimate the terms on the right-hand side of
\eqref{7.19} in groups.
Throughout the following estimates, $C$ denotes a positive
constant depending on the domain, the upper and lower bounds
of the density, the bounds and Lipschitz constants of $\eta$
and $m$ on $[-1,1]$, and
$\|\bar{\boldsymbol{u}}\|_{L^\infty(\Omega\times(0,T))}$.
The constant $C_\varepsilon$ may additionally depend on
$\varepsilon>0$.}

{Recall that
$0<\rho_\ast\leq\rho,\bar{\rho}\leq\rho^\ast$ and
$|\chi|,|\bar{\chi}|\leq1$.
Thus all polynomial factors in \eqref{7.19} are bounded, and
\[
|m(\chi)-m(\bar{\chi})|
+|\eta(\chi)-\eta(\bar{\chi})|
\leq C|\chi-\bar{\chi}|.
\]
Moreover, for $d=2,3$, the Sobolev, Poincar\'e's, and Korn's
inequalities imply
\[
\|\boldsymbol{u}-\bar{\boldsymbol{u}}\|_{L^6(\Omega)}
\leq
C\|\mathbb{D}\boldsymbol{u}
-\mathbb{D}\bar{\boldsymbol{u}}\|_{L^2(\Omega)}.
\] }
{First, H\"older's and Young's inequalities give
\begin{align*}
|R_1|
&\leq
\int_0^t
\|\Delta\bar{\chi}\|_{L^\infty(\Omega)}
\|\boldsymbol{u}-\bar{\boldsymbol{u}}\|_{L^2(\Omega)}
\|\nabla\chi-\nabla\bar{\chi}\|_{L^2(\Omega)}
\,\mathrm{d}s
\\
&\quad+
\int_0^t
\|\nabla\bar{\boldsymbol{u}}\|_{L^\infty(\Omega)}
\|\nabla\chi-\nabla\bar{\chi}\|_{L^2(\Omega)}^2
\,\mathrm{d}s
\\
&\leq
C\int_0^t
\|\boldsymbol{u}-\bar{\boldsymbol{u}}\|_{L^2(\Omega)}^2
\,\mathrm{d}s
\\
&\quad+
C\int_0^t
\left(
\|\Delta\bar{\chi}\|_{L^\infty(\Omega)}^2
+\|\nabla\bar{\boldsymbol{u}}\|_{L^\infty(\Omega)}
\right)
\|\nabla\chi-\nabla\bar{\chi}\|_{L^2(\Omega)}^2
\,\mathrm{d}s.
\end{align*}
The convection term satisfies
\[
|R_3|
\leq
\int_0^t
\|\nabla\bar{\boldsymbol{u}}\|_{L^\infty(\Omega)}
\int_\Omega
\rho|\boldsymbol{u}-\bar{\boldsymbol{u}}|^2
\,\mathrm{d}x\,\mathrm{d}s.
\]
For the acceleration and density transport terms, we use
H\"older's inequality with spatial exponents $2,3,6$:
\begin{align*}
|R_2|+|R_{20}|
&\leq
\int_0^t
\|\rho-\bar{\rho}\|_{L^2(\Omega)}
\|\boldsymbol{u}-\bar{\boldsymbol{u}}\|_{L^6(\Omega)}
\\
&\qquad\qquad\times
\left(
\|\partial_t\bar{\boldsymbol{u}}
+\bar{\boldsymbol{u}}\cdot
\nabla\bar{\boldsymbol{u}}\|_{L^3(\Omega)}
+\|\nabla\bar{\rho}\|_{L^3(\Omega)}
\right)\,\mathrm{d}s
\\
&\leq
\varepsilon\int_0^t
\|\mathbb{D}\boldsymbol{u}
-\mathbb{D}\bar{\boldsymbol{u}}\|_{L^2(\Omega)}^2
\,\mathrm{d}s
\\
&\quad+
C_\varepsilon\int_0^t
\left(
\|\partial_t\bar{\boldsymbol{u}}\|_{L^3(\Omega)}^2
+\|\nabla\bar{\boldsymbol{u}}\|_{L^\infty(\Omega)}^2
+\|\nabla\bar{\rho}\|_{L^3(\Omega)}^2
\right)
\|\rho-\bar{\rho}\|_{L^2(\Omega)}^2
\,\mathrm{d}s.
\end{align*}
The Lipschitz continuity of $\eta$ similarly yields
\begin{align*}
|R_4|
&\leq
C\int_0^t
\|\nabla\bar{\boldsymbol{u}}\|_{L^\infty(\Omega)}
\|\chi-\bar{\chi}\|_{L^2(\Omega)}
\|\mathbb{D}\boldsymbol{u}
-\mathbb{D}\bar{\boldsymbol{u}}\|_{L^2(\Omega)}
\,\mathrm{d}s
\\
&\leq
\varepsilon\int_0^t
\|\mathbb{D}\boldsymbol{u}
-\mathbb{D}\bar{\boldsymbol{u}}\|_{L^2(\Omega)}^2
\,\mathrm{d}s
+
C_\varepsilon\int_0^t
\|\nabla\bar{\boldsymbol{u}}\|_{L^\infty(\Omega)}^2
\|\chi-\bar{\chi}\|_{L^2(\Omega)}^2
\,\mathrm{d}s.
\end{align*}
The five terms containing $\mu-\bar{\mu}$ can be estimated
together. By the bounds stated above,
\begin{align*}
&|R_6|+|R_7|+|R_{10}|+|R_{13}|+|R_{17}|
\\
&\leq
C\int_0^t
\|\mu-\bar{\mu}\|_{L^2(\Omega)}
\Big[
\left(1+\|\bar{\mu}\|_{L^\infty(\Omega)}\right)
\|\chi-\bar{\chi}\|_{L^2(\Omega)}
\\
&\hspace{43mm}+
\left(
\|\bar{\mu}\|_{L^\infty(\Omega)}
+\|\Delta\bar{\chi}\|_{L^\infty(\Omega)}
\right)
\|\rho-\bar{\rho}\|_{L^2(\Omega)}
\Big]\,\mathrm{d}s
\\
&\leq
\varepsilon\int_0^t
\|\mu-\bar{\mu}\|_{L^2(\Omega)}^2
\,\mathrm{d}s
\\
&\quad+
C_\varepsilon\int_0^t
\left(
1+\|\bar{\mu}\|_{L^\infty(\Omega)}^2
+\|\Delta\bar{\chi}\|_{L^\infty(\Omega)}^2
\right)
\left(
\|\chi-\bar{\chi}\|_{L^2(\Omega)}^2
+\|\rho-\bar{\rho}\|_{L^2(\Omega)}^2
\right)\,\mathrm{d}s.
\end{align*}
To estimate $R_{15}$, the Neumann elliptic estimate,
$|\bar{\chi}|\leq1$, and
$W^{2,p}(\Omega)\hookrightarrow W^{1,\infty}(\Omega)$
for a fixed $p>d$ give
\[
\|\nabla\bar{\chi}\|_{L^\infty(\Omega)}
\leq
C\left(
1+\|\Delta\bar{\chi}\|_{L^\infty(\Omega)}
\right).
\]
Hence
\begin{align*}
|R_{15}|
&\leq
\int_0^t
\|\nabla\bar{\chi}\|_{L^\infty(\Omega)}
\|\boldsymbol{u}-\bar{\boldsymbol{u}}\|_{L^2(\Omega)}
\|\chi-\bar{\chi}\|_{L^2(\Omega)}
\,\mathrm{d}s
\\
&\leq
C\int_0^t
\|\boldsymbol{u}-\bar{\boldsymbol{u}}\|_{L^2(\Omega)}^2
\,\mathrm{d}s
+
C\int_0^t
\left(
1+\|\Delta\bar{\chi}\|_{L^\infty(\Omega)}^2
\right)
\|\chi-\bar{\chi}\|_{L^2(\Omega)}^2
\,\mathrm{d}s.
\end{align*}
The remaining terms contain only the density and phase
differences. Using the bounded polynomial factors
and the Lipschitz continuity of $m$, we obtain
\begin{align*}
&|R_5|+|R_8|+|R_{9}|+|R_{11}|+|R_{12}|
+|R_{14}|+|R_{16}|+|R_{18}|+|R_{19}|+|R_{21}|
\\
&\leq
C\int_0^t
\left(
1+\|\bar{\mu}\|_{L^\infty(\Omega)}^2
+\|\Delta\bar{\chi}\|_{L^\infty(\Omega)}^2
\right)
\left(
\|\rho-\bar{\rho}\|_{L^2(\Omega)}^2
+\|\chi-\bar{\chi}\|_{L^2(\Omega)}^2
\right)\,\mathrm{d}s.
\end{align*}
In particular, the term involving $|\bar{\mu}|^2$ and
the mobility difference is estimated by
\[
|R_{14}|
\leq
C\int_0^t
\|\bar{\mu}\|_{L^\infty(\Omega)}^2
\|\chi-\bar{\chi}\|_{L^2(\Omega)}
\|\rho-\bar{\rho}\|_{L^2(\Omega)}
\,\mathrm{d}s,
\]
so that only the second power of
$\|\bar{\mu}\|_{L^\infty(\Omega)}$ is required.
Taking
\[
0<\varepsilon
\leq\min\left\{\frac{\eta_\ast}{4},
               \frac{m_\ast}{2}\right\},
\]
we absorb the two velocity dissipation contributions
from $|R_2|+|R_{20}|$ and $|R_4|$, together with the
chemical-potential dissipation contribution, into the
left-hand side of \eqref{7.19}.
The terms involving
$\|\boldsymbol{u}-\bar{\boldsymbol{u}}\|_{L^2(\Omega)}^2$
in the estimates of $|R_1|$ and $|R_{15}|$
remain on the right-hand side.
Using $\rho\geq\rho_\ast$, we estimate them by
\[
C\int_0^t
\|\boldsymbol{u}-\bar{\boldsymbol{u}}\|_{L^2(\Omega)}^2
\,\mathrm{d}s
\leq
\frac{C}{\rho_\ast}\int_0^t\int_\Omega
\rho|\boldsymbol{u}-\bar{\boldsymbol{u}}|^2
\,\mathrm{d}x\,\mathrm{d}s.
\]
These contributions are included in the constant part
of the coefficient $C(s)$ below.
Combining the preceding estimates with \eqref{7.19},
we conclude that
\begin{align*}
&\frac12\int_\Omega
\left\{
\rho|\boldsymbol{u}-\bar{\boldsymbol{u}}|^2
+|\nabla\chi-\nabla\bar{\chi}|^2
+|\chi-\bar{\chi}|^2
+|\rho-\bar{\rho}|^2
\right\}\,\mathrm{d}x
\\
&\quad+
\frac{\eta_\ast}{2}\int_0^t
\|\mathbb{D}\boldsymbol{u}
-\mathbb{D}\bar{\boldsymbol{u}}\|_{L^2(\Omega)}^2
\,\mathrm{d}s
+
\frac{m_\ast}{2}\int_0^t
\|\mu-\bar{\mu}\|_{L^2(\Omega)}^2
\,\mathrm{d}s
\\
&\leq
\int_0^t C(s)\int_\Omega
\left\{
\rho|\boldsymbol{u}-\bar{\boldsymbol{u}}|^2
+|\nabla\chi-\nabla\bar{\chi}|^2
+|\chi-\bar{\chi}|^2
+|\rho-\bar{\rho}|^2
\right\}\,\mathrm{d}x\,\mathrm{d}s,
\end{align*}
where
\begin{align*}
C(s)=C\Big(
&1+\|\partial_t\bar{\boldsymbol{u}}(s)\|_{L^3(\Omega)}^2
+\|\nabla\bar{\rho}(s)\|_{L^3(\Omega)}^2
+\|\nabla\bar{\boldsymbol{u}}(s)\|_{L^\infty(\Omega)}^2
\\
&+\|\Delta\bar{\chi}(s)\|_{L^\infty(\Omega)}^2
+\|\bar{\mu}(s)\|_{L^\infty(\Omega)}^2
\Big).
\end{align*}
The assumptions of Theorem~\ref{Th7.1} imply that
$C\in L^1(0,T)$.
Since the initial data coincide, Gr\"{o}nwall's inequality
gives $\rho=\bar{\rho}$,
$\boldsymbol{u}=\bar{\boldsymbol{u}}$, and
$\chi=\bar{\chi}$ a.e. in $\Omega\times(0,T)$.
The chemical-potential dissipation then implies
$\mu=\bar{\mu}$ a.e. in $\Omega\times(0,T)$. }
\hfill$\Box$

\section*{Acknowledgments}
\noindent
The authors thank the referee for the careful reading of the
manuscript and the constructive comments, which have helped
improve the clarity and presentation of the paper.
Li's work is supported by the National Natural Science Foundation of China (No.12371205)
and the Natural Science Foundation of Guangdong Province (No.2025A1515012026, 2025A1515040001).

\vskip2mm
\noindent {\bf Data Availability}\quad The manuscript has no associated data.

\vskip1mm
\noindent {\bf Conflicts of Interest}\quad The authors declare that they have no Conflict of interest.

\appendix
\section{Existence of approximate solutions}\label{A-A}
\noindent
In this appendix, we use a fixed-point argument to prove the existence of the approximate solutions $( \rho _m,\boldsymbol{u}_m,\chi _m)$ to the problem (\ref{A.1}) that satisfy
\begin{align*}   
&\rho _m\in \mathcal{C} ^1( \overline{Q}_T ), ~~ \quad
\boldsymbol{u}_m\in \mathcal{C}( [ 0,T ] ;\mathbf{V}_m )\cap W^{1, p}(0, T; \mathbf{V}_m), \\
&\chi _m\in \mathcal{C} ([ 0,T] ;{H^1( \Omega )} )\cap{ L^2( 0,T;H^2\left( \Omega \right)) }
\end{align*}
for any $p>1$ and $T>0$, where $Q_T = \Omega \times (0,T)$.

Let $\{\boldsymbol{w}_i\} _{i=1}^{\infty}$ be the eigenfunctions of the Stokes operator and
$\mathbf{V}_m=\mathrm{span}\left\{ \boldsymbol{w}_1,\dots \boldsymbol{w}_m \right\}$.
For
$$
\boldsymbol{u}_m(x,t)=\displaystyle\sum_{i=1}^m{g_{m}^{i}}(t)\boldsymbol{w}_i(x),
$$
we seek $g_m=( g_{m}^{1},\cdots ,g_{m}^{m})^T$ in a closed ball
\begin{align*}
\mathbf{U}_m( T_0 ) =\bigg\{ ( g_{m}^{1},\cdots ,g_{m}^{m} ) ^T\in ( C\left[ 0,T_0 \right] ) ^m\Big|
\mathop {\mathrm{sup}}_{0\le t\le T_0}\sum_{i=1}^m{| g_{m}^{i}(t) |^2}\le M,\;g_{m}^{i}(0)=\left(\boldsymbol{u}_{0m},\boldsymbol{w}_i \right) \bigg\}
\end{align*}
with $T_{0}, M>0$ to be specified later.
Given $h_m\in \mathbf{U}_m( T_0 )$ with associated velocity $\boldsymbol{v}_m=\sum_{i=1}^m{h_{m}^{i}}(t)\boldsymbol{w}_i(x)$, we solve the auxiliary system
\begin{equation} \label{A1} 
\begin{cases}
\partial _t\rho _m+\boldsymbol{v}_m\cdot \nabla \rho _m=0,
\\
\left( \rho _m\partial _t\boldsymbol{u}_m,\boldsymbol{w}_j \right) +\left( \rho _m\left( \boldsymbol{v}_m\cdot \nabla \right) \boldsymbol{u}_m,\boldsymbol{w}_j \right) +\left( \eta \left( \chi _m \right) \mathbb{D} \boldsymbol{u}_m,\nabla \boldsymbol{w}_j \right)
=\left( \nabla \chi _m\otimes \nabla \chi _m,\nabla \boldsymbol{w}_j \right),
\\
\rho _m\partial _t\chi _m+\rho _m\boldsymbol{v}_m\cdot \nabla \chi _m =-m(\chi _m)\mu _m,
\\ 
\rho _m\mu _m=-\Delta \chi _m+\rho _m F^{\prime}\left( \chi _m \right),
\end{cases}
\end{equation}
$j=1, \cdots, m$, with initial and boundary conditions as in (\ref{A.1}).

\vskip2mm
\noindent
{\it\bfseries A.1 Successive solvability.}\quad
For fixed $\boldsymbol{v}_m\in\mathcal{C}( [ 0,T_0 ] ;\mathbf{V}_m)$, by the standard method of characteristics \cite{L-S}, the transport equation yields $\rho_m\in C^1(\overline{Q}_{T_0})$ defined by $\rho _m\left( x,t \right) =\rho _{0m}\left( \mathbf{X}_m\left( 0;t,x \right) \right)$, where
\begin{equation} \label{A.4}  
\mathbf{X}_m\left( s;t,x \right) =x+\int_t^s{\boldsymbol{v}_m\left( \mathbf{X}_m\left( \tau; t,x \right) ,\tau \right)}\ \mathrm{d}\tau, \quad
\forall s,t\in \left[ 0,T_0 \right].
\end{equation}
Moreover, $\rho_m$ satisfies
\begin{align} \label{A2}
0<\rho _{\ast}\le \rho _m(x,t)\le \rho ^{\ast}, \quad \forall \left( x,t \right) \in \overline Q_{T_0}.
\end{align}
With $\rho_m$ and $\boldsymbol{v}_m$ known, the Allen-Cahn equation (\ref{A1})$_{3,4}$ becomes a semilinear parabolic problem for $\chi_m$.
Classical parabolic theory gives a unique solution $\chi_m\in C([0,T_0];H^1(\Omega))\cap L^2(0,T_0;H^2(\Omega))$ satisfying $|\chi_m|\le 1$ a.e. in $\Omega\times[0,T_0]$. Moreover, testing the Allen-Cahn equation by $m(\chi_m)^{-1}\partial_t\chi_m$, we can obtain the estimate
\begin{align} \label{A3}
\sup_{0\le s\le t}\left\| \nabla \chi _m(s) \right\|^2
\leq e^{C_1Mt}\left( \left\| \nabla \chi_0 \right\|^2+C_1t \right),
\end{align}
where $C_{1}$ is a constant depending on $m^*, \rho _*$, $\rho^*$ and $\displaystyle\max _{1\leq i \leq m}\left\|\boldsymbol{w}_{i}\right\|_{L^{\infty}(\Omega)}$.
Finally, the definition of $\boldsymbol{u}_m$ yields an ODE system
\begin{equation} \label{A4}
\begin{cases}
\displaystyle \sum_{k=1}^m{A_{jk}}{\frac{\mathrm{d}}{\mathrm{d}t}g_{m}^{k}(t)}
+\sum_{k,l=1}^m{B_{kl}^{j}}{ h_{m}^{k}(t)g_{m}^{l}(t)}
+\sum_{k=1}^m{C_{jk}} {g_{m}^{k}(t)}
=\int_{\Omega}{\nabla}\chi _m\otimes \nabla \chi _m:\nabla \boldsymbol{w}_j\mathrm{d}x,
\\
g_{m}^{k}(0)=(\boldsymbol{u}_{0m}, \boldsymbol{w}_k), \quad k=1, \cdots, m,
\end{cases}
\end{equation}
where
$$
A_{jk}=\int_{\Omega}{\rho _m}\boldsymbol{w}_j\cdot \boldsymbol{w}_k\mathrm{d}x,
\quad
B_{kl}^{j}=\int_{\Omega}{\rho _m}\left( \boldsymbol{w}_k\cdot \nabla \right) \boldsymbol{w}_l\cdot \boldsymbol{w}_j\mathrm{d}x,
\quad
C_{jk}=\int_{\Omega}{\eta}\left( \chi _m \right) \mathbb{D}\boldsymbol{w}_j: \mathbb{D}\boldsymbol{w}_k\mathrm{d}x.
$$
The matrix $\{A_{i j} \}$ is continuous on $[0, T_0]$ and uniformly positive definite.
Thus standard ODE theory guarantees that there is $0<t_{0} \leq T_{0}$ depending only on $\displaystyle\max _{1 \leq i \leq m}\left|( \boldsymbol{u}_{0m}, \boldsymbol{w}_i)\right|$ such that \eqref{A4}) has a solution $g_{m}(t)=(g_{m}^{1}, \ldots, g_{m}^{m})^{T} \in\left(C\left[0, t_{0}\right]\right)^{m}$ and $g_{m}^{\prime} \in\left(L^{p}\left(0, t_{0}\right)\right)^{m}$ for any $p>1$.

\vskip2mm
\noindent
{\it\bfseries A.2 Fixed point.}\quad
We now show that, for suitable choices of $M$ and $T_m$, the solution $g_m$ remains in the ball $\mathbf{U}_{m}(T_m)$.
Testing (\ref{A4}) with $g_m$ and using (\ref{A2}), (\ref{A3}), we can deduce the following differential inequality
\begin{align} \label{A5}
\frac{\mathrm{d}}{\mathrm{d}t}&\sum_{i=1}^m{\left| g_{m}^{i}(t) \right|^2}
\leq C_{2} \sum_{i=1}^{m} \left|g_{m}^{i}(t)\right|^{2}\left(M^{\frac{1}{2}}+\|\nabla \chi_{m}\|_{L^2(\Omega)}^{2}+1\right)+C_{2}\| \nabla \chi_{m}\|_{L^2(\Omega)}^{2},
\end{align}
where $C_{2}$ depends on $m, \rho^*, \eta^*$ and $\left\|\boldsymbol{w}_{j}\right\|_{C^{1}(\overline{\Omega})}$. Choose
$M=2 e\big(1+\displaystyle\sum_{i=1}^{m}\left|(\boldsymbol{u}_{0m}, \boldsymbol{w}_i)\right|^{2}\big)$ and $T_m\in(0,\min\{t_0,1\})$ satisfying
\begin{align*}
C_2T_m\left( M^{\frac{1}{2}}+1 \right) +C_2T_m\left( \left\| \nabla \chi_0 \right\|^2+C_1T_m \right) e^{C_1MT_m}\le 1,
\end{align*}
we have $g_{m}=(g_{m}^{1}, \ldots, g_{m}^{m})^{T} \in \mathbf{U}_{m}\left(T_{m}\right)$ by Gr\"{o}nwall's inequality.
Consequently, the map $\mathcal{L}: h_m\mapsto g_m$ defined by solving (\ref{A4}) maps $\mathbf{U}_{m}(T_m)$ into itself. Moreover, the ODE system (\ref{A4}) provides a uniform bound on $g_m'$ in $L^p(0,T_m)$ for any $p>1$; by the Arzel\`{a}-Ascoli theorem, $\mathcal{L}$ is compact.
Standard energy estimates imply the continuity of $\mathcal{L}$. Schauder's fixed-point theorem therefore yields a fixed point $g_m=\mathcal{L}(g_m)$, which together with the corresponding $\rho_m, \chi_m$ gives a local solution to (\ref{A.1}) on $[0,T_m]$.

\vskip2mm
\noindent
{\it\bfseries A.3 Global extension.}\quad
To extend the local solution to any $T>0$, we derive a priori bounds independent of the time horizon. Multiplying the momentum equation in (\ref{A1}) by
$g_m^j$ (and summing up), the Allen-Cahn equation by $\mu_m$, integrating over $\Omega$, and adding the resulting identities, we deduce that
\begin{align*}  
\frac{\mathrm{d}}{\mathrm{d}t}\int_{\Omega}{\left( \frac{1}{2}\rho _m\left| \boldsymbol{u}_m \right|^2
+\frac{1}{2}\left| \nabla \chi _m \right|^2 +\rho _mF\left( \chi _m \right) \right)}\mathrm{d}x
+\int_{\Omega}\left( \eta \left( \chi _m \right) \left| \mathbb{D} \boldsymbol{u}_m \right|^2+m(\chi_m)|\mu_m|^2\right)\mathrm{d}x=0.
\end{align*}
Recalling that $\left\|\boldsymbol{u}_m(0) \right\| _{L^2\left( \Omega \right)}^{2}\le \left\| \boldsymbol{u}_0 \right\| _{L^2\left( \Omega \right)}^{2}$ and that $\rho_m$ is bounded from above and below, we infer that
\begin{align*}    
\mathop {\mathrm{sup}}_{0\le t\le T_m}\left( \left\| \boldsymbol{u}_m \right\| _{L^2(\Omega )}^{2}+\left\| \nabla \chi _m \right\| _{L^2(\Omega )}^{2} \right)
+\iint_{Q_{T_m}}{\left( \eta \left( \chi _m \right) \left| \mathbb{D} \boldsymbol{u}_m \right|^2+m(\chi_m)|\mu_m|^2\right)}\mathrm{d}x\mathrm{d}\tau\le C,
\end{align*}
where $C$ depends only on $\rho_0,\chi_{0}, \boldsymbol{u}_{0}$, $\Omega$ and is independent of $m$.
Consequently, no finite-time blow-up can occur, and the local solution of the problem $(\ref{A.1})$ can be extended step by step to the time interval $[0,T]$ for any $T>0$.
\hfill$\Box$

\vskip-4cm
\section{Extension to the Flory-Huggins logarithmic potential}\label{B-F-H}
\noindent
In this appendix we outline the modifications needed to extend the global weak-solution existence result (Theorem \ref{Th4.1}) from the Landau double-well potential to the physically relevant Flory-Huggins logarithmic potential
\begin{align*}
F(s) = \frac{\theta}{2}\big[(1+s)\ln(1+s) + (1-s)\ln(1-s)\big] - \frac{\theta_0}{2}s^2, \qquad s\in[-1,1].
\end{align*}
where $0 < \theta < \theta_0$, and the convex part $F_{\rm conv}(s) = \frac{\theta}{2}[(1+s)\ln(1+s) + (1-s)\ln(1-s)]$. The proof follows the same Galerkin scheme as in Section 3, with the following replacements and uniform estimates that exploit the specific algebraic structure of the Allen-Cahn equation.

\vskip1mm
\noindent
{\it\bfseries B.1  Regularization of the singular potential.}\quad
Choose $\varepsilon_*\in(0,\frac12)$ sufficiently small that
\[
F_{\rm conv}''(1-\varepsilon_*)>\theta_0,
\qquad
F_{\rm conv}'(1-\varepsilon_*)>\theta_0.
\]
For any $\varepsilon\in(0,\varepsilon_*)$, we approximate $F$
by a family of $C^2(\mathbb R)$ potentials,
$$
\Psi_\varepsilon(s) = F_\varepsilon(s) - \frac{\theta_0}{2}s^2, \qquad s\in\mathbb R,
$$
where
\begin{align}\label{F-H}
F_\varepsilon(s) =
\begin{cases}
\displaystyle\sum_{j=0}^{2} \frac{1}{j!}F_{\rm conv}^{(j)}(1-\varepsilon)[s-(1-\varepsilon)]^j, & s \geq 1-\varepsilon, \\[6pt]
F_{\rm conv}(s), & s \in [-1+\varepsilon, 1-\varepsilon], \\[4pt]
\displaystyle\sum_{j=0}^{2} \frac{1}{j!}F_{\rm conv}^{(j)}(-1+\varepsilon)[s-(-1+\varepsilon)]^j, & s \leq -1+\varepsilon.
\end{cases}
\end{align}
It is easy to see that $\Psi^{''}_\varepsilon(s)\ge-\theta_0$ uniformly in $\varepsilon$, $\Psi'_\varepsilon(s)\rightarrow F'(s)$ on $(-1,1)$, and
$\Psi_\varepsilon(s)\ge-\frac{\theta_0}{2}$ for $s\in\mathbb R$. Replacing $F$ by $\Psi_\varepsilon$ in (\ref{1.1})$_5$, the existence of global weak solutions $(\rho_\varepsilon, u_\varepsilon, \chi_\varepsilon, \mu_\varepsilon)$ for each fixed $\varepsilon>0$ follows from the same Galerkin argument as in the proof of Theorem \ref{Th4.1}; the bounds obtained below are uniform in $\varepsilon$.
By the choice of $\varepsilon_*$ and symmetry,
$\Psi_\varepsilon'(1)>0$ and
$\Psi_\varepsilon'(-1)<0$. Hence the comparison principle
gives $|\chi_\varepsilon|\le1$.

\vskip1mm
\noindent
{\it\bfseries B.2  Uniform estimates independent of $\varepsilon$.}\quad
The energy balance (\ref{4.7}) gives, uniformly in $\varepsilon$
\begin{align}\label{B1}
\sup_{0\le t\le T}\left(\|u_\varepsilon\|_{L^2(\Omega)}^2+\|\nabla\chi_\varepsilon\|_{L^2(\Omega)}^2\right)
+\int_0^T\int_{\Omega}\left(|{\mathbb D}u_\varepsilon|^2+|\mu_\varepsilon|^2\right){\rm d}x{\rm d}t\le C,
\end{align}
where $C$ is independent of $\varepsilon$. Hence $\|\mu_\varepsilon\|_{L^2(0,T; L^2(\Omega))}\le C$.

From the relation
\begin{align}\label{B3}
\rho_\varepsilon\mu_\varepsilon = -\Delta\chi_\varepsilon + \rho_\varepsilon \Psi_\varepsilon'(\chi_\varepsilon),
\end{align}
we obtain the identity
\begin{align*}
\int_\Omega \rho_\varepsilon\mu_\varepsilon^2{\rm d}x
+ 2\theta_0\int_\Omega |\nabla\chi_\varepsilon|^2{\rm d}x
= \int_\Omega \frac{(\Delta\chi_\varepsilon)^2}{\rho_\varepsilon}{\rm d}x
+ 2\int_\Omega F_\varepsilon''(\chi_\varepsilon)|\nabla\chi_\varepsilon|^2{\rm d}x
+ \int_\Omega \rho_\varepsilon(\Psi_\varepsilon'(\chi_\varepsilon))^2{\rm d}x,
\end{align*}
where $F_\varepsilon''(s) \geq 0$ for all $s \in \mathbb{R}$ by construction. Since the left-hand side is bounded in $L^1(0,T)$ uniformly in $\varepsilon$ by (\ref{B1}), we obtain
\begin{align}\label{B2}
\int_0^T\!\!\int_\Omega \frac{(\Delta\chi_\varepsilon)^2}{\rho_\varepsilon}\,dxdt \leq C, \qquad \int_0^T\!\!\int_\Omega \rho_\varepsilon(\Psi_\varepsilon'(\chi_\varepsilon))^2\,dxdt \leq C,
\end{align}
with $C$ independent of $\varepsilon$. Thus $\|\Delta\chi_\varepsilon\|_{L^2(0,T; L^2(\Omega))}$ and $\|\Psi_\varepsilon'(\chi_\varepsilon)\|_{L^2(0,T; L^2(\Omega))}$ are bounded uniformly.

\vskip1mm
\noindent
{\it\bfseries B.3  Passage to the limit and confinement $|\chi|<1$.}\quad
From (\ref{B1}), (\ref{B2}), and the uniform $L^\infty$ bound on $\rho_\varepsilon$, there exists a subsequence (not relabelled) such that
$$
\chi_\varepsilon\rightharpoonup\chi \quad \text{weakly~in}~L^2(0,T; H^2(\Omega)),
\qquad
\Psi'_\varepsilon(\chi_\varepsilon)\rightharpoonup\Psi^* \quad \text{weakly~in}~L^2(0,T; L^2(\Omega)).
$$
To identify $\Psi^*=F'(\chi)$ and to prove $|\chi|<1$ a.e., we use the standard measure argument. For any $\delta>0$, define
$$
A^+_{\varepsilon, \delta}=\big\{(x,t)\in\Omega\times(0,T);~ \chi_\varepsilon(x,t)>1-\delta\big\}.
$$
{For sufficiently small $\delta>0$ and $0<\varepsilon<\delta$,
the monotonicity of $F_\varepsilon'$ and
$|\chi_\varepsilon|\le1$ imply that, on
$A^+_{\varepsilon,\delta}$,
\[
\Psi_\varepsilon'(\chi_\varepsilon)
=F_\varepsilon'(\chi_\varepsilon)-\theta_0\chi_\varepsilon
\ge F_{\rm conv}'(1-\delta)-\theta_0>0.
\]
Consequently,
\[
\bigl(F_{\rm conv}'(1-\delta)-\theta_0\bigr)^2
\,{\rm meas}(A^+_{\varepsilon,\delta})
\le
\int_{A^+_{\varepsilon,\delta}}
|\Psi_\varepsilon'(\chi_\varepsilon)|^2
\,{\rm d}x\,{\rm d}t
\le C,
\]
where $C$ is independent of $\varepsilon$ and $\delta$.
Since $F_{\rm conv}'(1-\delta)\to+\infty$ as
$\delta\to0$, we obtain}
$$
\lim_{\delta\rightarrow0}\limsup_{\varepsilon\rightarrow0}{\rm meas}(A^+_{\varepsilon, \delta})=0.
$$
The same argument applies to $\big\{(x,t)\in\Omega\times(0,T);~ \chi_\varepsilon(x,t)<-1+\delta\big\}$. Hence, up to a subsequence, $\chi_\varepsilon\rightarrow\chi$ a.e. and $|\chi(x,t)|<1$ a.e. in $\Omega\times(0,T)$. Moreover, by Egorov's theorem and the uniform bound (\ref{B2}), the convergence $\Psi'_\varepsilon(\chi_\varepsilon)\rightharpoonup F'(\chi)$ holds weakly in $L^2(0,T;L^2(\Omega))$.
Passing to the limit in (\ref{B3}) gives
$$
\rho\mu=-\Delta\chi+\rho F'(\chi).
$$

\vskip1mm
The rest of the proof (passage to the limit in the Galerkin scheme, time continuity, and weak formulation) follows as in Section 3, noting that the singular nature of $F'$ near $\pm 1$ is handled by the $L^2(L^2)$ bound and the a.e. constraint $|\chi| < 1$.
Thus Theorem \ref{Th4.1} also holds for the Flory-Huggins potential \eqref{1.5}, with the additional property $|\chi(x,t)|<1$ a.e. in $\Omega\times(0,T)$.
\hfill$\Box$

\begin{Remark}
The key simplification in the NSAC system compared with the NSCH system (cf. \cite{G-T}) is that the Allen-Cahn equation provides an $L^2(L^2)$ bound on $\mu$ itself, and the algebraic relation (\ref{B3}) directly yields the $L^2(L^2)$ bound on $\Psi'_\varepsilon(\chi_\varepsilon)$ without requiring the intricate monotonicity estimates needed in the Cahn-Hilliard case to control $\|F_\varepsilon'(\chi_\varepsilon)\|_{L^1(\Omega)}$.
\end{Remark}



\end{document}